\documentclass[12pt]{article}

\usepackage[small,compact]{titlesec}
\usepackage{setspace}
\usepackage{wrapfig,caption,epsfig,subcaption,sidecap}
\usepackage{verbatim,algorithmic}
\usepackage[ruled,boxed]{algorithm}
\usepackage[T1]{fontenc}

\usepackage{authblk}

\usepackage[margin=1in]{geometry}

\usepackage[charter,uppercase=italic,lowercase=italic,greekuppercase=upright,greeklowercase=italic]{mathdesign}
\usepackage{inconsolata}
\renewcommand{\leq}{\leqslant}
\renewcommand{\geq}{\geqslant}

\newcommand{\keqref}[1]{Eq.~(\ref{#1})}
\newcommand{\meqref}[1]{Eq.~(\ref{#1})}
\newcommand{\applabel}[1]{\label{app-#1}}
\newcommand{\appref}[1]{\ref{app-#1}}

\usepackage{amsmath,amsfonts,amsthm}
\usepackage{tikz-cd}
\usepackage{mathtools}

\numberwithin{equation}{section}

\renewenvironment{itemize}{
  \begin{list}{-}
    {\setlength{\parsep}{3pt}
      \setlength{\labelwidth}{24pt}
      \setlength{\itemsep}{1pt}
      \setlength{\topsep}{3pt}}}{\end{list}}

\usepackage{graphicx,color}
\usepackage{url,hyperref}
\usepackage{soul}

\newcommand{\C}{\mathbb C}
\newcommand{\E}{\mathbb E}
\newcommand{\cov}{\mathrm{cov}}
\newcommand{\Ecal}{\mathcal E}

\newcommand{\R}{\mathbb R}
\newcommand{\Z}{\mathbb Z}

\def\wv{u}
\def\wu{ \widetilde{u}}  

\newcommand{\spaceX}{\mathbb X}
\newcommand{\spaceH}{\mathbb H}
\newcommand{\spaceV}{\mathbb V}

\newcommand{\ad}[1]{{#1}^\dagger}
\newcommand{\braket}[1]{\langle{#1}\rangle}
\newcommand{\dt}{{\mbox{\scriptsize$\Delta$}}{t}}
\newcommand{\dtheta}{{\mbox{\scriptsize$\Delta$}}{\theta}}

\newcommand{\floor}[1]{\lfloor{#1}\rfloor}

\newcommand{\terms}{~+~\cdots~+~}

\newcommand{\ub}[2]{\underbrace{{#1}}_{\mbox{\footnotesize{#2}}}}

\renewcommand{\bar}[1]{{\overline{#1}}}
\renewcommand{\hat}[1]{{\widehat{#1}}}
\renewcommand{\phi}{\varphi}
\renewcommand{\tilde}[1]{{\widetilde{#1}}}

\definecolor{kcolor}{rgb}{0.0,0.6,0.0}

\newcommand{\heading}[1]{\smallskip\noindent{\em{#1}}}
\newcommand{\ubrace}[2]{\underbrace{~{#1}~}_{\mbox{{\footnotesize{#2}}}}}

\newcommand{\fcommentout}[1]{} 

\begin{document}

\title{Data-driven model reduction, Wiener projections, and the
  Koopman-Mori-Zwanzig~formalism}

\author[1]{Kevin K.~Lin}

\affil[1]{\footnotesize Department of Mathematics, University of
  Arizona, Tucson, AZ 85721, USA
  (\href{mailto:klin@math.arizona.edu}{klin@math.arizona.edu})}

\author[2]{Fei Lu}

\affil[2]{\footnotesize Department of Mathematics, Johns Hopkins
  University, Baltimore, MD 21218, USA
  (\href{mailto:feilu@math.jhu.edu}{feilu@math.jhu.edu})}

\date{\today}

\maketitle

\begin{abstract}
  Model reduction methods aim to describe complex dynamic phenomena
  using only relevant dynamical variables, decreasing computational
  cost, and potentially highlighting key dynamical mechanisms. In the
  absence of special dynamical features such as scale separation or
  symmetries, the time evolution of these variables typically exhibits
  \emph{memory effects}. Recent work has found a variety of
  data-driven model reduction methods to be effective for representing
  such non-Markovian dynamics, but their scope and dynamical
  underpinning remain incompletely understood. Here, we study
  data-driven model reduction from a dynamical systems
  perspective. For both chaotic and randomly-forced systems, we show
  the problem can be naturally formulated within the framework of
  Koopman operators and the Mori-Zwanzig projection operator
  formalism. We give a heuristic derivation of a NARMAX (Nonlinear
  Auto-Regressive Moving Average with eXogenous input) model from an
  underlying dynamical model. The derivation is based on a simple
  construction we call \emph{Wiener projection}, which links
  Mori-Zwanzig theory to both NARMAX and to classical Wiener
  filtering.  We apply these ideas to the Kuramoto-Sivashinsky model
  of spatiotemporal chaos and a viscous Burgers equation with
  stochastic forcing.
\end{abstract}

\begin{footnotesize}
  \tableofcontents
\end{footnotesize}

\section{Introduction}

Unsteady fluid flow, fluctuations in power grids, neural activity in the
brain: these and many other complex dynamical phenomena arise from the
interaction of a large number of degrees of freedom across many orders
of magnitude in space and time.  But, in these and many other systems,
only a relatively small subset of the dynamical variables are of direct
interest or even observable.  Reduced models, i.e., models that use only
relevant dynamical variables to reproduce dynamical features of interest
on relevant timescales, are thus of great potential utility, especially
in tasks requiring repeated model runs like uncertainty quantification,
optimization, and control.  Moreover, relevant dynamical mechanisms are
often easier to glean and understand in reduced models.

Many approaches to model reduction --- also known as the closure problem
in physics and reduced-order modeling in engineering --- have been
proposed.  On one hand, a variety of analytical and computational
methods have been proposed based on dynamical systems theory and
statistical mechanics.  These have been especially successful in
situations with special dynamical features like sharp scale separation,
low dimensional attractors, or symmetries; see,
e.g.,~\cite{pavliotis2008multiscale, kevrekidis, roberts2014model,
  abdulle2012heterogeneous}.  However, not all scientific and
engineering applications exhibit these features, and in such cases
reduced models must account for memory and noise effects (see, e.g.,
\cite{CH13, zwanzig} as well as Sect.~\ref{sect:bkgnd}).
On the other hand, while purely data-driven approaches, i.e., those
based on fitting generic statistical models to simulation data or
physical measurements, have been quite successful in a variety of
settings without sharp scale separation (see,
e.g.,~\cite{kondrashov_Datadriven2015a, CL15,
  harlimParametricReduced2015, leiDatadrivenParameterization2016,
  xieDataDrivenFiltered2018, chekroun_Dataadaptive2017a}, the dynamical
basis for these methods is often unclear, and as a result their scope of
applicability remain incompletely understood.  In addition, a systematic
understanding from the nonlinear dynamical systems point of view would
provide a framework for analyzing and improving these methods.

This paper is the first step in our effort to bridge this gap; for
different perspectives and approaches to similar questions, see, e.g.,
\cite{berry2020bridging,kutz2016dynamic,mezic2005spectral,williams2015data}.  First,
using Koopman operators, the Mori-Zwanzig formalism, and Wiener
filtering, we propose a simple mathematical formulation of data-driven
model reduction.  The resulting framework links dynamical systems theory
and data-driven modeling, and can serve as a starting point for
systematic approximations in model reduction.  In particular, we show
that a variant of the NARMAX (Nonlinear Auto-Regressive Moving Average
with eXogenous input) representation of stochastic processes, widely
used in time series analysis and data-driven modeling (see
\cite{Ham94,Bil13,CL15} and references therein), can be derived via a
construction we call ``Wiener projections,'' which is equally applicable
to either deterministic chaotic or random dynamical systems.  Another
consequence of our work is that for problems with time-stationary
statistics, classical Wiener filtering can provides an alternative to
Mori-Zwanzig as a framework for model reduction.

\heading{Organization.}  In Sect.~\ref{sect:bkgnd}, we recall relevant
dynamical systems theory background, including a discrete-time version of the
Mori-Zwanzig formalism and the NARMAX representation of stochastic
processes.  We also formulate the problem of data-driven driven model
reduction considered in this paper.  Sect.~\ref{sect:wiener-projections}
describes the Wiener projection and its basic properties, and shows how
it can be used to derive a variant of NARMAX.  Sect.~\ref{sect:numerics}
is concerned with numerical implementation details, and
Sect.~\ref{sect:numAll} examines the application of these ideas to the
Kuramoto-Sivashinsky partial differential equation (PDE) and to a
stochastic Burgers equation.  For the convenience of readers, we have
included appendices on an alternate derivation of the Mori-Zwanzig
equation (which sheds some light on its interpretation); a summary of
classical Wiener filtering and the $z$-transform; and detailed numerical
results on our two examples.

\section{Data-driven model reduction in discrete time}
\label{sect:bkgnd}

\subsection{Problem formulation and dynamical systems setting}
We assume the full system of interest is a discrete-time dynamical
system
\begin{equation}
  \label{eq:system}
  X_{n+1} = F(X_n).
\end{equation}
The states $X_n$ are points in a space $\spaceX$, which can be a vector
space, a manifold, or a more general space.  We refer to
\keqref{eq:system} as the \emph{full model}.  The dynamical variables of
interest, or {\em relevant variables}, are defined by $x=\pi(X)$, $\pi$
being a given function mapping points in $\spaceX$ to points in
$d$-dimensional Euclidean space $\R^d$, generally with
$d\ll\dim(\spaceX)$.  (The choice of $\pi$ is dictated in part by the
application, in part by dynamical considerations such as scale
separation.)  \keqref{eq:system} can accommodate continuous time systems
by letting $F$ be the time-$\dt$ solution map (for some $\dt>0$) or a
Poincar\'e map.  We focus on discrete-time reduction because (i)
observations are always discrete in time, and (ii) discrete-time reduced
models avoid the numerical errors that come from integrating continuous
time reduced models, which can be significant in chaotic
systems~\cite{CL15,LLC16}.

By data-driven model reduction, we mean using data to construct reduced
models that use only the relevant variables.  We are interested in
reduced models that can (i) forecast $x_n$ given its past history, and
(ii) reproduce long-time statistics, e.g., correlations and marginal
distributions.  In general, parametric model reduction methods begin
with a family of models with unknown parameters
and observations
$\tilde{x}_n=\pi(\tilde{X}_n)$, where $(\tilde{X}_n)_{n=0}^N$ is a
trajectory (or multiple trajectories) of the full model.  One then estimates the parameters 
by fitting the model
to the data, usually by minimizing a suitable loss function.  Methods
differ in their choice of models and loss functions, which can impact
both model fitting and the performance of the reduced model.

\medskip

A useful approach to the statistical properties of dynamical systems is
to view the space of observables on $\spaceX$ as forming a Hilbert space
$\spaceH = L^2(\mu)$ with inner product $\braket{f,g}=\int fg~d\mu$.
The probability distribution $\mu$ describes the long-time statistics of
typical solutions of \keqref{eq:system}, and is \emph{invariant}, i.e.,
if $X_0$ has distribution $\mu$, then so does $X_n$ for all $n>0$.
Inner products are thus naturally interpreted as steady-state
correlations.  Many dynamical systems of interest possess multiple
``natural'' invariant measures; the choice of a suitable measure is
dictated in part by the application, in part by computational
tractability.  For example, in molecular dynamics, one may consider
microcanonical, canonical, or grand canonical ensembles; for dissipative
chaotic systems, relevant invariant probability measures are often
singular distributions supported on strange attractors that nevertheless
reflect the statistics of a set of initial conditions with positive
phase space volume.

In principle, the measure $\mu$ need not be invariant.  But invariance
significantly simplifies the problem of data-driven model reduction, and
in addition guarantees many convenient mathematical properties.  Without
the invariance assumption, far more data would be needed.  For these
reasons, we focus on stationary processes in this paper.  Equivalently,
we assume $X_0$ has distribution $\mu$, so that $(X_n)$ is
stationary.

\medskip

Recall that the \emph{Koopman operator} is the operator $M$ defined by
$M\phi(X) = \phi(F(X))$.  The Koopman operator advances observables
forward in time: $M\phi(X)$ gives the value of $\phi$ at the next step
if the current state is $X$.  The Koopman operator and its adjoint, the
Perron-Frobenius transfer operator, describes the dynamics from the
function space point of view.  Much is known about their properties as
operators on $\spaceH$ and on other relevant Banach spaces, see, e.g.,
\cite{walters2000introduction} or \cite{reed1980methods}.  Both the
Koopman and Perron-Frobenius operators have been used extensively in
computational nonlinear dynamics; see, e.g.,
\cite{froyland2009almost,klus2018data}.

We will use extensively two properties:
\begin{enumerate}

\item\label{koopman1} The Koopman operator is invertible when $F$ is
  invertible, and $M^{-1}\phi = \phi\circ F^{-1}$.

\item\label{koopman2} With the inner product $\braket{\cdot,\cdot}$
  above, Koopman operators are Hilbert space isometries (i.e.,
  $\braket{Mu,Mv} = \braket{u,v}$) and unitary ($MM^*=M^*M=I$) when $F$
  is invertible~\cite{reed1980methods,walters2000introduction}.
  
\end{enumerate}
We note that property~\ref{koopman2} relies on the invariance of $\mu$.

One of the uses of Koopman operators (and the Mori-Zwanzig formalism
introduced in the next section) is to turn nonlinear dynamics questions
into questions involving linear operators, for which mathematical
analysis and formal manipulation are often easier.  We will take
advantage of this in Sect.~\ref{sect:wiener-projections}.

\subsection{Discrete-time Mori-Zwanzig formalism}
\label{sect:mz-review}

The MZ formalism originally arose in classical statistical
mechanics~\cite{zwanzig,CH13}, and has been used in physical
applications ranging from fluid dynamics to materials science and
molecular dynamics (see, e.g.,
\cite{CH13,chorin2002optimal,leiDatadrivenParameterization2016,ma2019_CoarsegrainingLangevin,cho2014statistical,li2017computing,li2015incorporation,panchenko,venkataramani,stinis2004stochastic,parish2017non,wang2019implicit}).
As we will discuss in Sect.~\ref{sect:wiener-rds}, it also applies to systems
with random forcing and/or (bounded) delays.  Here we review a discrete
MZ theory \cite{DSK09}.

The starting point of Mori-Zwanzig formalism is the \emph{Mori-Zwanzig
  equation}, which asserts that there exists a sequence of functions
$\xi_1,\xi_2,\cdots:\spaceX\to\R^d$ such that for $n\geq0$,
\begin{subequations}
  \label{eq:mz} 
  \begin{equation}
    x_{n+1} = PF(x_{n}) + \sum_{k=1}^n \Gamma_k(x_{n-k}) + \xi_{n+1}(X_0)
  \end{equation}
  with
  \begin{equation}
    \label{eq:mz-aux} 
    \Gamma_k = P(\xi_k\circ F)\qquad\mbox{and}\qquad P\xi_n=0.
  \end{equation}
\end{subequations}
In Eq.~(\ref{eq:mz}), $P$ can be any projection operator.  The first,
``Markov'' term is then the ``best'' approximation of $F$ by functions
in the range $\spaceV$ of the projection $P$ (more on this below).  The
second, ``memory'' term captures all non-Markovian effects representable
in $\spaceV$.  The last, ``noise'' term represents errors at each step,
and are orthogonal to functions in $\spaceV$.  We present a derivation
below, and an alternate derivation via a dual equation in
\ref{app:dual}.

To make use of the Mori-Zwanzig equation, one must first choose a
projection operator $P$.  A common choice is the conditional expectation
$(P\varphi)(x) = \E_\mu[\varphi(X)|\pi(X)=x]$.  Another is {\it finite
  rank projection:} fix a collection of linearly independent functions
$\psi_1(x),\cdots,\psi_\nu(x)$ of $x$, then take $P$ to be orthogonal
projection onto their linear span, i.e.,
\begin{equation}
  \label{eq:projection}
  P\phi(x) = \Psi(x)\cdot\braket{\Psi,\Psi}^{-1}\cdot\braket{\Psi,\phi},
\end{equation}
where $\braket{f,g}=\int f^T\cdot g~d\mu$ for matrix-valued $f$ and $g$,
and the columns of $\Psi(x) = [\psi_1(x)~~\cdots~~\psi_\nu(x)]$ span
$\spaceV$.  With $P$ as in \keqref{eq:projection}, we can write $PF =
\Psi\cdot c_0$ and $\Gamma_k = \Psi\cdot c_k$ for coefficient vectors
$c_k.$ \keqref{eq:mz} then becomes
\begin{subequations}
  \label{eq:mz-filter}
  \begin{equation}
    \label{eq:mz-filter-main}
    x_{n+1} = \sum_{k\geq0}\Psi(x_{n-k})\cdot c_k +
    \xi_{n+1}.
  \end{equation}
  Eq.~(\ref{eq:mz-aux}) now take the form
  \begin{equation}
    \label{eq:mz-filter-fdt}
    c_k = \braket{\Psi,\Psi}^{-1}\cdot\braket{\Psi,\xi_k\circ F}
  \end{equation}
  and
  \begin{equation}
    \label{eq:mz-filter-ortho}
    \braket{\xi_n,\Psi\circ\pi}=0.
  \end{equation}
\end{subequations}
In this paper, we will mainly consider finite rank projections and a closely related ``Wiener projection'' 
 in Sect.~\ref{sect:wiener-projections}.  See, e.g., \cite{CH13,
  grabert, zwanzig, ma2019_CoarsegrainingLangevin} for discussions of
the conditional expectation and other choices.

The Mori-Zwanzig equation is an exact description of the dynamics of
$x_n$.  Without further approximation, it does not represent a reduction
in model complexity.  The equation does, however, highlight the
interdependence of the projection $P$ and the noise $(\xi_n)$.  To
arrive at closed equations of motion for the relevant variables $x_n$,
it is necessary to choose $P$ so that the noise terms $(\xi_n)$ can be
effectively modeled.  A common approach is to choose $P$ to be a
projection onto the slow variables.  One then appeals to scale
separation and other physical considerations to justify modeling
$(\xi_n)$ by a stochastic process $\eta_n$, e.g., a stationary Gaussian
process.  The coefficients $(c_k)$ can be approximated by, e.g.,
perturbation techniques.  The power spectrum of the noise and the memory
kernel are related by so-called fluctuation-dissipation relations, of
which Eq.~(\ref{eq:mz-filter-fdt}) is an example~\cite{forster,zwanzig}.

As a physical example, one may consider the motion of heavy particle
suspended in a fluid, a problem originally studied by Smoluchowski and
Einstein~\cite{einstein1956investigations}.  The ``system'' consists of
the heavy particle and the water molecules making up the surrounding
fluid.  Projecting onto the particle degrees of freedom, the Markov term
is given by equation of motion for a free particle, the memory term
gives rise to drag due to the fluid, and the noise term represents
random forces due to thermal fluctuation of the surrounding fluid.

Orthogonality conditions (e.g., \keqref{eq:mz-filter-ortho}) play a key
role in MZ theory and in Wiener filtering: they are equivalent to
optimality in the least squares sense.  In using reduced models to
generate predictions, one often assumes the driving noise (i.e., the
$\xi_n$ in Eq.~(\ref{eq:mz-filter})) is independent of $x_m$ for $n >
m$.  Orthogonality conditions provide partial justification for this
(standard) procedure.  \keqref{eq:mz-filter-ortho} comes from
$P\xi_n=0$, but does \emph{not} imply $\Psi(x_m)$ is uncorrelated with
$\xi_n$ for $n > m$.  More on this in
Sect.~\ref{sect:wiener-projections}.

\heading{Derivation of the Mori-Zwanzig equation.}  The MZ equation can
be driven as follows. We start with the \emph{Dyson formula}
\begin{equation}
  \label{eq:dyson}
  M^{n+1} = \sum_{k=0}^n M^{n-k}PM(QM)^k + (QM)^{n+1} ,
\end{equation}
where, as before, $M$ is the Koopman operator and $P$ is a projection on
$\spaceH$ whose range $\spaceV$ are functions that depend \emph{only} on
the relevant variables $x$; and $Q=I-P$ is the orthogonal projection.
Eq.~(\ref{eq:dyson}) is readily proved by induction.  To see how
Eq.~(\ref{eq:mz-filter}) follows from Eq.~(\ref{eq:dyson}), apply both
sides of Eq.~(\ref{eq:dyson}) to the observation function $\pi$ and
evaluate at $X_0$, yielding
\begin{equation}
  \ubrace{(M^{n+1}\pi)(X_0)}{(I)} = \ubrace{\sum_{k=0}^n
    (M^{n-k}PM(QM)^k\pi)(X_0)}{(II)} +
  \ubrace{((QM)^{n+1}\pi)(X_0)}{(III)}.
\end{equation}
Define $\xi_n = (QM)^n\pi$, so that $P\xi_n=0$ for $n\geq1$.  For Term
(I), the definition of the Koopman operator $M$ gives $\pi(F^{n+1}(X_0))
= \pi(X_{n+1}) = x_{n+1}.$ For Term (III), we have (by definition)
$\xi_{n+1}(X_0)$.  For Term (II), we have
\begin{displaymath}
  (M^{n-k}PM(QM)^k\pi)(X_0) = (PM(QM)^k\pi)(X_{n-k})
\end{displaymath}
as before.  Since $M(QM)^k\pi = M\xi_k = \xi_k\circ F$ and the range of
$P$ consists of functions of $x=\pi(X)$, we get
\begin{displaymath}
  (M^{n-k}PM(QM)^k\pi)(X_0) = P(\xi_k\circ F)(x_{n-k}).
\end{displaymath}
Combining all these and $PQ=0$ yields Eq.~(\ref{eq:mz}).

\subsection{NARMAX modeling}
\label{sect:narmax}

Whereas MZ theory seeks systematic derivations of reduced models, NARMAX
(Nonlinear Auto-Regressive Moving Average with eXogenous input) is a
generic approach to parametric data-driven modeling of time
series~\cite{Bil13,FY03,Ham94}. 
A common version of the NARMAX model is
\begin{subequations}
  \label{eq:narmax}
  \begin{align}
    \label{eq:narmax-predictor}
    x_{n+1} =& ~f(x_n) + z_n , \\
    \label{eq:narmax-corrector}
    z_n + a_{p-1}z_{n-1} \terms a_0z_{n-p}
    =& ~d_q w_n \terms d_0w_{n-q}\\
    &+\Psi(x_n)\cdot b_1 \terms \Psi(x_{n-r})\cdot b_r ,\nonumber
  \end{align}
\end{subequations}
where $f$ and $\Psi$ are given functions, and the $w_i$ are independent
identically distributed (IID) random variables, usually assumed to be
Gaussian (as we do here).  One can view $x_{n+1}=f(x_n)$ as a crude
predictor of $x_{n+1}$, and \keqref{eq:narmax-corrector} a corrector
based on a model of the residuals $z_n$.  Note that like the MZ
equation, \keqref{eq:narmax} is non-Markovian.

In applications of NARMAX, the main task of the would-be modeler is to first
choose the forms of $f$, $\Psi$ and the orders $p, q, r$, then
determine $a_i$, $b_i$, and $d_i$ by minimizing a suitable loss
function.  One common approach to parameter estimation is least squares
regression: let $\tilde{x}_n$ denote time series obtained from the full
model (either by simulation or physical measurement), and define
\begin{subequations}
  \begin{align}
    \hat{x}_{n+1} =& ~f(\tilde{x}_n) + \tilde{z}_n , \\
    \tilde{z}_n + a_{p-1}\tilde{z}_{n-1} \terms a_0\tilde{z}_{n-p}
    =& \Psi(\tilde{x}_n)\cdot b_1 \terms \Psi(\tilde{x}_{n-r})\cdot b_r
  \end{align}
\end{subequations}
The $\hat{x}_{n+1}$ is the one-step prediction based on $\tilde{x}_n,
\cdots, \tilde{x}_{n-r}$.  One then tunes $(a_i,b_i)$ to minimize the mean
squared error $\sum_n\|x_n-\hat{x}_n\|^2$, possibly in combination with
regularization techniques, e.g., Tikhonov regularization or a
sparsity-inducing $\ell^1$ term.  The moving average coefficients $d_n$
are determined by fitting a stochastic process of the form $d_qw_n
\terms d_0w_{n-q}$ to the residual.  Another approach to parameter estimation is based on maximum likelihood
estimation (MLE).  In this approach, one assumes the statistics of the
noise $(w_n)$, e.g., independent $N(0,I)$ random vectors, and infer the
$(a_i,b_i)$ and $d_i$ jointly by maximum likelihood methods and variations
thereof.

Whatever the method, we emphasize that the form of
Eqs.~(\ref{eq:narmax}) does not, by itself, determine a reduced model or
a model reduction procedure.  One must either specify the statistics of
the noise term, or the loss function to be minimized.  (And, for
non-convex loss functions, the optimization procedure.)  These choices
can have a significant impact on the usefulness of the model so
obtained.

\section{Wiener projections}
\label{sect:wiener-projections}

\subsection{Definition and basic properties}
\label{sect:wiener-stuff}

We now set aside Mori-Zwanzig for a moment, and consider another way to
conceptualize memory effects in model reduction based on Wiener
filters~\cite{hannan,kailath}.  Let $u_n$ and $v_n$ be two zero-mean
wide-sense stationary processes.  The Wiener filter is the sequence
$(h_n)$ that minimizes the mean-squared error (MSE):
\begin{equation}
  \label{eq:mse}
  \E\big(\big\|u_n-(v\star h)_n\big\|^2\big),  
\end{equation}
where $(v\star h)_n =\sum_{k}v_{n-k}\cdot h_k$ denotes convolution, with
$h_n=0$ for $n<0$.  (See \ref{app:wiener-and-z} for more
details.)  It satisfies the orthogonality condition
\begin{equation}
  \label{eq:wiener-ortho}
  \cov(v_m,\overline{r}_n)=0~,~~~n\geq m,
\end{equation}
where $r_n$ is the residual $u_n-\sum_{k}v_{n-k}\cdot h_k$, i.e., filter
errors are uncorrelated with the data on which the filter output is
based.  \keqref{eq:wiener-ortho} is equivalent to the minimum-MSE
criterion.

We observe that the Wiener filter can be applied to model reduction as
well: with $X_n$ as in \keqref{eq:system} and $\Psi$ as before, let
$h_n$ be the causal Wiener filter for $u_n=x_{n+1}=\pi(X_{n+1})$ and
$v_n=\Psi(x_n)$.  We then obtain $x_{n+1} =
\sum_{k\geq0}\Psi(x_{n-k})\cdot h_k + r_{n+1}$ with $\cov(\Psi(
x_m),r_n)=0$ for $n>m$ with $r_n$ playing the role of the residual $r_n$
in Eq.~(\ref{eq:wiener-ortho}).

How is this Wiener filter view related to the MZ formalism?  We now
sketch an argument showing that Wiener-based model reduction is in fact
a special case of the MZ equation, one with some attractive properties.
Let $\Psi_n=\Psi(x_n)$, and assume $F$ is invertible
so that $M$ is invertible and unitary.  Let $P_W$ be orthogonal
projection onto the subspace
\begin{equation} \label{spaceW}
  W~=~{\rm span}(\Psi \cup M^{-1}\Psi \cup M^{-2}\Psi \cup
    \cdots),
\end{equation}
where $M^{-k}\Psi$ is a short-hand for $\{M^{-k}\psi_1, \cdots,
M^{-k}\psi_\nu\}.$ Note $P_W=P_W^*$, i.e., $P_W$ is self-adjoint.  Since
$M^{-1}v\in W$ for all $v\in W$, we have
\begin{equation}
  \label{eq:tower}
  M^{-\ell}P_W = P_WM^{-\ell}P_W~,~~\ell\geq0.
\end{equation}
This implies $Q_WM^{-\ell}P_W=0$ or, upon taking adjoints,
$P_WM^{\ell}Q_W=0$.  The Dyson formula~(\ref{eq:dyson})
for $P_W$ thus simplifies:
\begin{subequations}
  \begin{align}
    \label{eq:wiener-dyson}
    M^{n+1} &= \sum_{k=0}^n M^{n-k}P_WM(Q_WM)^k + (Q_WM)^{n+1}\\
    \label{eq:pathProj}
    &= M^nP_WM + (Q_WM)^{n+1}
  \end{align}
\end{subequations}
since $PM(QM)^k=0$ for $k\geq1$.  Applying both sides of
\keqref{eq:pathProj} to $\pi$, we obtain
\begin{subequations}
  \label{eq:wiener-projection}
  \begin{align}
    \label{eq:wiener-projection-main}
    x_{n+1} &= \sum_{k\geq0}\Psi(x_{n-k})\cdot h_k + \xi_{n+1},\\
    \label{eq:wiener-projection-ortho}
    &\braket{\xi_n,\Psi(x_m)}=0~,~~n > m.
  \end{align}
\end{subequations}
Though \keqref{eq:wiener-projection-main} and \keqref{eq:mz-filter-main}
are formally identical, the orthogonality
relation~(\ref{eq:wiener-projection-ortho}) is strictly stronger than
\keqref{eq:mz-filter-ortho}.  The reason is that for the finite rank
projection in Sect.~\ref{sect:mz-review}, the orthogonality relation
means $\int\xi_n(X)^T\cdot\Psi(\pi X)~d\mu(X)=0$, i.e., the noise
functions $\xi_n$ are orthogonal to a finite dimensional subspace of
$\spaceH$.  In contrast, in Eq.~(\ref{eq:wiener-projection-ortho}), the
orthogonality relation $\mathbb{E}_{X_0\sim \mu}[\xi_n(X_0)^T\cdot\Psi(\pi
  X_m)]=0$ means the $\xi_n$ is (in general) orthogonal to an infinite
dimensional subspace of $\spaceH$, and is analogous to
Eq.~(\ref{eq:wiener-ortho}), where the expectation is with respect to
the stationary measure $\mu$ on a suitably defined path space.  The
orthogonality~(\ref{eq:wiener-projection-ortho}) is significant for two
reasons.  First, in stochastic models like NARMAX~(\ref{eq:narmax}), one
typically assumes the driving noise $w_m$ is independent of $x_n$ for
$m>n$.  While natural, this is not guaranteed by the MZ equation.
Eq.~(\ref{eq:wiener-projection-ortho}) does not imply such independence,
either, but comes a step closer.\footnote{For the analogous construction
  with $P$ being conditional expectation (rather than finite rank
  projection), one can show that the $(\xi_n)$ are martingale
  differences.}  Second, orthogonality relations like
(\ref{eq:wiener-projection-ortho}) are equivalent to optimality in the
sense of least squares.  The MZ equation does not guarantee the stronger
orthogonality (\ref{eq:wiener-projection-ortho}) because it does not
guarantee optimal estimation of $x_{n+1}$ using $\Psi(x_n),
\Psi(x_{n-1}), \cdots$.

We refer to the projection $P_W$ and the associated
decomposition~(\ref{eq:wiener-projection}) as the \emph{Wiener
  projection}.  Two comments: first, the lack of (explicit) memory terms
in Eq.~(\ref{eq:pathProj}) is not surprising because we have simply
incorporated all relevant memory effects in the definition of $P_W$
itself, and also assumed the availability of that entire past history at
the initial time $n=0$, so there is nothing more for a memory term to
capture.  Second, though the subspace $W$ is defined in terms of
$M^{-1}$ and its powers, in practice one does not need to compute
$M^{-1}$ or $F^{-1}$ in working with $W$ as one can simply keep track of
the (recent) history in stepping forward the reduced model.  So our
formalism can be safely applied to dissipative dynamical systems, for
which $F^{-1}$ may be extremely unstable.

\medskip

In addition to the orthogonality~(\ref{eq:wiener-projection-ortho}), the
Wiener projection has the following properties:
\begin{enumerate}

\item \keqref{eq:tower} implies the existence of $h_0, h_1,\cdots$ such
  that~Eqs.~(\ref{eq:wiener-projection}) hold, and if the vectors
  $\cup_{k\geq0}M^{-k}\Psi$ are linearly independent, then the
  coefficients $(h_k)$ are unique.  (The coefficients $h_n$ may be
  ill-conditioned functions of the data if the basis functions are
  nearly degenerate.  This is an important but nontrivial issue, which
  we plan to explore in future work.)

\item The correlation matrices $\braket{\xi_m,\Psi_n}$ and
  $\braket{\xi_m,\xi_n}$ are functions of $m-n$, i.e., $\xi_m$ and
  $\Psi_n$ are jointly wide sense stationary.  (The process $(\Psi_n)$
  is stationary by assumption.)

\end{enumerate}

The first claim is a direct consequence of the preceding discussion.
For the second claim, first we show that $\xi_n = (Q_WM)^n\pi$ is wide
sense stationary: by taking adjoints in Eq.~(\ref{eq:tower}), we get
$P_WM^{\ell} = P_WM^{\ell}P_W$.  A short calculation\footnote{Since
  $P_WM^{\ell} = P_WM^{\ell}(P_W+Q_W) = P_WM^{\ell}P_W+P_WM^{\ell}Q_W$.
  Combined with $P_WM^{\ell} = P_WM^{\ell}P_W$, we have
  $P_WM^{\ell}Q_W=0$.  From this, we get $M^{\ell}Q_W =
  (P_W+Q_W)M^{\ell}Q_W = P_WM^{\ell}Q_W+Q_WM^{\ell}Q_W =
  Q_WM^{\ell}Q_W$.}  yields
\begin{equation}
  \label{eq:q-tower}
  M^{\ell}Q_W = Q_WM^{\ell}Q_W~,~~\ell\geq0.
\end{equation}
Repeated application of Eq.~(\ref{eq:q-tower}) yields
\begin{subequations}
  \label{eq:wiener-noise-is-nice}
  \begin{align}
    (Q_WM)^n\pi
    &= Q_WMQ_WMQ_W\cdots Q_WMQ_WM\pi\\
    &= M^{n-1}Q_WM\pi\\
    \label{eq:wiener-noise-is-nice-b}
    &= M^{n-1}\xi_1\\
    &= \xi_1\circ F^{n-1}~,~~n=1,2,\cdots.
  \end{align}
\end{subequations}
Thus, $\braket{\xi_m,\xi_n} = \braket{\xi_1\circ F^{m-1},\xi_1\circ
  F^{n-1}}$.  Since the probability distribution $\mu$ is $F$-invariant,
we have
\begin{align*}
  \braket{\xi_1\circ F^{m-1},\xi_1\circ F^{n-1}} 
  &= \int\xi_1(F^{m-1}(x))\cdot\xi_1(F^{n-1}(x))^T~d\mu(x)\\
  &= \int\xi_1(F^{m-n}(x))\cdot\xi_1(x)^T~d\mu(x)\\
  &= \braket{\xi_1\circ F^{m-n},\xi_1},
\end{align*}
i.e., $\xi_1,\xi_2,\cdots$ is wide sense stationary.

To see that $\braket{\xi_m,\Psi_n}$ is also a function of $m-n$, observe
\begin{subequations}
  \begin{align}
    \braket{\xi_m,\Psi_n}
    &= \braket{\xi_1\circ F^{m-1},\Psi_0\circ F^n}\\
    &=\braket{\xi_1\circ F^{m-n-1},\Psi_0},
  \end{align}
\end{subequations}
using Eq.~(\ref{eq:wiener-noise-is-nice}) and the invariance of $\mu$.
This can also be established by a more ``operator-theoretic'' argument:
observe
\begin{subequations}
  \begin{align}
    P_WM^{-m}(Q_WM)^n
    &= P_WM^{-m}M^{n-1}Q_WM\label{eq:12a}\\
    &= P_WM^{n-m-1}Q_WM.
  \end{align}
\end{subequations}
(Eq.~(\ref{eq:12a}) follows by repeated use of Eq.~(\ref{eq:q-tower})
with $\ell=1$.)  Using $\xi_n = M^{n-1}\xi_1$ (see
Eq.~(\ref{eq:wiener-noise-is-nice-b})) and the definition of $P_w$, we
see that $\braket{\Psi_m,\xi_n}$ is a function of $m-n$.

\subsection{Deriving NARMAX via rational approximations}
\label{sect:rat-narmax}

\keqref{eq:wiener-projection} would not reduce computational cost unless
the sum in $k$ can be truncated.  Simply keeping a small number of
terms, however, may not provide a good approximation.  Put another way,
to use \keqref{eq:wiener-projection} as the basis for model reduction,
it is necessary to find an effective way to parametrize the space of
filters $(h_n)$.  To do this, we use an idea from MZ
theory~\cite{zwanzig}.  Let
\begin{equation}
  \label{eq:z-transform}
  H(z)=\sum_{n\geq0}h_nz^{-n}
\end{equation}
denote the {\em $z$-transform} of $(h_n)$.  This is the discrete-time
analog of the Laplace transform; its properties are summarized in
\ref{app:wiener-and-z}.  The $z$-transforms $X(z)$, $\Psi(z)$, and
$\Xi(z)$ of $(x_n)$, $(\Psi_n)$, and $(\xi_n)$, respectively, are
similarly defined.  Then using the convolution property of the
$z$-transform (see \ref{app:wiener-and-z}), we have the formal relation
\begin{equation}
  \label{eq:z-wiener}
  X(z) = \Psi(Z)\cdot H(z) + \Xi(z).
\end{equation}
In applications of MZ theory to, e.g., statistical physics, rational
approximations of the transfer function $H(z)$ are frequently
effective~\cite{forster,zwanzig}.  This suggests the (uncontrolled)
approximation 
\begin{subequations}
  \begin{equation}
    \label{eq:A(z)}
    H(z)\approx B(z)/A(z),
  \end{equation}
  with
  \begin{equation}
    A(z) = z^p + a_{p-1}z^{p-1} \terms a_0
    \qquad\mbox{and}\qquad
    B(z) = b_r z^r \terms b_0~.
  \end{equation}
\end{subequations}
Neglecting convergence and other mathematical issues for now, if we
substitute the \emph{ansatz} $H(z) = B(z)/A(z)$ into
Eq.~(\ref{eq:z-wiener}), we obtain $X(Z) = \Psi(z)\cdot B(z)/A(z) +
\Xi(z)$.  This relation among $z$-transforms is equivalent to a
recurrence relation.  To see this, define $y_n = \sum_{n\geq0}
\Psi_{n-k}\cdot h_k$.  Then $Y(z) = \Psi(z)\cdot H(z)$, so that
$A(z)Y(z) = \Psi(z)\cdot B(z)$.  Inverting the $z$-transform yields
$y_{n} + a_{p-1}y_{n-1} \terms a_0y_{n-p} = \Psi_{n-p+r}\cdot b_r \terms
\Psi_{n-p}\cdot b_0.$ Summarizing, this suggests
Eq.~(\ref{eq:wiener-projection-main}) with the {\em ansatz} $H(z) =
B(z)/A(z)$ can be written
\begin{subequations}
  \label{eq:recursion}
  \begin{align}
    x_{n+1} = & ~y_n + \xi_{n+1},\\
    y_{n} + a_{p-1}y_{n-1} \terms a_0y_{n-p} 
    = & ~\Psi_{n-p+r}\cdot b_r \terms \Psi_{n-p}\cdot b_0.\label{eq:14b}
  \end{align}
\end{subequations}
If we set one column of $\Psi$ to be $f$ in \keqref{eq:narmax},
\keqref{eq:recursion} is essentially \keqref{eq:narmax}.

Modulo transients, Eq.~(\ref{eq:14b}) will correctly compute $y_n$
provided the recursion is stable in the sense that bounded $\Psi_n$ lead
to bounded $y_n$.  This holds if and only if the roots of the polynomial
$A(z)$ all lie strictly within the unit circle.  In this paper, we refer
to the condition $h_n\to0$ as the {\em decaying memory condition.}
Decaying memory is necessary for \keqref{eq:wiener-projection} to be
meaningful, for otherwise the reduced model would be sensitive to
information in the distant past.  We note decaying memory is necessary
but not sufficient for the overall numerical stability of the reduced
model.

If the decaying memory condition can be enforced, Eq.~(\ref{eq:14b})
provides an efficient way to compute the convolution in
Eq.~(\ref{eq:wiener-projection-main}), at a cost of not satisfying
Eq.~(\ref{eq:pathProj}) exactly.  As a result, there may be additional
memory-like corrections.  A detailed analysis of this is left for future
work.

\medskip

Eq.~(\ref{eq:z-wiener}) is purely formal in that in our context, where
the $(x_n)$, $(\Psi_n)$, and $(\xi_n)$ are stationary time series, the
$z$-transforms do not converge for any $z\in\C$.  A more careful
treatment uses the idea of power spectra.  As this is useful later in
the paper, we recall the notion here.

For a stationary stochastic process $(u_n)$, its {\em spectral power
  density} (or simply {\em power spectrum}) is the function
$S_{uu}(\theta) = \sum_{n=-\infty}^\infty C_{uu}(n)e^{in\theta}$, where
$C_{uu}(n)$ is the autocovariance function (ACF)
$\cov(u_n,\overline{u_0})$.  Similarly, for two stationary stochastic
processes $(u_n)$ and $(v_n)$, their cross power spectrum
$S_{uv}(\theta)$ is defined by $\sum_{n=-\infty}^\infty
C_{uv}(n)e^{in\theta}$, where $C_{uv}(n) = \cov(u_n,\overline{v_0})$ is
the cross correlation function (CCF).  In our context, we can view
$x_n$, $\Psi_n = \Psi(x_n)$, and $\xi_n$ are (possibly matrix-valued)
zero-mean (wide-sense) stationary time series satisfying $x_{n+1} =
\sum_{k\geq0}\Psi_{n-k}\cdot h_k + \xi_{n}$ with $\cov(x_m,\xi_n) =
\cov(\Psi_m,\xi_n)=0$ for all $n>m$.  Then, using the properties of
power spectra and $z$-transforms, one can show
\begin{equation}
  \label{eq:spectral}
  S_{xx}(\theta) = H^*(e^{-i\theta})S_{\Psi\Psi}(\theta)H(e^{-i\theta})
  + H^*(e^{i\theta})S_{\Psi\xi}(\theta)
  + S_{\xi\Psi}(\theta)H(e^{-i\theta})
  + S_{\xi\xi}(\theta).
\end{equation}
Eq.~(\ref{eq:z-transform}) typically does not converge for all $z\in\C$;
we assume the domain of convergence contains the unit circle, so that
Eq.~(\ref{eq:spectral}) makes sense.

\heading{Loss function and nonlinear regression.}
Eq.~(\ref{eq:recursion}) does not, by itself, fully specify a dynamical
model: to have a well-defined model, one needs to specify, e.g., the
statistics of the $(\xi_n)$.  For example, we can approximate $\xi_n$ by
a moving average of the form $d_qw_n \terms d_0w_{n-q}$, where the $w_n$
are independent $N(0,I)$ random vectors; this then gives a NARMA(X)
representation~(\ref{eq:narmax}).  Alternatively, one can prescribe the
properties of the $(\xi_n)$ implicitly by specifying the loss function
to be minimized, which we now discuss.  We observe that the rational
approximation above implies $p=q$, simplifying order selection.

Since Mori-Zwanzig aims to minimize the difference between the full and
reduced models with respect to the $L^2$ norm, a natural choice is to
minimize the mean squared error
\begin{equation}
  \label{eq:cost}
  \Ecal(a,b) = \frac1N\sum_{n=0}^{N-1}\Big\|\tilde{x}_{n+1} -
  \hat{x}_{n+1}\big(\tilde{\Psi}_1, \cdots, \tilde{\Psi}_n; a,
  b\big)\Big\|^2
\end{equation}
$a=(a_{p-1},\cdots,a_0)$ and $b=(b_q,\cdots,b_0)$ are the coefficients
of $A(z)$ and $B(z)$ in \keqref{eq:A(z)}, $(\tilde{x}_n)$ are data obtained from the full
model (say by simulation), and $\tilde{\Psi}_n = \Psi(\tilde{x}_n)$, and
where the one-step prediction $\hat{x}_n$ is here defined by
\begin{equation}
  \label{eq:one-step-prediction}
  \hat{x}_{n+1}(\tilde{\Psi}_1,\cdots,\tilde{\Psi}_n) =
  \sum_{k\geq0}\tilde{\Psi}_{n-k}\cdot h_k~.
\end{equation}
Because of the parametrization $H(z) = B(z)/A(z)$, the mean squared error
$\Ecal(a,b)$ depends {\it nonlinearly} on $a$ and $b$.  This leads to
two\footnote{In standard approaches to Wiener filtering, one makes use
  of the power spectra $S_{xx}$, $S_{x\psi}$, and $S_{\psi\psi}$ and
  their meromorphic continuations and solves the filtering problem by
  Wiener-Hopf techniques (see, e.g., \cite{kailath}).  In the context of
  data-driven modeling, direct minimization of $\Ecal(a,b)$ is more
  attractive because of the various sources of statistical error.}
possible approaches:
\begin{itemize}
\item {\em Nonlinear regression,} i.e., tuning $a$ and $b$ to minimize $ \Ecal(a,b)$ in \keqref{eq:cost}.
\item Finding $h_n$ directly by solving a (potentially very large)
  linear programming problem, then finding a good rational approximation
  $H(z) \approx B(z)/A(z)$.
\end{itemize}
In either case, we then fit a noise model to the residuals from the
nonlinear regression.

For high dimensional problems, the second approach is computationally
more challenging.  In this paper, we use the nonlinear regression
approach.  Numerical details are described in Sect.~\ref{sect:numerics}.

\heading{Multistep form and linear regression.}  Modulo transients,
\keqref{eq:recursion} is equivalent (see \ref{app:wiener-and-z}) to the
multistep recursion
\begin{equation}
  \label{eq:multistep}
  x_{n+p+1} + a_{p-1}x_{n+p}+\cdots+a_0x_{n+1} = \Psi(x_{n+r})\cdot b_r
  + \cdots+\Psi(x_n)\cdot b_0 + \bar{\xi}_{n+p+1},
\end{equation}
where $\bar{\xi}_{n+p+1} = \xi_{n+p+1} +
a_{p-1}\xi_{n+p}+\cdots+a_0\xi_{n+1} $.  Unlike \keqref{eq:recursion},
this formulation does not introduce any auxiliary variables.  The noise
$(\bar{\xi}_n)$ in \keqref{eq:multistep} is related to the $(\xi_n)$ in
\keqref{eq:recursion} by $S_{\bar{\xi}\bar{\xi}}(\theta) =
|A(e^{i\theta})|^2S_{\xi\xi}(\theta)$.  This means there is no simple
orthogonality relation between $\bar{\xi}_n$ and $\Psi_n$.  For these
reasons, \keqref{eq:multistep} is less convenient than
\keqref{eq:recursion} for model fitting.  Both require $p$ vectors
$x_1,\cdots,x_p \in \R^d$ as initial conditions.  In practice, these
initial conditions can have a measurable impact on noise models; we
discuss this and other implementation issues in
Sect.~\ref{sect:numerics}.

Eq.~(\ref{eq:multistep}) suggests an alternative loss function: compute
the one-step predictions using
\begin{equation} 
  \label{eq:narma-one-step-prediction}
  \hat{x}_{n+p+1} + a_{p-1}\tilde{x}_{n+p}+\cdots+a_0\tilde{x}_{n+1} =
  \Psi(\tilde{x}_{n+r})\cdot b_r + \cdots+\Psi(\tilde{x}_n)\cdot b_0,
\end{equation}
and minimizing the left and right hand sides, i.e.,
\begin{equation}
  \Ecal_*(a,b) = \frac1N\sum_{n=1}^{N-1-p}\Big\|\tilde{x}_{n+p+1} -
  \sum_{j=0}^{p-1} a_{j} \tilde{x}_{n+j+1} - \sum_{j=0}^{r}
  b_j\tilde{\Psi}_{n+j}\Big\|^2.
\end{equation}
One can then fit the residual by a noise model, e.g., by a power
spectrum method (see Sect.~\ref{sect:noise}) or a moving average model.
The difference between minimizing $\Ecal(a,b)$ in Eq.~(\ref{eq:cost})
and $\Ecal_*(a,b)$ above is that the latter entails only linear
regression, which can be computed very quickly when the number of time
lags is not large.  Also, whereas Eq.~(\ref{eq:one-step-prediction})
depends on the all available past history,
Eq.~(\ref{eq:narma-one-step-prediction}) depends only on the past $r$
steps.  However, minimizing $\Ecal_*(a,b)$ may produce such effective
models because it neglects long-range correlations in the data.

Finally, we observe that in Eq.~(\ref{eq:multistep}), if the sequence
$(\xi_n)$ is assumed to be IID Gaussian, the resulting model is what is
often referred to as the NARMA model in time series analysis (see, e.g.,
\cite{BD02,FY03}).  In this case, one can infer the coefficients $a$ and
$b$ by the conditional maximal likelihood method, which entails
minimizes the cost function
\begin{equation}
  \label{eq:cost1}
  \Ecal(a,b\mid \xi_1,\cdots,\xi_p ) = \frac1N\sum_{n=1}^{N-1-p}\Big\|\tilde{x}_{n+p+1} -
  \sum_{j=0}^{p-1}  a_{j} (\tilde{x}_{n+j+1}-\tilde \xi_{n+j+1}) - \sum_{j=0}^{r} b_j\tilde{\Psi}_{n+j}\Big\|^2~.
\end{equation}
In the above, the sequence $(\tilde \xi_{n})_{n>p}$ can be computed
recursively from data for each given pair of $(a,b)$.  This cost
function is similar to $\Ecal(a,b)$, and the optimization is similar to
the nonlinear regression above: instead of using $(h_n)$ above, one
computes the sequence $(\tilde \xi_{n})_{n>p}$ in each optimization step
(see \cite{CL15} for more details).

\subsection{Random dynamical systems and systems with delays}
\label{sect:wiener-rds}

Model reduction techniques are routinely applied to both deterministic
and random dynamical systems, as well as systems with delays.  The MZ
formalism applies to both random dynamical systems and to discrete-time
systems with bounded delays, as we now explain.  Our construction here
is related to the ``shift operator'' discussed in
\cite{berry2015nonparametric}.

We first explain how the MZ formalism applies to random dynamical
systems.  Consider the Euler-Maruyama discretization\footnote{The ideas
  we introduce here are quite general; we focus on Euler-Maruyama for
  the sake of simplicity.}  of a stochastic differential equation (SDE)
of the form $\dot{u}_t = f(u_t) + \dot{w}_t$:
\begin{equation}
  \label{eq:stochastic-burgers-scheme}
  u_{n+1} = u_n + f(u_n)\dt + \sqrt{\dt}~w_n~,
\end{equation}
where the $w_n$ are independent $N(0,I)$ random vectors.  The above has
the general form
\begin{equation}
  \label{eq:rds}
  u_{n+1} = F(u_n, w_n).
\end{equation}
Let $\underline{w} = (\cdots, w_{-1}, w_0, w_1, \cdots)$ denote the
entire history of the forcing.  A standard way to rewrite
Eq.~(\ref{eq:stochastic-burgers-scheme}) as an autonomous dynamical system
(\keqref{eq:system} above) is to augment the state $u_n$
with the history of the forcing $\underline{w}$.  In dynamical systems
language, such constructions are known as ``skew products.''  Here we
sketch the key ideas, and refer interested readers to, e.g.,
\cite{ledrappier-young,kifer2012ergodic,arnold} for mathematical
details (see also \cite{baxendale,kunita} for extensions to
stochastic differential equations).

Given a forcing sequence $\underline{w}$, we define
$\sigma(\underline{w})$ to be the sequence whose $n$th entry is
$w_{n+1}$, i.e., $\sigma(\underline{w})_n = w_{n+1}.$ In other words,
$\sigma(\underline{w})$ is sequence $\underline{w}$ shifted by 1 in
time.  If we shift $n$ times, then $w_n$ is moved into position 0, so
that $\pi_0(\sigma^n(\underline{w}) = w_n$, where $\pi_0(\underline{w})
= w_0.$

Using this notation, we can rewrite Eq.~(\ref{eq:rds}) as $u_{n+1} =
F(u_n, \pi_0(\sigma^n(\underline{w})))$, where $\underline{w}$ is a
given realization of the forcing sequence.  Now denote
$\underline{w}^{(n)} = \sigma^n(\underline{w})$; then
$\{\underline{w}^{(n)}~|~n\in{\mathbb Z}\}$ is a sequence of forcing
sequences, all related to each other by time shifts.  Then
\begin{subequations}
  \label{eq:skew-product}
  \begin{align}
    u_{n+1} &= F\Big(u_n, \pi_0(\underline{w}^{(n)})\Big) , \\
    \underline{w}^{(n+1)} &= \sigma(\underline{w}^{(n)}).
  \end{align}
\end{subequations}
Let $\spaceX$ be the space of all pairs $(u,\underline{w})$, i.e.,
$\spaceX$ is the state space of the discretized SDE augmented with its
forcing history.  Then Eq.~(\ref{eq:skew-product}) is a dynamical system
of the form \keqref{eq:system}, albeit one with an infinite-dimensional
state space $\spaceX$.  This does not prevent one from applying the
Mori-Zwanzig formalism.  In practice, one does not need to (and
generally cannot) keep track of the entire forcing history
$\underline{w}$, and a fragment of it is often sufficient.  Note that
within this framework, observation functions $\Psi$ can depend on both
the state $u_n$ and the forcing history $\underline{w}^{(n)}$.

Finally, we note that an invariant probability distribution $\mu$,
related in a natural way to the stationary distribution of
Eq.~(\ref{eq:stochastic-burgers-scheme}), can be constructed on this
augmented state space; see, e.g.,
\cite{ledrappier-young,kifer2012ergodic}.

\medskip

As for general delay terms, for example terms of the form
$\Psi(x_k,x_{k-\ell})$ for $\ell\leq L$ (which appear in our model for
the Burgers equation later in the paper), one can use a standard
construction: as in Eq.~(\ref{eq:system}), let $F$ be a given dynamical
system with state space $\spaceX$, and replace the state space $\spaceX$
by the $(L+1)$-fold cartesian product $\overline{\spaceX} =
\spaceX^{L+1}$, and replace $F$ by a map $\overline{F}$ on
$\overline{\spaceX}$ with
\begin{equation}
  \overline{F}(\overline{X}) ~~=~~ \overline{F}(X_0,\cdots,X_L) ~~=~~
  (F(X_0), X_0,\cdots,X_{L-1})
\end{equation}
for $\overline{X} = (X_0,\cdots,X_L) \in \overline{\spaceX}.$ This
constructions can be combined with the skew product construction
described earlier to handle stochastic systems with delays.

\section{Numerical implementation}
\label{sect:numerics}

This section addresses the problem of fitting models of the
form~(\ref{eq:recursion}) to data.  We take a two-step approach: we
first tune the coefficients $a$ and $b$ of the polynomials $A(z)$ and
$B(z)$, respectively, to minimize $\Ecal(a,b)$ in Eq.~(\ref{eq:cost});
we then use a stationary Gaussian process to model the residuals.
Sects.~\ref{sect:cascade} and \ref{sect:init-and-run} concern the
decaying memory constraint.  Sect.~\ref{sect:fitting} discusses other
details of optimization, and Sect.~\ref{sect:noise} noise modeling.

We have implemented the algorithms described here and the examples of
Sect.~\ref{sect:numAll} in Julia version~1.4~\cite{bezanson}.  For
numerical optimization, we used the {\tt NLopt.jl}
package~\cite{johnson}.  The source code is being prepared for public
release, and will be available at
\href{https://github.com/kkylin}{https://github.com/kkylin}~.

\subsection{Decaying memory constraint and the second-order cascade}
\label{sect:cascade}

To fit a model of the form Eq.~(\ref{eq:wiener-projection}) to data, we
will need to enforce the decaying memory condition $h_k\to0$ for two
reasons.  First, the decaying memory condition is necessary for the
reduced model to be meaningful.  Second, while we can compute one-step
predictions directly using Eq.~(\ref{eq:one-step-prediction}), either
directly or by the fast Fourier transform, the computational cost will
be quite high for high dimensional problems.  It would be much more
efficient if we can implement the convolution indirectly by making use
of Eq.~(\ref{eq:recursion}), i.e., compute the one-step prediction by
\begin{align}
  \label{eq:one-step-prediction-rational}
  \hat{x}_{n+1} =&   ~y_n, \\
  y_{n} + a_{p-1} y_{n-1} \terms
  a_0y_{n-p}   %
  =& ~\Psi(\tilde{x}_{n-p+r})\cdot b_r\terms
  \Psi(\tilde{x}_{n-p})\cdot b_0\nonumber.
\end{align}
But as discussed earlier, we need the decaying memory condition to
ensure these recursions will correctly compute $y_n~.$ The challenge is
that the loss function $\Ecal(a,b)$ is  highly nonlinear in $a$
and $b$.  Because the decaying memory condition involves the roots of
$A(z)$ in \keqref{eq:A(z)}, it exacerbates the problem.  Our general approach is to
reformulate Eq.~(\ref{eq:recursion}) so that the decaying memory
constraint becomes easier to implement, at the cost of making the cost
function highly non-convex.  We then fit reduced models to data using
this representation by numerical optimization.  We have found this to be
sufficient for the examples in this paper, though more work needs to be
done to ensure its robustness and efficiency for more general problems.

Consider a model of the form Eq.~(\ref{eq:recursion}) given
coefficients, and suppose for simplicity that $A(z)$ has real scalar
coefficients.
  We begin with the the observation that for a quadratic polynomial
$z^2+\alpha z+\beta$, its roots lie inside the unit disc if and only if
$(\alpha,\beta)$ lies inside the triangle in the $\alpha\beta$-plane
with vertices $(\pm2,1)$ and $(0,-1)$.  That is to say, for such an
$A(z)$, the decaying memory condition consists of three \emph{linear}
inequalities.  To make use of this observation for non-quadratic $A(z)$,
we factor $A(z)$ into a product of quadratic factors when $p=\deg(A)$ is
even, and quadratic factors and one linear factor if $p$ is odd, i.e.,
\begin{equation}
  \label{eq:z-transform-A}
  A(z) = \prod_{i=1}^{p/2} (z^2+\alpha_iz+ \beta_i)\qquad\mbox{or}\qquad
  A(z) = (z+\alpha_0)\prod_{i=1}^{\floor{p/2}} (z^2+\alpha_iz+ \beta_i).
\end{equation}
In this form, the decaying memory condition is naturally expressed as a
system of linear inequalities, which are easily imposed when performing
numerical optimization.

In view of the convolution theorem for $z$-transforms, the quadratic
factorization of $A(z)$ is equivalent to representing the linear filter
with transfer function $1/A(z)$ as a \emph{cascade of second-order
  filters}.  To see this, suppose (for simplicity) that $p=2s$.  We
introduce auxiliary variables $(z_i^n)$ for $i=1,\cdots,s$ (these
variables $z_i^n$ differ from the $z$ in $z$-transforms), and suppose
they satisfy
\begin{subequations}
  \label{eq:cascade}
  \begin{equation}
    \label{eq:cascade-main}
    \begin{array}{crcl}
      \mbox{\small Stage $1$}\qquad & z_{1}^{n}
      +\alpha_1z_{1}^{n-1}+\beta_1z_{1}^{n-2} &=& \Psi_{n-p+r}\cdot b_r
      + \cdots +\Psi_{n-p}\cdot b_0\\[2ex]
      \mbox{\small Stage $2$}\qquad &
      z_{2}^{n}+\alpha_2z_{2}^{n-1}+\beta_2z_{2}^{n-2} &=&
      z_{1}^n\\[1ex]
      \vdots\qquad&&\vdots&\\[1ex]
      \mbox{\small Stage $s$}\qquad
      &z_{s}^{n}+\alpha_{s}z_{s}^{n-1}+\beta_{s}z_{s}^{n-2} &=&
      z_{s-1}^n\\
    \end{array}.
  \end{equation}
  We claim that if
  \begin{equation}
    \label{eq:cascade-output}
    x_{n+1} = z_{s}^n + \xi_{n+1}
  \end{equation}
\end{subequations}
then Eq.~(\ref{eq:cascade}) is equivalent to Eq.~(\ref{eq:recursion}),
modulo transients.  To see this, observe that (neglecting initial
conditions) we have
\begin{align*}
  (1+\alpha_iz^{-1} +\beta_iz^{-2})Z_i(z) &= Z_{i-1}(z)~,\qquad i=2,\cdots,s\\
  (1+\alpha_1z^{-1} +\beta_1z^{-2})Z_1(z) &= z^{-p}\Psi(z)\cdot B(z).
\end{align*}
Putting it all together (and remembering $p=2s$) gives $Z_s(z) =
\Psi(z)\cdot B(z) / \Pi_{i=1}^s(z^2+\alpha_iz +\beta_i) = \Psi(z)\cdot
B(z)/A(z)$, and inverting $z$-transforms yields Eq.~(\ref{eq:cascade}).
The equivalence is up to transients because we have neglected initial
conditions in this discussion, and the argument is valid only if the
decaying meory condition holds.  The recursion in Eq.~(\ref{eq:cascade})
is explicit when $p\geq q$.  In the notation of
Eq.~(\ref{eq:recursion}), the output of the last stage gives $y_n$,
i.e., $y_n=z_{s}^n$.

\heading{Example.}  For $p=r=4$, we have two stages:
\begin{equation}
  \begin{array}{crcl}
    \mbox{\small Stage $1$}\qquad & z_{1}^{n}
    +\alpha_1z_{1}^{n-1}+\beta_1z_{1}^{n-2}  &=& 
    \Psi_{n}\cdot b_4 + ~\cdots~ + \Psi_{n-4}\cdot b_0\\[1ex]
    \mbox{\small Stage $2$}\qquad &
    z_{2}^{n}+\alpha_2z_{2}^{n-1}+\beta_2z_{2}^{n-2} &=&  z_{1}^n\\
  \end{array}
\end{equation}
In this case, it is easy to show directly that
\begin{equation}
  y_{n} + a_3y_{n-1} + ~\cdots~ + a_0y_{n-4} = \Psi_{n}\cdot b_4 +
  ~\cdots~ + \Psi_{n-4}\cdot b_0
\end{equation}
where $y_n = z_2^n$ and
\begin{equation}
  z^4 + a_3z^{3} + a_2z^{2}+a_1z+ a_0 =
  (z^2+\alpha_1z+\beta_1)\cdot(z^2+\alpha_2z+\beta_2)~,~~z\in\C.
\end{equation}
The corresponding reduced model can be written as a system
\begin{align*}
  x_{n+1} =& y_n + \xi_{n+1}\\[1ex]
  y_{n}  =& -\big(a_3y_{n-1}+~\cdots~+ a_0y_{n-4}\big) ~~+~~\big(\Psi_{n}\cdot b_4 +~\cdots~+ \Psi_{n-4}\cdot b_0\big)~.
\end{align*}
With $p=r=0$, we have a one-step (Galerkin) recursion $x_{n+1} =
\Psi_n\cdot b_0 + \xi_{n+1}$.  Similarly, with $p=r=1$, we have $x_{n+1}
= y_n+\xi_{n+1}$ and $y_{n} = -a_0y_{n-1} + \Psi_{n}\cdot b_1 +
\Psi_{n-1}\cdot b_0~,$ and setting $a_0=0$ yields $x_{n+1} =
\Psi_{n}\cdot b_1 + \Psi_{n-1}\cdot b_0+\xi_{n+1}$.

\subsection{Initializing and running cascade-form  models}
\label{sect:init-and-run} 

We use the cascade-form model to impose the decaying memory condition.
This is needed both for running fitted reduced models and, as we explain
later, for fitting models to data.  Here, we discuss how to initialize
and run such models.

Running the model to produce predictions entails carrying out the
recursions in Eq.~(\ref{eq:cascade}), at each point computing the predictors
$\Psi_n = \Psi(x_n)$ with $x_{n}=y_{n-1} + \xi_n = z_{s}^{n-1} + \xi_n$.
Though derived from Eq.~(\ref{eq:recursion}), Eq.~(\ref{eq:cascade}) is quite
different in form.  Here we examine Eq.~(\ref{eq:cascade}) more closely, to
clarify the flow of information in the algorithm and other details.

It is useful to first visualize Eq.~(\ref{eq:cascade}) as a computation
graph, a fragment of which is shown here:

\begin{tikzcd}
  &&\mbox{Step $n-2$}\arrow[dotted]{d}&\mbox{Step
    $n-1$}\arrow[dotted]{d}&\mbox{Step $n$}\arrow[dotted]{d}\\
  \mbox{Stage $s-1$} &\cdots\arrow[dotted]{r}
  &z_{s-1}^{n-2}\arrow[dotted]{r}\arrow[dotted]{d}\arrow[bend
    right=20,dotted]{rr}
  &z_{s-1}^{n-1}\arrow[dotted]{r}\arrow[dotted]{d}
  &\boxed{z_{s-1}^{n}}\arrow{d}\\
  \mbox{Stage $s$} &\cdots\arrow[dotted]{r}
  &\boxed{z_{s}^{n-2}}\arrow[dotted]{r}\arrow[bend
    right=20,swap]{rr}{\beta_s} &\boxed{z_{s}^{n-1}}
  \arrow{r}{\alpha_s} &\boxed{z_{s}^{n}}
\end{tikzcd}

\noindent
(For legibility, we have drawn the edges going into $z_s^n$ as solid
lines; all others are dotted.)  The variable $z^n_s$ at time $n$ and
stage $s$ depends on the corresponding variable $z^n_{s-1}$ in the
previous stage, as well as the two previous steps ($z^{n-1}_s$ and
$z^{n-2}_s$) in the same stage.

Once we have initial conditions, Eq.~(\ref{eq:cascade}) can be iterated to
generate sample paths. The first thing is then to find the initial values
$z_{i}^{p-1}$ and $z_{i}^{p-2}$ for $i=1,\cdots,r$ from the given data
$\tilde{x}_1,\cdots,\tilde{x}_N$.  An effective procedure is suggested
by the computation graph: we set
\begin{equation}
  \tilde{y}_0 = \tilde{x}_1~,~~ \tilde{y}_1 = \tilde{x}_2~,~~
  \cdots~,~~ \tilde{y}_{p-1} = \tilde{x}_p
\end{equation}
in the notation of Eq.~(\ref{eq:recursion}) and Eq.~(\ref{eq:cascade}).
Assuming the coefficients $\alpha_i$ and $\beta_i$ have already been
determined, the computation graph shows that knowing the values at stage
$s$ for $n=1,2,\cdots,p$ (which is the same as knowing $y_1,\cdots,y_p$
in Eq.~(\ref{eq:recursion})) allows one to solve for the values at stage
$s-1$ for $n=3,4,\cdots,p$.  Iterating, this means we can determine
$z_i^{p-1}, z_i^{p}$ for all stages $i$.  From this, it is also
straightforward to see that if $y_0=\cdots=y_{p-1}=0$, then
$z_{i}^{p-1}=z_{i}^{p-2}=0$ for $i=1,\cdots,r$, so that the initial
conditions for Eq.~(\ref{eq:cascade}) are uniquely determined by those of
Eq.~(\ref{eq:recursion}).

Once the initial data have been determined and noise generated (as
described in Sect.~\ref{sect:wiener-projections}), the recurrence
relations~(\ref{eq:cascade}) can be iterated to generate predictions
from the reduced model.

\subsection{Fitting models to data}
\label{sect:fitting}

We now describe our overall optimization strategy:
\begin{enumerate}
\item From the time series $\tilde{x}_1,\cdots,\tilde{x}_N$, compute the
  observations $\tilde{\Psi}_n = \Psi(\tilde{x}_n)$.
\item For given parameter vectors $\alpha, \beta, b$, use the initial
  values $\tilde{x}_1,\cdots,\tilde{x}_p$ to determine the initial
  values $z^{p-2}_i$, $z^{p-1}_i$, $i=1,\cdots, r$, for
  Eq.~(\ref{eq:cascade}).
\item\label{optim-step3} Generate one-step predictions $\hat{x}_{n+1}$
  by Eq.~\eqref{eq:one-step-prediction-rational} for $n=p,\cdots,N$,
  where $H(z) = B(z)/A(z)$.

\end{enumerate}
In the cascade representation, the MSE has the form
\begin{subequations}
  \begin{equation}
    \label{eq:eprime}
  \Ecal'(\alpha,\beta,b) = \frac1N\sum_{n=p+1}^N\Big\|\tilde{x}_{n+1} -
  \hat{x}_{n+1}\big(\tilde{\Psi}_1, \cdots, \tilde{\Psi}_n; \alpha,
  \beta,b\big)\Big\|^2
  \end{equation}
  with the decaying memory constraints
  \begin{equation}
    \beta_i\leq1\qquad\mbox{and}\qquad\beta_i\geq\pm\alpha_i-1~.
\end{equation}
  (This says that $(\alpha_i,\beta_i)$ lies within a triangle in the
  $\alpha$-$\beta$ plane with vertices $(\pm2,1),(0,-1)$.  As asserted
  in Sect.~\ref{sect:cascade}, one can check that this is equivalent to
  the roots of $z^2+\alpha_iz +\beta_i$ lying in the unit disc.)
\end{subequations}
This can be minimized by direct optimization.  One then finds the
residuals
\begin{equation}
  \tilde{\xi}_{n} = \tilde{x}_{n+1} - \hat{x}_{n+1}\big(\tilde{\Psi}_1, \cdots,
  \tilde{\Psi}_n; \alpha, \beta,b\big)
\end{equation}
and fit a noise model as before.  

One can actually further reduce the dimensionality of the optimization
problem; this is described below.  But first, we note that Step (iii)
above is more efficiently implemented by iterating
\begin{equation}
  \label{eq:cascade-one-step-prediction}
  \begin{array}{crcl}
    \mbox{\small Stage $1$}\qquad & z_{1}^{n}
    +\alpha_1z_{1}^{n-1}+\beta_1z_{1}^{n-2} &=&
    \tilde{\Psi}_{n-p+r}\cdot b_r + \cdots +\tilde{\Psi}_{n-p}\cdot
    b_0\\[2ex]
    \mbox{\small Stage $2$}\qquad &
    z_{2}^{n}+\alpha_2z_{2}^{n-1}+\beta_2z_{2}^{n-2}
    &=& z_{1}^n\\[1ex]
    \vdots\qquad&&\vdots&\\[1ex]
    \mbox{\small Stage $s$}\qquad
    &z_{s}^{n}+\alpha_{s}z_{s}^{n-1}+\beta_{s}z_{s}^{n-2}
    &=& z_{s-1}^n\\[2ex]
    \mbox{\small Output}\qquad&\hat{x}_{n+1} &=& z_{s}^n\\
  \end{array}
\end{equation}
Modulo transients (see ``initial conditions'' below), this computes the
convolutions in Eq.~(\ref{eq:one-step-prediction}).  Note this iteration
can only be carried out if $\alpha,\beta$ satisfy the decaying memory
condition.

To further reduce the dimensionality of the nonlinear optimization
problem, we observe that for given $(\alpha,\beta)$, the function
$b\mapsto\Ecal'(\alpha,\beta,b)$ can be minimized by linear
regression.
For a given $(\alpha,\beta)$, we thus define $\hat{b}(\alpha,\beta)$ to
be the (unique) minimizer of $b\mapsto\Ecal'(\alpha,\beta,b)$.  We then
minimize $\Ecal'(\alpha,\beta,\hat{b}(\alpha,\beta))$ by nonlinear
optimization.  For the examples in this paper, this is done using the
BOBYQA algorithm\cite{powell} as implemented in the NLopt
package~\cite{johnson}.

\heading{Initial conditions.} The recursions in
Eq.~(\ref{eq:cascade-one-step-prediction}), viewed as a system of
non-autonomous linear 
  recurrences with $\tilde{\Psi}_n$ as time-dependent forcing, have their own initial
conditions.  Neglecting these ``internal'' initial conditions during
fitting, for example by setting them all to zero, can lead to worse
fits.  Without accounting for initial conditions, the residuals also
exhibit longer transients before approaching stationarity, which can
complicate the construction of noise models.

To estimate initial conditions for
Eq.~(\ref{eq:cascade-one-step-prediction}), we exploit the linearity of
Eq.~(\ref{eq:cascade-one-step-prediction}) in the variables $z_i^n$ by
decomposing the $z_i^n$ into the sum of a homogeneous solution
$z_i^{h,n}$ and a particular solution $z_i^{p,n}$, with $z_i^{p,n}$
satisfying Eq.~(\ref{eq:cascade-one-step-prediction}) with zero initial
conditions and $z_i^{h,n}$ solving
Eq.~(\ref{eq:cascade-one-step-prediction}) with $\tilde{\Psi}_n\equiv0$.
(In linear systems theory, these are the ``zero state response'' and
``zero input response,'' respectively.)  This leads to a linear
regression problem for the initial values $\{z_i^{h,0},z_i^{h,1}\mid
i=1,\cdots,s\}$, which can be solved jointly with the computation of
$\hat{b}(\alpha,\beta)$ via linear regression.

\heading{Remarks on optimization and related issues.}
\begin{itemize}

\item\emph{Cascade-form models, decaying memory, and optimization.}
  Eq~(\ref{eq:cascade}) enables us to impose the decaying memory
  constraint reliably during optimization.  However, the decomposition
  of $A(z)$ into quadratic factors introduces a symmetry: the value of
  the loss function is invariant when the quadratic factors are
  permuted.  This means there are many equivalent global minima, which
  introduce many saddles into the landscape.  While any of the symmetric
  global minima will give equivalent reduced models, the presence of the
  saddles can potentially slow down optimizers.

\item\emph{Other optimization strategies.}  For simplicity, we have
  opted for direct nonlinear minimization of
  $\Ecal'(\alpha,\beta,\hat{b}(\alpha,\beta))$ in this paper.  It may be
  possible to improve the efficiency of the optimization by exploiting
  the structure of Eq.~(\ref{eq:cascade}) or the multistep
  representation (\keqref{eq:multistep} above) by using, e.g., iterative
  least squares.

\item\emph{An implementation detail.}  For interested readers, we
  describe the computation of $\hat{b}(\alpha,\beta)$ by linear
  regression.  We run the matrix version of the recursion
\begin{equation}
  \label{eq:matrix-cascade}
  \begin{array}{crcl}
    \mbox{\small Stage $1$}\qquad & Z_{1}^{n}
    +\alpha_1Z_{1}^{n-1}+\beta_1Z_{1}^{n-2} &=& \tilde{\Psi}_{n}\\[2ex]
    \mbox{\small Stage $i>1$}\qquad &
    Z_{i}^{n}+\alpha_2Z_{i}^{n-1}+\beta_2Z_{i}^{n-2} &=& Z_{i-1}^n,\\
  \end{array}
\end{equation}
for $i=2,\cdots,r$ and setting $Y_n = Z^n_s$; this is a matrix version
of Eq.~(\ref{eq:cascade-one-step-prediction}) with $q=0$ and $b_0=I$.
The resulting $Y_n$ and $Z_i^n$ are matrix-valued, with the same shape
as $\tilde\Psi_n$.  By exploiting the commutativity of the convolution
operators defined by $B(z)$ and $1/A(z)$, one can show that the desired
one step prediction is given by
\begin{equation}
  \hat{x}_{n+1} = Y_{n-p+q}\cdot b_q + \cdots + Y_{n-p}\cdot b_0~.
\end{equation}
Combining this with the definition of $\Ecal'(\alpha,\beta,b)$ lets us
compute $\hat{b}(\alpha,\beta)$ via linear regression.

\end{itemize}

\subsection{Noise model}
\label{sect:noise}

To construct a stochastic process $\eta_n$ to model the residuals
$\xi_n$, there are a few standard options:
\begin{enumerate}

\item moving average representation, i.e., $\eta_n = d_qw_n \terms
  d_0w_{n-q}$ with independent $w_i\sim N(0,I)$;

\item estimating the power spectrum of $\xi_n$ and generating a
  stationary Gaussian process matching the power spectrum;

\item constructing a linear SDE and fitting it to $\xi_n$ by, e.g.,
  maximum likelihood.

\end{enumerate}
In earlier work, we have used a moving average representation together
with a MLE to infer the coefficients $a$ and $b$ simultaneously with the
coefficients of the moving average.  In this paper, because we want to
compare nonlinear regression with other approaches, the power spectrum
method was found to be simpler.

After finding optimal values for $a_i$, $b_i$, and the initial $y_i$, we
fit a stationary Gaussian process $\eta_n$ to the residuals
$\tilde{\xi}_{n+1} = \hat{x}_{n+1}-\tilde{x}_{n+1}$, by a random Fourier
series approximating a Wiener integral:
\begin{equation}
  \label{eq:random-fourier}
  \eta_n = \frac1{\sqrt{2\pi}}\sum_{j=0}^{M-1}
  f(j\dtheta)~e^{-inj\dtheta}~w_j\sqrt{\dtheta}
  \qquad\xrightarrow[M\to\infty]{{\mathcal D}}\qquad
  \frac1{\sqrt{2\pi}}\int_{0}^{2\pi}
  f(\theta)~e^{-in\theta}~\dot{W}_\theta~d\theta,
\end{equation}
where $\dtheta = 2\pi/M$, the $w_j$ are independent standard normal
random variables (in the complex case, $Re(w_j)$ and $Im(w_j)$ are
independent with variance $1/2$), $\dot{W}_\theta$ is white noise on the
circle $S^1$, and $f(\theta)$ is a square root of the spectral power
density, i.e., $S_{\xi\xi}(\theta) = f(\theta)f(\theta)^*$.  When
$\eta_n$ takes on values in $\R^d$, then $S_{\xi\xi}$ and $f$ are
$d\times d$ matrices and $\dot{W}_\theta$ is $d$-dimensional.  The power
spectrum can be estimated from data by the periodogram method (see,
e.g.,~\cite{press2007numerical} and references therein).  More efficient
and accurate sampling methods are available~\cite{cameron2003relative},
but we have found the random Fourier series above to be sufficient the
residuals $(\xi_n)$ are relatively small, as occurs in many examples
(including ours).  Whatever the method, the resulting reduced models
will only satisfy the orthogonality conditions approximately.

\section{Examples}
\label{sect:numAll}

We now consider two concrete examples.  In addition to illustrating the
methods described in earlier sections, there are two specific questions
we would like to address:
\begin{itemize}

\item How effective is the model reduction method based on nonlinear
  regression (as described in Sect.~\ref{sect:rat-narmax})?

\item How does the nonlinear regression compare to the linear regression
  described in Sect.~\ref{sect:rat-narmax}?

\end{itemize}
We would also like to see how the least squares based nonlinear regression
compare to the MLE used in \cite{LLC17}.

\subsection{Kuramoto-Sivashinsky (KS) PDE}
\label{sect:ks}
The KS equation  
\begin{equation}
  \label{eq:ks}
  U_t + UU_x + U_{xx} + U_{xxxx} = 0
\end{equation}
is a prototypical model of spatiotemporal chaos.  Here, we consider
\keqref{eq:ks} with $0\leq x\leq L$ and periodic boundary conditions.
In Fourier variables $u_k(t)$, \keqref{eq:ks} is
\begin{equation}
  \label{eq:ks-fourier}
  \dot{u}_k = -\frac{i\lambda_k}2\sum_\ell u_\ell u_{k-\ell} +
  (\lambda_k^2-\lambda_k^4)u_k~,~~\lambda_k = \frac{2\pi k}{L}.
\end{equation}
The lowest $\approx L/2\pi$ modes are linearly unstable.  This long-wave
instability and its interaction with the quadratic nonlinearity lead to
sustained chaotic behavior, with positive Lyapunov exponents and
exponential decay of correlations~\cite{hyman}.  NARMAX modeling of
\keqref{eq:ks} was studied in~\cite{LLC17}, using likelihood-based
parameter estimation and a slightly different form of NARMAX.  Here, we
use the least squares procedure.  Following \cite{LLC17}, we set
$L\approx 21.55$, leading to $3$ linearly unstable modes and a maximum
Lyapunov exponent of $\approx 0.04$ (Lyapunov time $\approx 25$).  In
this regime, time correlation functions exhibit complex oscillations
instead of the simple exponential decay often seen in strongly chaotic
systems (Fig.~\ref{fig:ks-stats}(a)), providing a nontrivial testbed for
model reduction.

\begin{figure}[tb!]
  \begin{center}
    \begin{tabular}{c@{\hskip -0.4in}c@{\hskip -0.2in} c}
      \begin{tabular}{l}
        \resizebox{!}{2in}{\includegraphics*[bb=0in 0.1in 5in 3in]{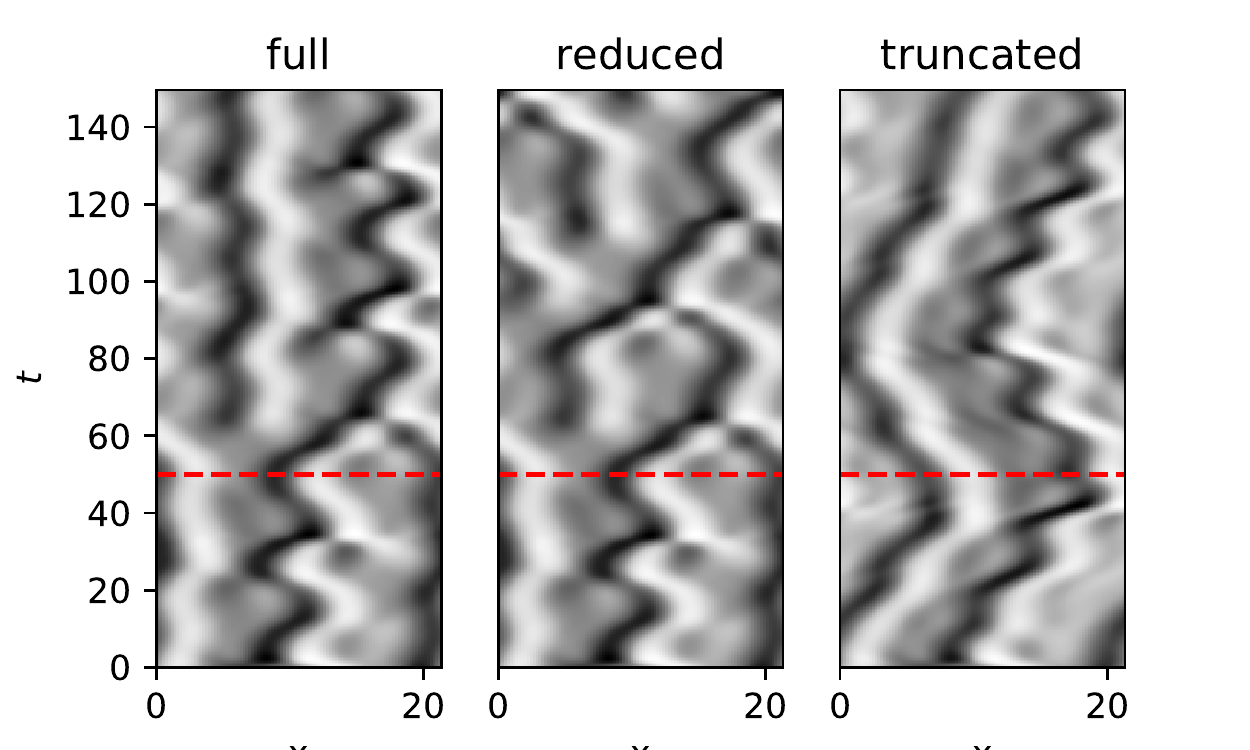}}\\[-2ex]
        \hspace{0.77in}$x$\hspace{0.9in}$x$\hspace{0.9in}$x$\\[-2ex]
      \end{tabular} &

      \begin{tabular}{r@{\hskip 0pt}c}
        \rotatebox{90}{\footnotesize\hspace{0.4in}$k=2$}&
        \resizebox{!}{1in}{\includegraphics{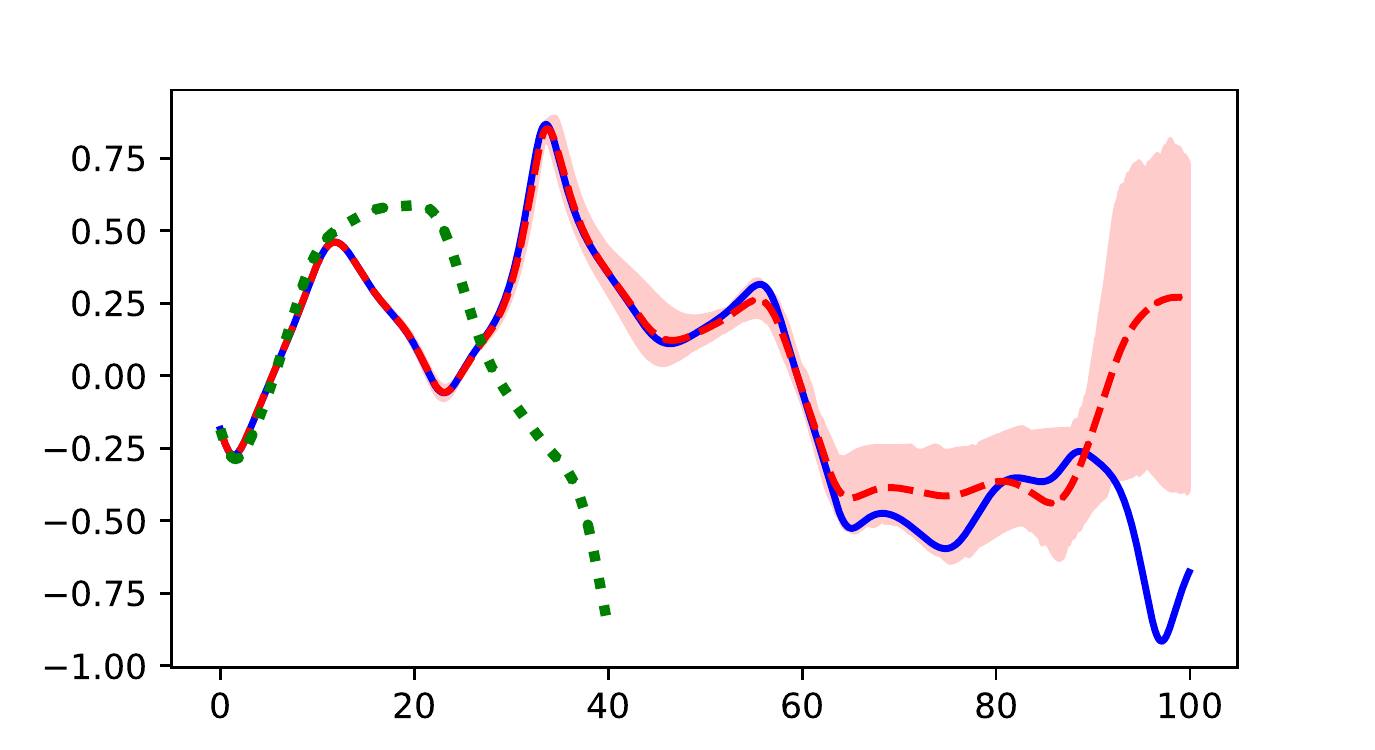}}\\[-1ex]
        \rotatebox{90}{\footnotesize\hspace{0.4in}$k=5$}&
        \resizebox{!}{1in}{\includegraphics{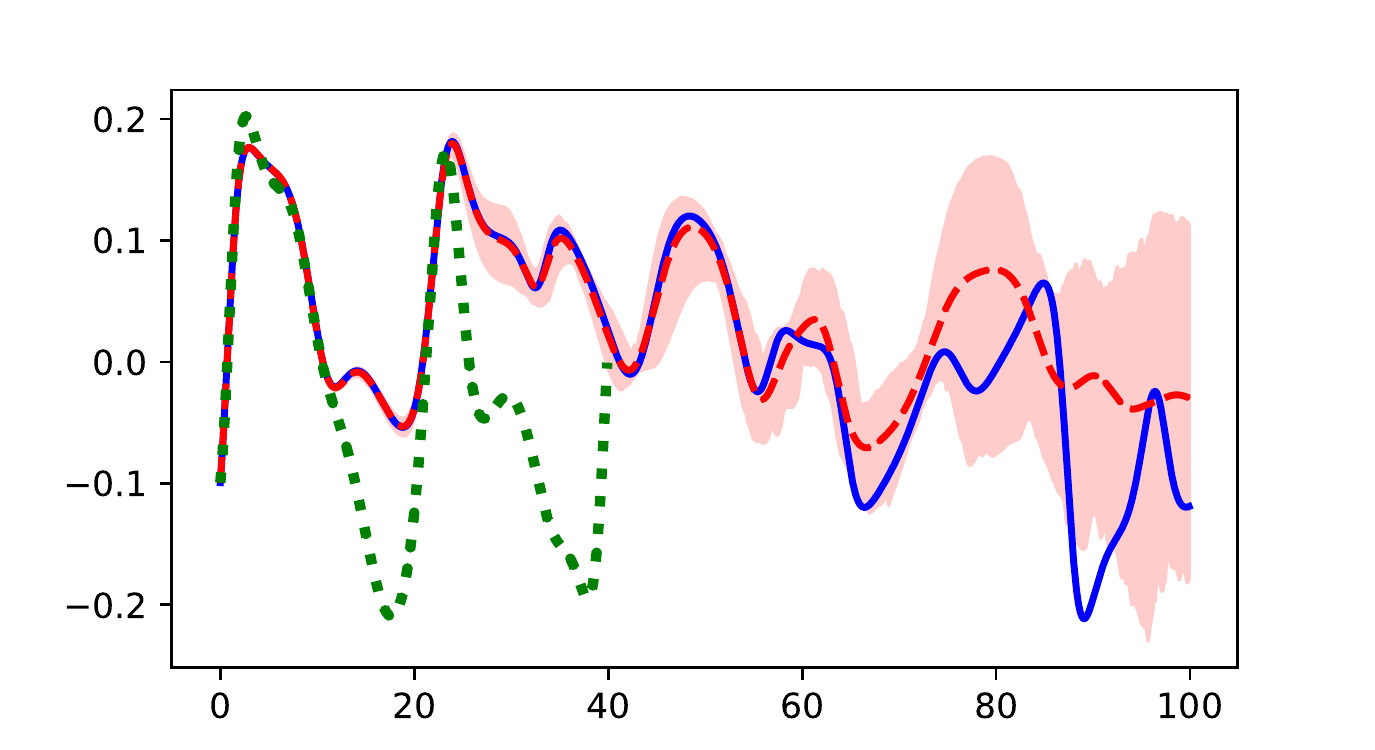}}\\[-1ex]
        &{\footnotesize $t$}\\
      \end{tabular} &
      \begin{tabular}{c}
        \resizebox{!}{0.6in}{\includegraphics*[bb=1.4in 1in 3in 2in]{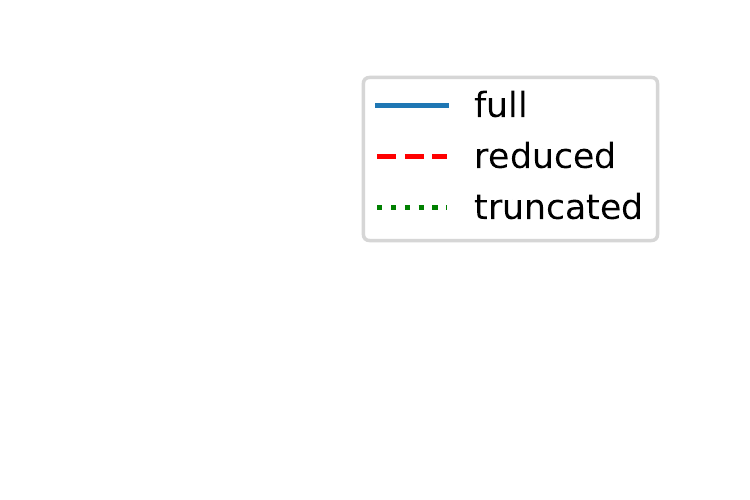}}\\
        \resizebox{!}{1.3in}{\includegraphics{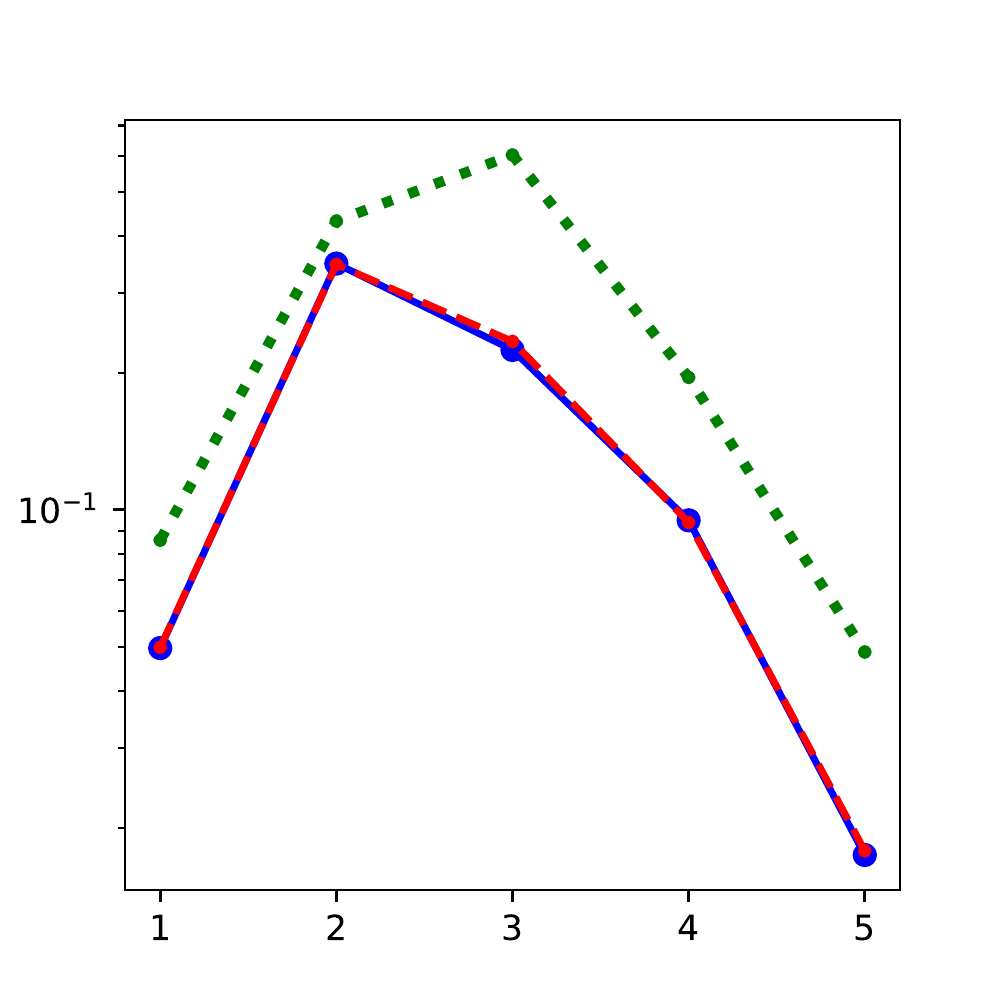}}\\[-2ex]
        {\footnotesize $k$}\\
      \end{tabular}\\
      (a) Spacetime views of KS solutions &
      (b) Trajectories of $Re(u_k(t))$& \hspace{6mm}
      (c) Energy $\braket{|u_k|^2}$ \\ 
    \end{tabular}~\vspace{-0.2in}



  \end{center}
  \caption{KS solutions.  Panel (a) shows results computed using the
    108-mode truncation ($\dt=10^{-3}$) {\em (left)}, the 5-mode reduced
    model ($\dt=0.1$) {\em (middle)}, and the 5-mode truncation
    ($\dt=10^{-3}$) {\em (right)}.  In (b), we plot two Fourier modes as
    functions of time, with 90\% confidence intervals for the reduced
    model.  Panel (c) shows the energy spectrum.} 


  \label{fig:ks-solutions}
\end{figure}

\keqref{eq:ks} is readily solved by truncating the Fourier series,
provided the cutoff is large enough.  Here, we take as full model the
108-mode truncation; numerical tests show that KS statistics are
insensitive to the cutoff beyond this.  Fig.~\ref{fig:ks-solutions}(a)
shows a sample solution of \keqref{eq:ks} using this 108-mode truncation
(``full'').  By comparison, the 5-mode truncation with the same initial
conditions (``truncated'') diverges rapidly, and fails to reproduce the
energy spectrum (Fig.~\ref{fig:ks-solutions}(c)).

\heading{Reduced model.}  To construct a reduced model using the lowest
$K=5$ Fourier modes, we follow the procedure outlined in
Sect.~\ref{sect:wiener-projections}.  The first step
is to generate data from the full model, which we do by numerically
integrating the 108-mode truncation using a 4th-order exponential
time-differencing Runge-Kutta (ETDRK4) method~\cite{CM02,KT05} with
timestep $\dt=10^{-3}$, for $10^8$ steps.  We observe the first $K=5$
Fourier modes at every 100 steps; the observation interval $\delta =0.1$
is the timestep for the reduced model.  We drop the first half of the
data to ensure stationarity.

We use the form of the reduced model in \keqref{eq:recursion} with $x_n$
corresponding to $u^n= (u^n_1,u^n_2,\dots, u^n_K)$; see
\appref{eq:ks-ansatz} for a detailed description of the model.  To
select the orders $p$ and $r$, we tried a variety of small values until
a combination is found that produces a stable reduced model. For the
function $\Psi(u)$, we use three groups of functions:
\begin{equation} \label{eq:Psi_KSE}
\begin{aligned}%
  \Psi^a_{n-j} &= u^{n-j}~,~~\\
  \Psi^b_{n-j} &= R^{\dt}(u^{n-j})~,\\
  \Psi^c_{n-j,k} &=  \sum_{\substack{ |k-l|\leq K, K< |l| \leq 2K \\ \text{ or }  |l|\leq K, K< |k-l| \leq 2K} }\widetilde u^{n-1}_l \overline{\widetilde u^{n-j}_{k-l}} \text{ for} ~~k=1,\cdots,K.
  \end{aligned}
 \end{equation}
Here the first two groups $\Psi^a$ and $\Psi^b$ come from the Galerkin
truncation. The third group in form of $\Psi^c$ represents that
interaction between the unresolved high modes and the resolved low
modes, in which the high modes $\widetilde u$, defined in
\appref{eq:ks-ansatz-d}, is motivated by the theory of approximate
inertial manifolds.  In terms of the formalism of
Sect.~\ref{sect:wiener-stuff}, the observation function $\Psi(u)$ is a
$K\times(2K+K^2)$ matrix whose entries consist of the terms given above,
where $K$ is the number of relevant Fourier modes.
Here, we use $K=5$.

Finally, the reduced model is fit to data by the procedure outlined in
Sect.~\ref{sect:numerics}.  As was found in~\cite{LLC17}, not all
combinations of $p$ and $r$ lead to stable reduced models.  Indeed, we
have experimented with ``replaying'' the residuals, i.e., compute the
residuals $\tilde{\xi}_n$ as in Sect.~\ref{sect:numerics}, then running
the reduced model with $\tilde{\xi}_{n+1}$ in place of the noise term.
In the absence of round-off, one would simply obtain $x_n =
\tilde{x}_n$, i.e., reconstruct the original time series.  Instead, for
some choices of $(p,r)$, round-off errors were rapidly amplified.  Here,
we use the pair $p=r=3$, which is found to strike a balance between
accuracy and efficiency.  As measured by the product of the mode and
step counts, the reduced model represents an over $100$-fold reduction
in computational cost.

\heading{Results.}  Fig.~\ref{fig:ks-solutions}(a) compares the full
model (``full''), the reduced model with $p=r=3$ (``reduced''), and the
5-mode truncation with $\dt=10^{-3}$ (``truncated'').  As one can see,
the reduced model reproduces the full solution up to $t\gtrsim50$, about
$1.8\times$ the Lyapunov time, consistent with~\cite{LLC17}.  In
contrast, the 5-mode truncation is accurate for a fraction of that time.
Fig.~\ref{fig:ks-solutions}(b) takes a closer look at selected Fourier
modes.  For the reduced model, 100 independent realizations are run, and
the resulting ensemble is used to estimate confidence intervals.  Shown
is the mean (dashed, red), and 90\% confidence intervals.  Though the
noise terms have amplitudes $\leq10^{-4}$, they are rapidly amplified by
exponential separation of trajectories due to the long-wave instability
in KS.  Consistent with Fig.~\ref{fig:ks-solutions}(a), the mean follows
the true trajectory up to $t\approx40$, at which point they begin to
diverge.  In contrast, the 5-mode truncation diverges by $t\approx20$.
Moreover, even when the confidence interval starts to widen, it
continues to provide useful bounds for some time.  Eventually the
ensemble approaches statistical steady state, and the ensemble mean
converges toward its expected value.  Fig.~\ref{fig:ks-solutions}(c)
compares the energy spectra $\braket{|u_k|^2}$: while the reduced model
correctly predicts the spectrum, the 5-mode truncation produces
fluctuations that are too large.

We note that while the noise terms are small in amplitude (see
Fig.~\ref{fig:ks-powerspectra}), we could not have constructed the
confidence intervals in Fig.~\ref{fig:ks-solutions} without them.
Moreover, we conducted numerical experiments without the noise terms.
The results (data not shown) show that the reduced models do
considerably worse at all tasks, and at least for some choices of
$(p,r)$ the solutions converge quickly to 0.

\begin{figure}[tb!]
  \begin{center}
    \begin{tabular}{r@{\hskip 0pt}c@{\hskip -6pt}c}
      &{\footnotesize $k=2$} & {\footnotesize $k=5$}\\
      \rotatebox{90}{\footnotesize\hspace{0.5in}ACF}&
      \resizebox{!}{1.5in}{\includegraphics{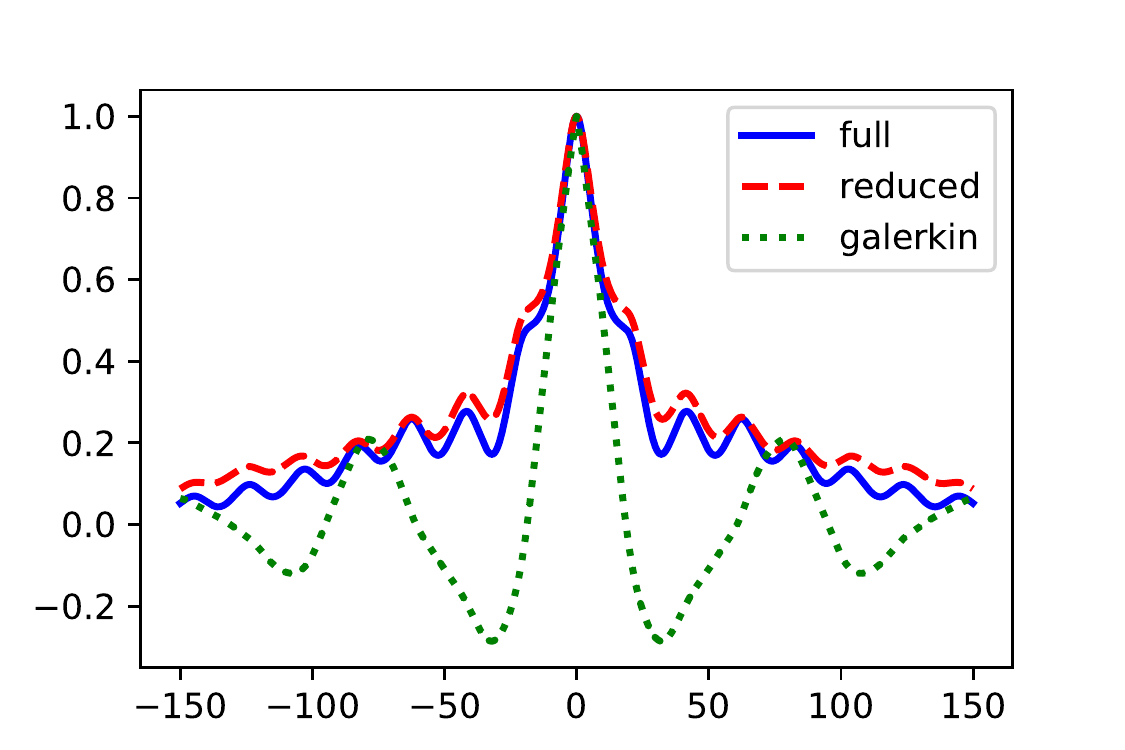}}&
      \resizebox{!}{1.5in}{\includegraphics{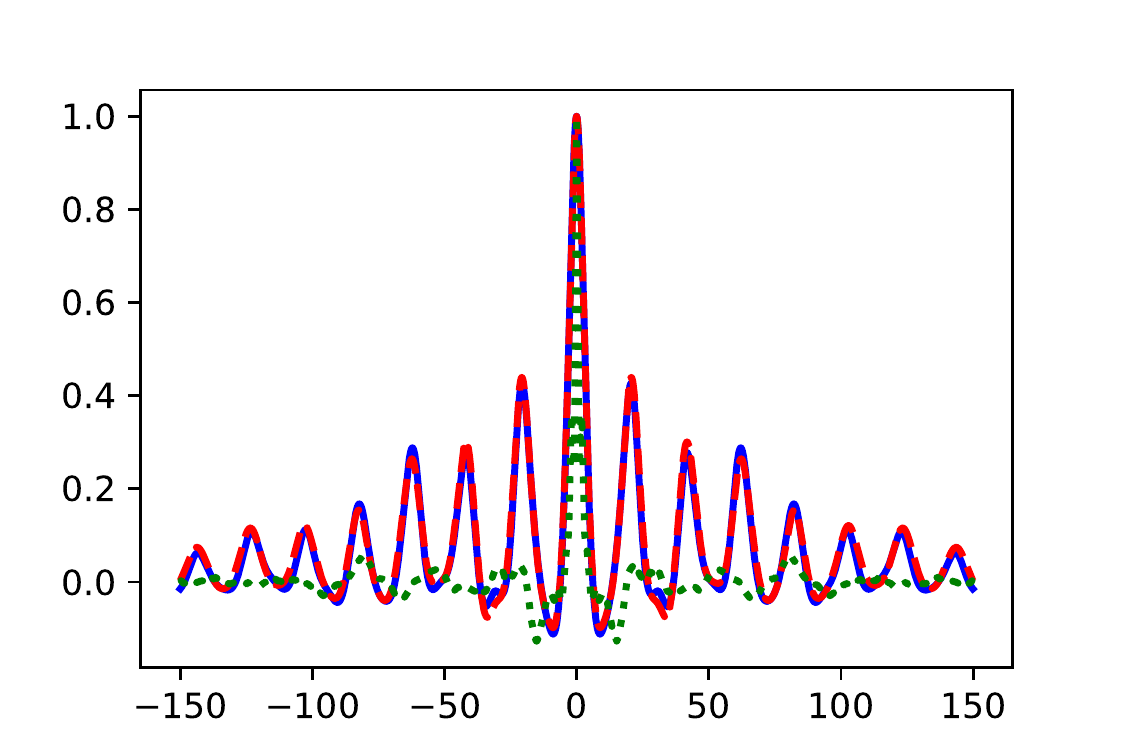}}\\[-1ex]
      &{\footnotesize Time lag}&{\footnotesize Time lag}\\
      \multicolumn{3}{c}{(a) Autocovariance functions $Re(u_k(t))$}\\[1ex]
      \rotatebox{90}{\footnotesize\hspace{0.5in}CCF}&
      \resizebox{!}{1.5in}{\includegraphics{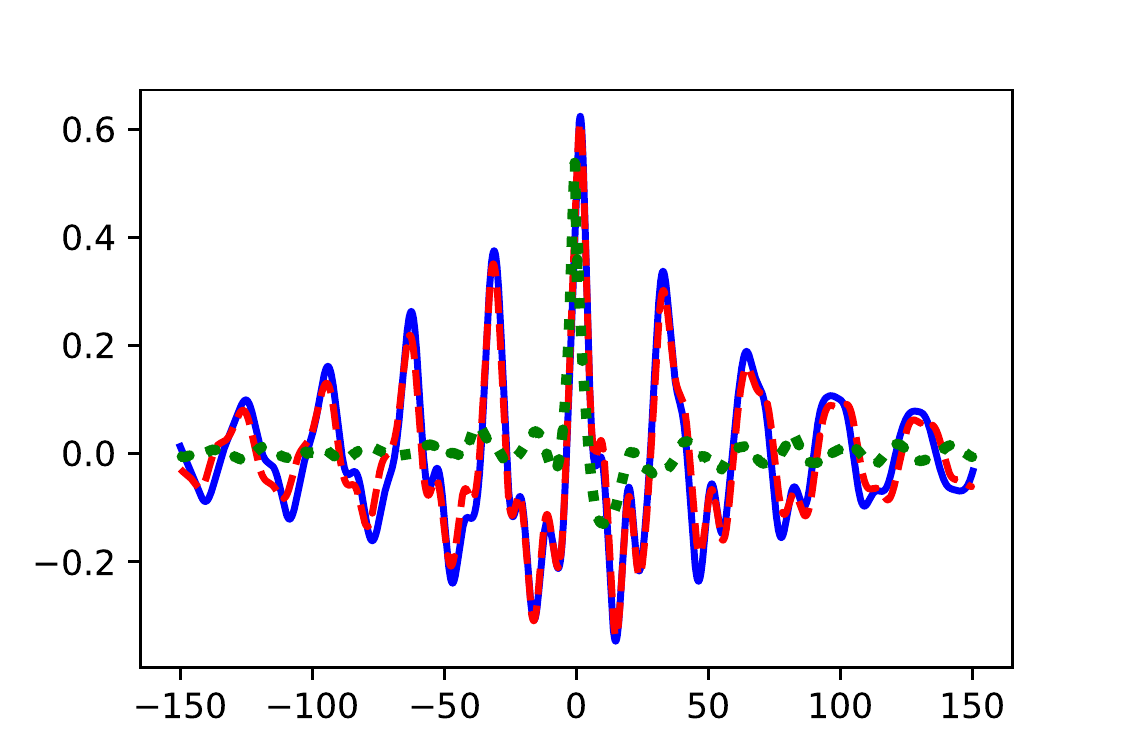}}&
      \resizebox{!}{1.5in}{\includegraphics{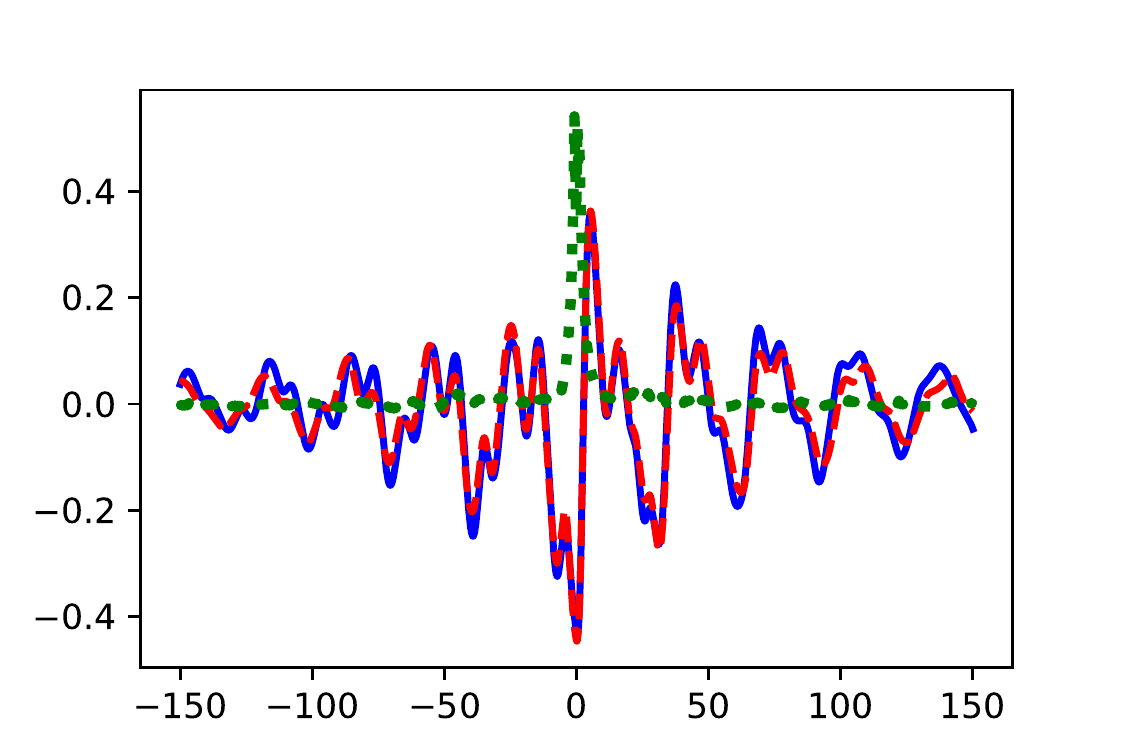}}\\[-1ex]
      &{\footnotesize Time lag}&{\footnotesize Time lag}\\
      \multicolumn{3}{c}{(b) Energy cross-correlations $\cov(|u_k(t)|^2,|u_4(0)|^2)$}\\[1ex]
      \rotatebox{90}{\footnotesize\hspace{0.5in}PDF}&
      \resizebox{!}{1.5in}{\includegraphics{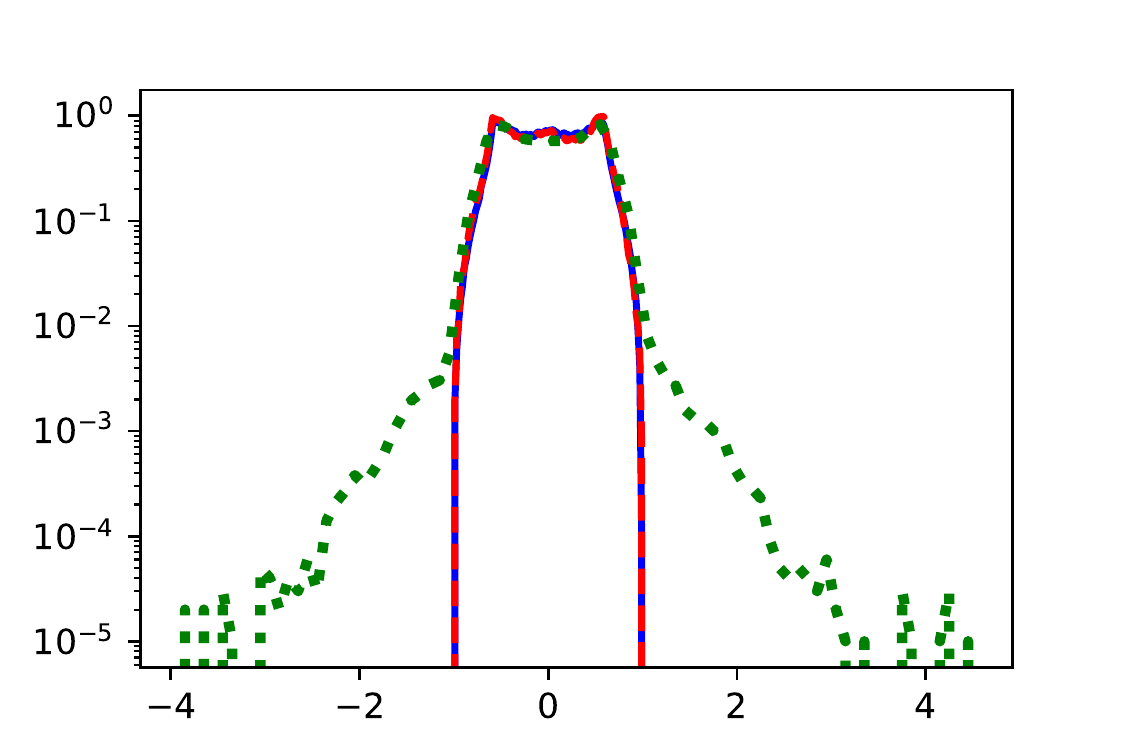}}&
      \resizebox{!}{1.5in}{\includegraphics{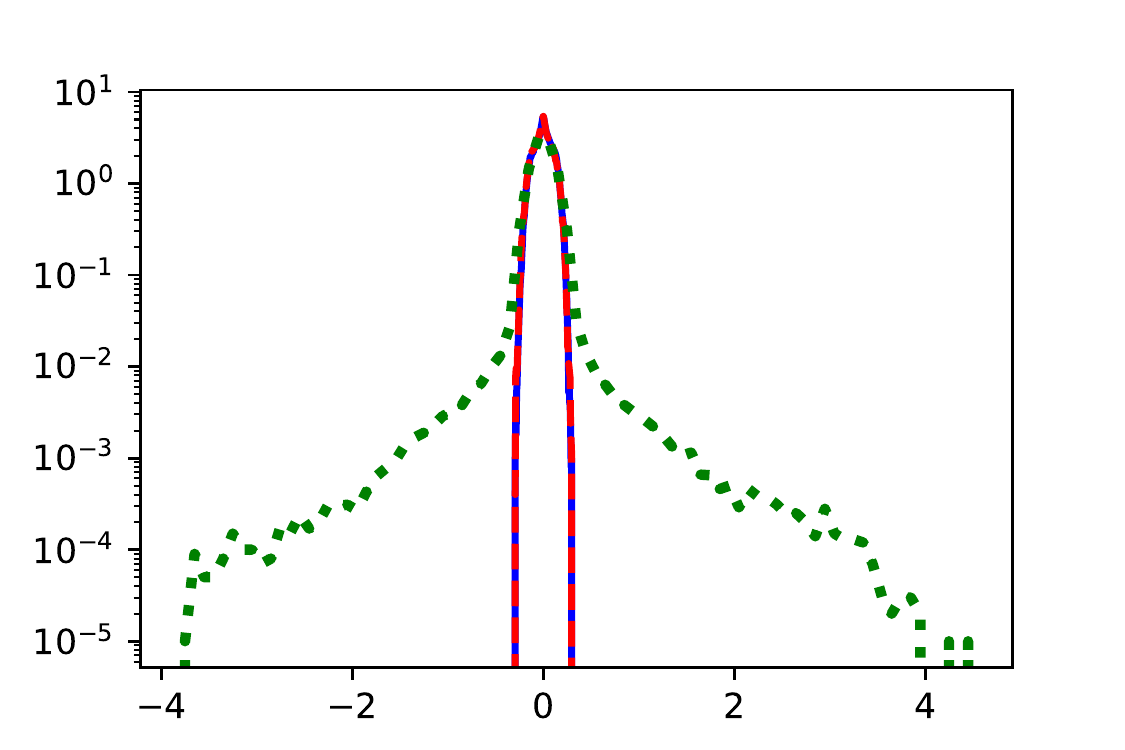}}\\[-1ex]
      &{\footnotesize $Re(u_k)$}&{\footnotesize $Re(u_k)$}\\
      \multicolumn{3}{c}{(c) Marginal distributions}\\
    \end{tabular}~\vspace{-0.2in}
  \end{center}
  \caption{KS statistics.  In all panels, solid blue is the full model,
    dashed red is the reduced model, and dotted green the 5-mode
    truncation.  Panel (a) shows autocovariance functions for two
    Fourier modes $Re(u_k(t))$.  In (b), we show cross correlation
    functions for the energies $|u_k(t)|^2$ and $|u_4(0)|^2$ for
    $k=2,5$.  In (c), distributions of $Re(u_k)$ are shown.}
  \label{fig:ks-stats}  
\end{figure}

In Fig.~\ref{fig:ks-stats}, we examine long-time statistics.  In (a), we
compare the autocovariance functions of selected Fourier modes.  Unlike
the 5-mode truncation, the reduced model is able to reproduce quite
complex features in the ACFs.  Fig.~\ref{fig:ks-stats}(b) shows cross
correlation functions for the energy of the $k$th mode with the energy
of the 4th mode, i.e., $\cov(|u^n_k|^2,|u^0_4|^2)$ as a function of the
time lag $n\dt$; such cross correlation functions can be viewed as a
measure of energy transfer between modes.  The reduced model correctly
predicts these 4th moments, showing that the reduced model captures
genuinely nonlinear effects in KS dynamics.  Panel (c) shows the reduced
model is able to reproduce marginal distributions, whereas the 5-mode
truncation produces marginals that are too wide (compare with
Fig.~\ref{fig:ks-solutions}(c)).
We conclude that both in terms of short-time forecasting and long-time
statistics, the reduced model effectively captures KS dynamics.  These
findings are consistent with~\cite{LLC17}, suggesting the
likelihood-based estimator used in~\cite{CL15,LLC17} and the least
squares estimator above are comparable, and the NARMAX model in
\cite{LLC17} nearly optimal in the least squares sense. Numerical tests
show that slightly different models (with different time lags $p$ and
$r$) may have similar statistical properties (such as consistency) and
comparable performance in prediction. This suggests that there may be
multiple reduced models fitting the data.

\begin{figure}[tb!]
  \begin{center}
    \begin{tabular}{r@{\hskip 0pt}c@{\hskip -6pt}c}
      &Nonlinear, $p=r=1$ & Linear, $p=1,r=0$\\
      \rotatebox{90}{\footnotesize\hspace{0.5in}$Re(u_k(t))$}&
      \resizebox{!}{1.5in}{\includegraphics{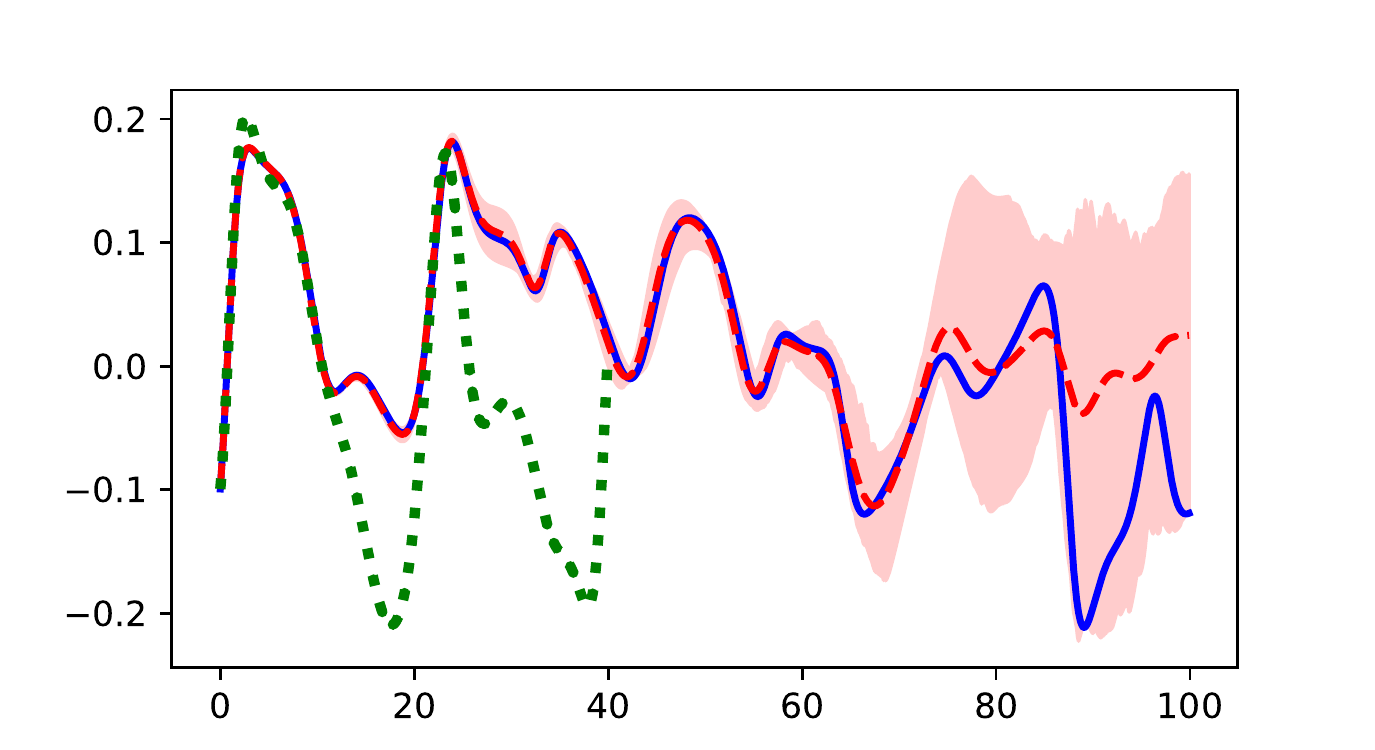}}&
      \resizebox{!}{1.5in}{\includegraphics{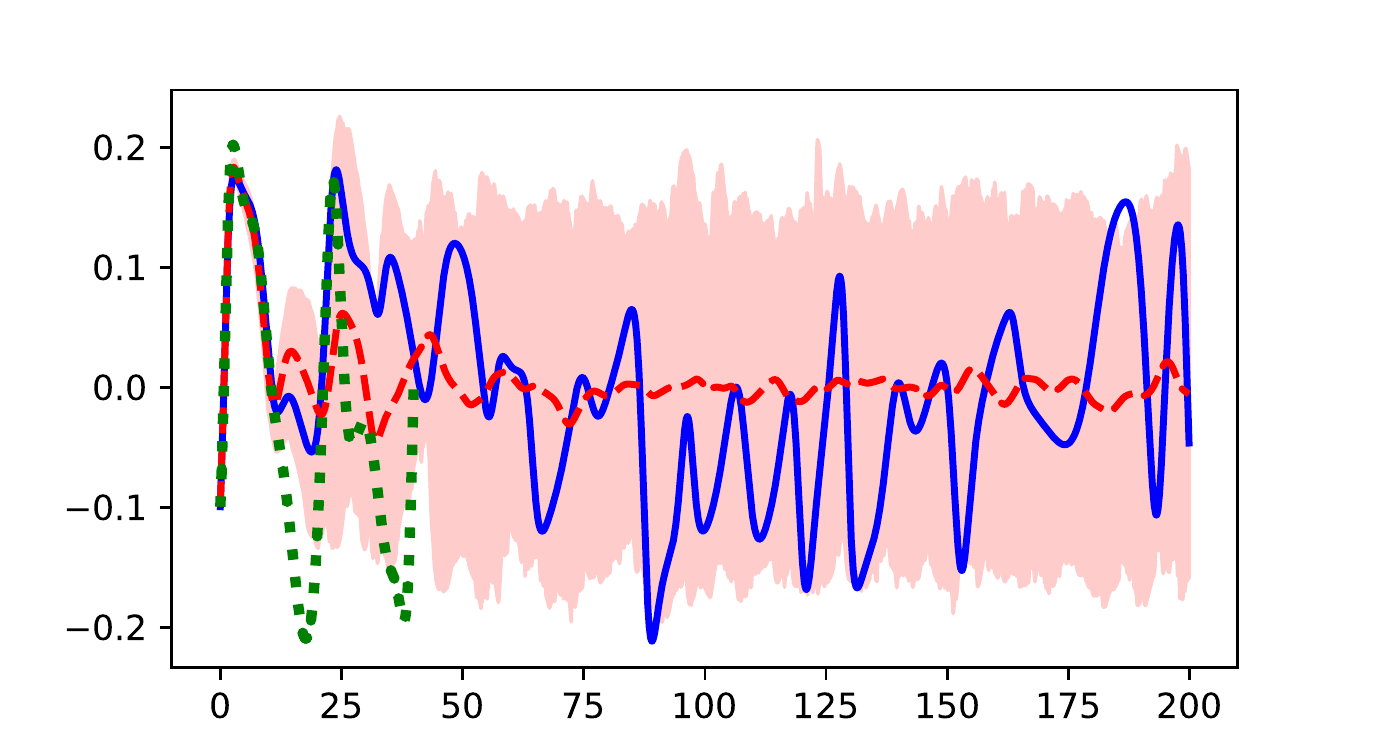}}\\[-1ex]
      &{\footnotesize Time lag}&{\footnotesize Time lag}\\
      \multicolumn{3}{c}{(a) Forecasting}\\[1ex]
      \rotatebox{90}{\footnotesize\hspace{0.5in}ACF}&
      \resizebox{!}{1.5in}{\includegraphics{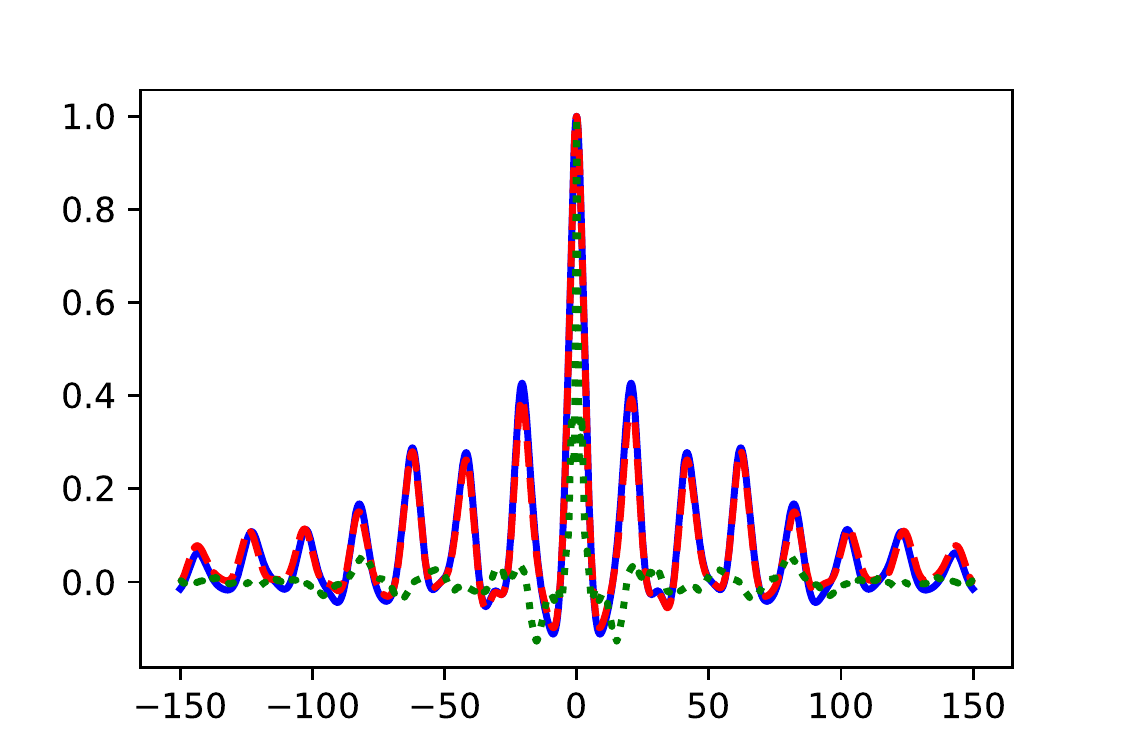}}&
      \resizebox{!}{1.5in}{\includegraphics{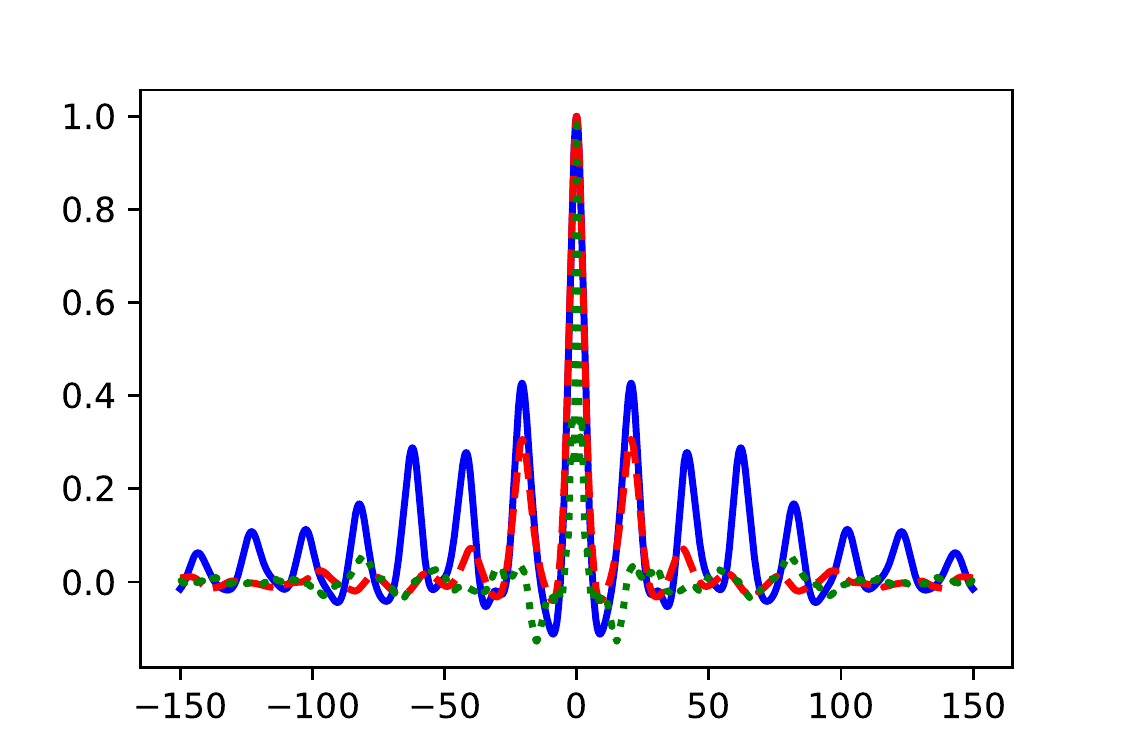}}\\[-1ex]
      &{\footnotesize Time lag}&{\footnotesize Time lag}\\
      \multicolumn{3}{c}{(b) Autocovariance functions $Re(u_k(t))$}\\[1ex]
      \rotatebox{90}{\footnotesize\hspace{0.5in}CCF}&
      \resizebox{!}{1.5in}{\includegraphics{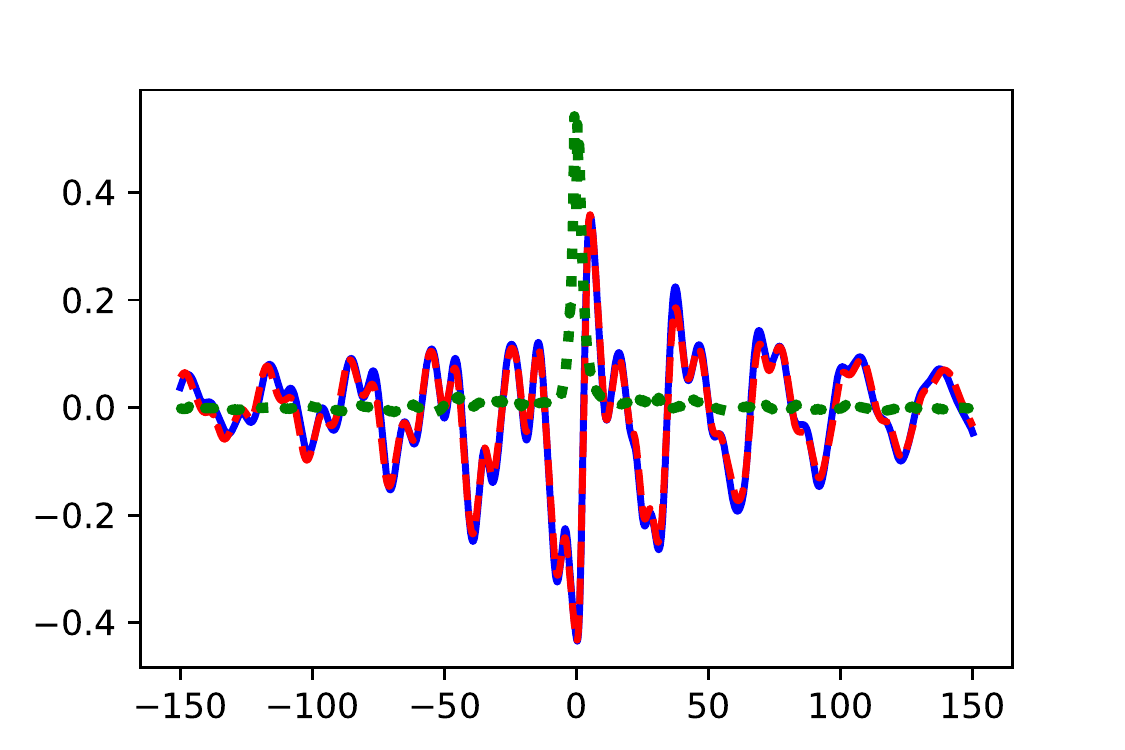}}&
      \resizebox{!}{1.5in}{\includegraphics{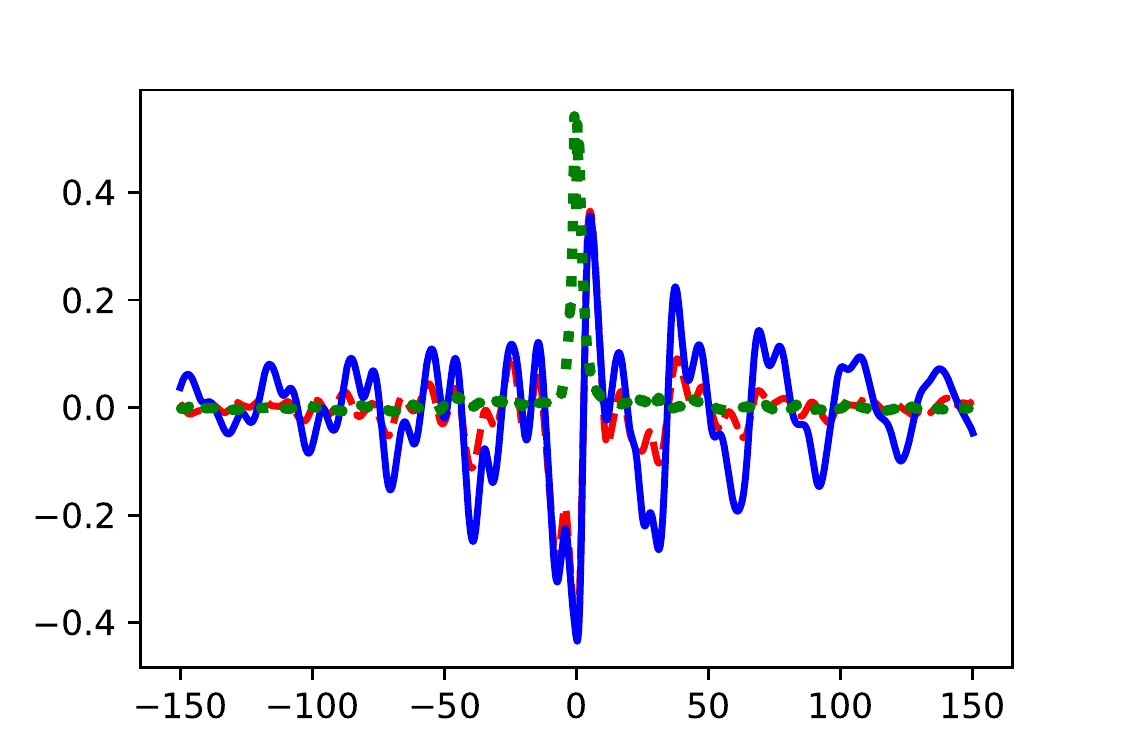}}\\[-1ex]
      &{\footnotesize Time lag}&{\footnotesize Time lag}\\
      \multicolumn{3}{c}{(c) Energy cross-correlations $\cov(|u_k(t)|^2,|u_4(0)|^2)$}\\[1ex]
      \rotatebox{90}{\footnotesize\hspace{0.5in}PDF}&
      \resizebox{!}{1.5in}{\includegraphics{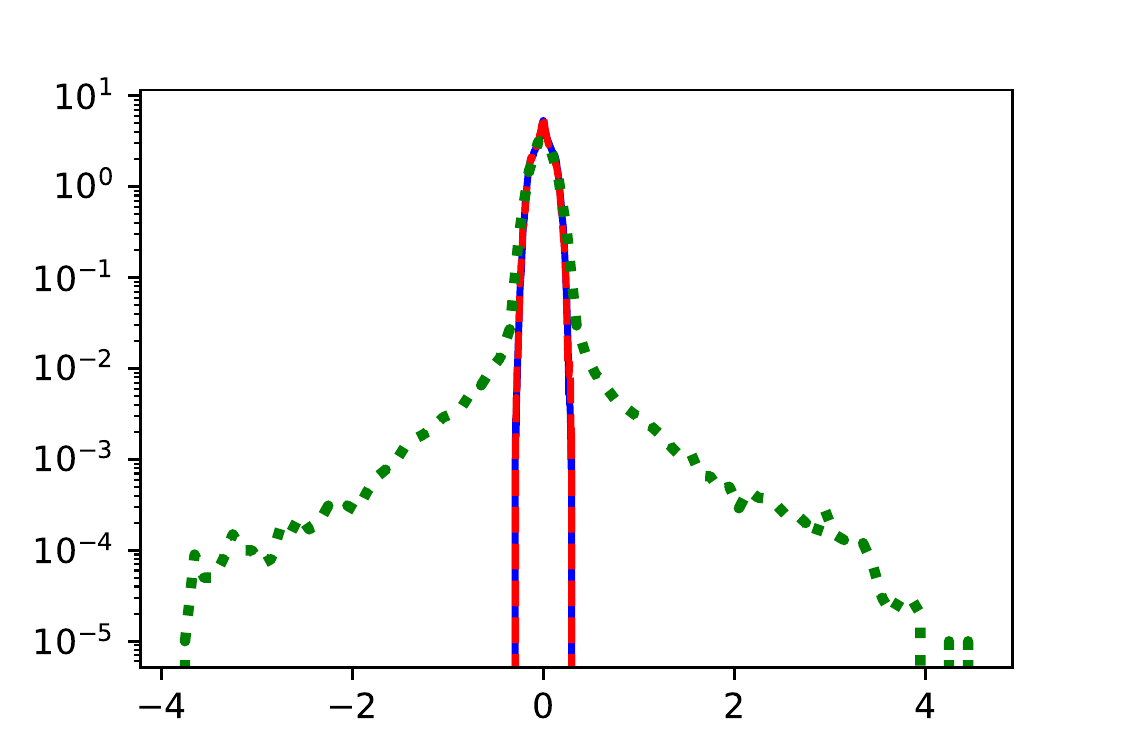}}&
      \resizebox{!}{1.5in}{\includegraphics{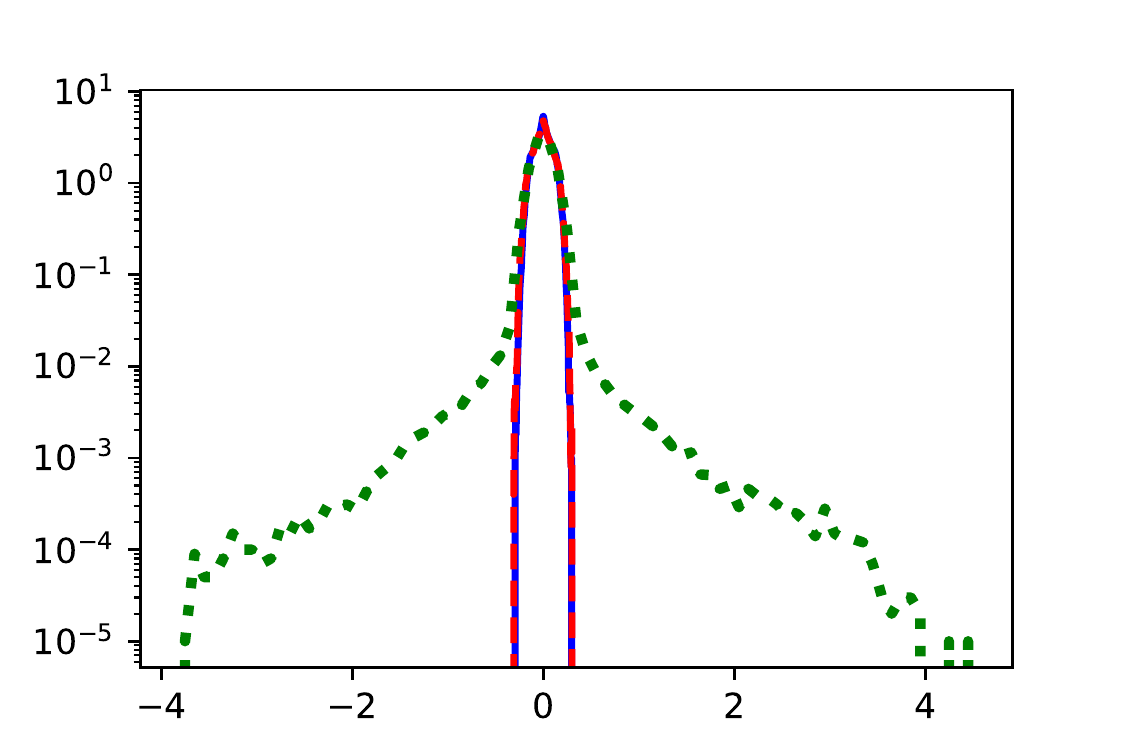}}\\[-1ex]
      &{\footnotesize $Re(u_k)$}&{\footnotesize $Re(u_k)$}\\
      \multicolumn{3}{c}{(d) Marginal distributions}\\
    \end{tabular}~\vspace{-0.2in}
  \end{center}
  \caption{Linear vs nonlinear regression.  Left: results from nonlinear
    regression with $p=r=1$.  Right: results from linear regression with
    $p=1, r=0$.  Note as explained in the text, we did not find any
    orders $(p,r)$ for which both procedures produced useful models.}
  \label{fig:ks-nonlin-comparison}  
\end{figure}

\heading{Linear vs. nonlinear regression.}  Sect.~\ref{sect:rat-narmax}
emphasizes that choice of loss function should be viewed as part of the
model reduction procedure.  In particular, for our {\it ansatz}, the MZ
formalism suggests a least squares approach leading to nonlinear
regression (the nonlinearity arising from the way we parametrize the
transfer function $H(z)$ by a rational approximation).  An alternative
is to infer the coefficients by linear regression, by minimizing the one
step predictions in Eq.~(\ref{eq:narma-one-step-prediction}).  Though
the resulting reduced model Eq.~(\ref{eq:multistep}) is formally
equivalent to Eq.~(\ref{eq:recursion}), the coefficients and the
statistics of the residuals are different.  We emphasize that both
models are nonlinear and share the same functional form, and differ only
in how model coefficients are inferred.  For both models, the residuals
are computed via Eq.~(\ref{eq:cascade}) and a stationary Gaussian
process fitted using the procedure outlined in Sect.~\ref{sect:noise}.

Overall, using linear regression, we found far fewer combinations of
$(p,r)$ for which the reduced models is stable.  Unfortunately, for the
range of relatively low order models we tested ($0\leq p,r\leq 3$), we
did not find any combinations of $p$ and $r$ for which both procedures
produced stable models.  Thus, we did not conduct a direct comparison
between the two.  The closest pair of parameters we found were $p=r=1$
using nonlinear regression, and $p=1,r=0$ using linear regression.  This
means our nonlinear regression example uses a reduced model of the form
$x_{n+2}+a_0x_{n+1} = \Psi(x_{n+1})\cdot b_1 + \Psi(x_{n})\cdot b_0 +
\bar\eta_n$, while our linear regression model has the form
$x_{n+2}+a_0x_{n+1} = \Psi(x_{n})\cdot b_0 + \bar\eta_n$.

Fig.~\ref{fig:ks-nonlin-comparison} shows the results.  Though both
models are fairly low order, the nonlinear regression model has
performance comparable to the higher-order ($p=r=3$) model discussed
above.  In contrast, the linear regression model has significantly worse
forecasting performance, and was unable to reproduce the
auto-correlation or cross correlation functions accurately.  However, it
does reproduce marginal distributions and energy spectra (not shown)
reasonably well.

\begin{figure}[tb!]
  \begin{center}
    \begin{tabular}{c}
      \resizebox{!}{2in}{\includegraphics{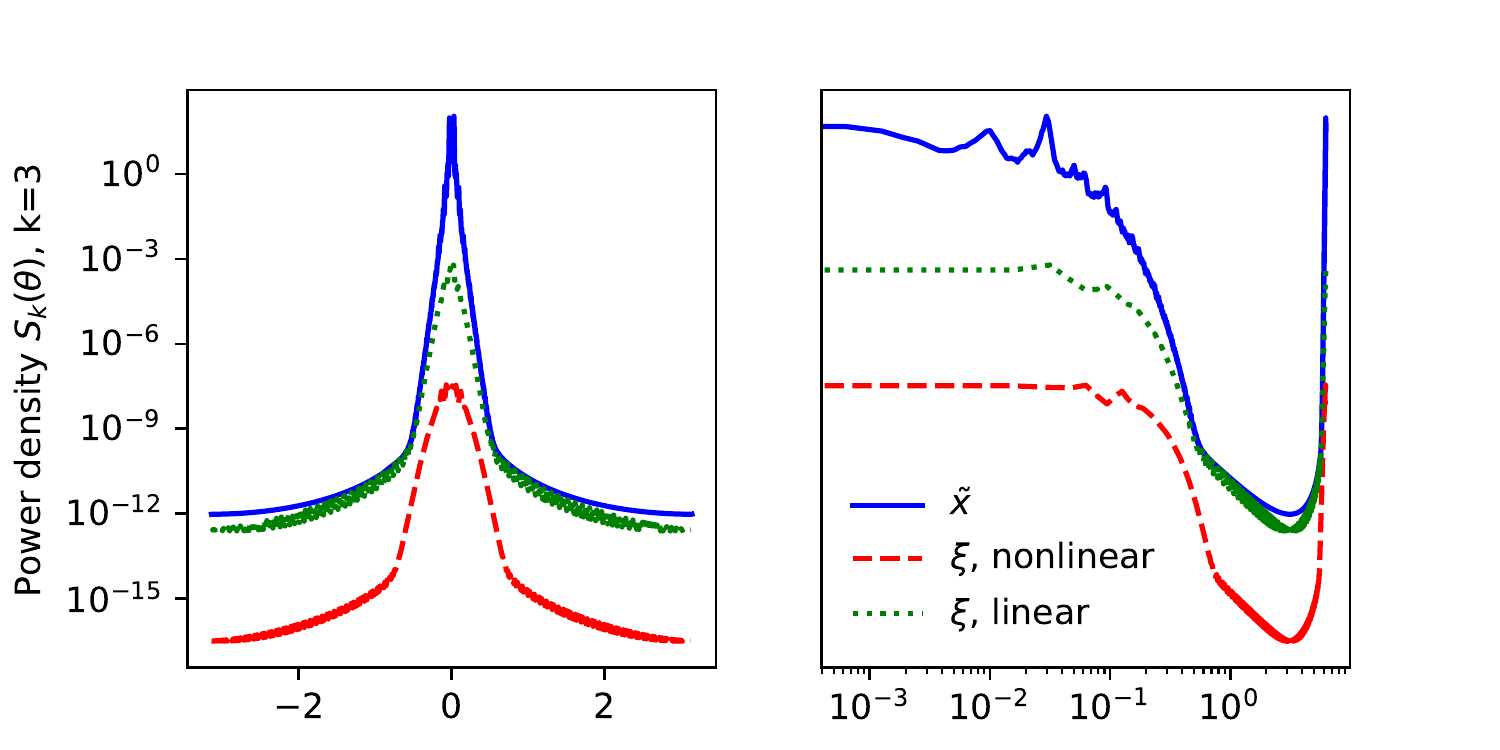}}\\
      $\theta$\hspace{1.6in}$\theta$\\
    \end{tabular}
  \end{center}
  \caption{KS power spectra.  The left panel shows the spectral power
    density $S_{xx}(\theta)$ for $x=u_3,$ the $k=3$ Fourier mode of the
    KS equation.  The right panel shows the same power spectrum on a
    log-log scale to better exhibit the structure near $\theta=0$.  The
    solid blue curve is the power spectrum of the Fourier mode $u_3$
    from the full model, the dashed red curve is the power spectrum of
    the residuals $\xi$ resulting from nonlinear regression with
    $p=r=1$; the dotted green curve is the power spectrum of the
    residual resulting from linear regression with $p=1,r=0$.  Modes
    with $k\neq3$ behave similarly and are not shown.}
  \label{fig:ks-powerspectra}  
\end{figure}

To compare the statistical properties of the reduced models produced by
linear and nonlinear regression, Fig.~\ref{fig:ks-powerspectra} shows the
power spectra $S_{xx}(\theta)$ and $S_{\xi\xi}(\theta)$ for the relevant
variables $x_n$ and the residuals $\xi_n$, for the linear regression
model and the two nonlinear regression model, for the $k=3$ Fourier
mode.  (The other modes show similar trends.)  As far as these power
spectra are concerned, the two nonlinear regression models have nearly
identical behavior.  For nonlinear regression, the residuals $(\xi_n)$
have broader and flatter power spectra than that of $(x_n)$, indicating
that the effect of the approximate Wiener projection here is to capture
the relatively slower dynamics.  The residual is, however, far from
white, suggesting the need for more refined noise models than white
noise forcing.

In contrast, linear regression produces much larger residuals, with a
flat but less broad power spectrum.  It appears that linear regression
could not fit the data nearly as well, but the addition of a suitable
noise model was able to correct for some of the defects of the reduced
model, e.g., marginal distributions and energy spectra.  Temporal
statistics appear to be more delicate, however, and the linear
regression model did not faithfully capture the details of
autocovariance functions.

Overall, these results suggest that for the KS equation, linear
regression results in considerably worse performance than nonlinear
regression.  This is consistent with our expectation (see
Sect~\ref{sect:rat-narmax}) that linear regression may have worse
performance because it neglects longer-range correlations in the data.
For ``static'' quantities like energy spectra and marginal
distributions, it appears that the noise model was able to compensate
for this, but unable to generate correct temporal statistics.

\subsection{Stochastically-forced Burgers equation}
\label{sect:burgers}

Now consider a stochastically-forced viscous Burgers equation
\begin{equation}   
  U_t + UU_x = \nu U_{xx} + \zeta
\end{equation}
with $\zeta(t,x)$ white in $t$ and smooth in $x$, and $U(t,x)$
$2\pi$-periodic in $x$.  More precisely, in Fourier variables,
\begin{equation}
  \label{eq:stochastic-burgers}
  \dot{u}_k = -\frac{i\lambda_k}2\sum_\ell u_\ell u_{k-\ell} -
  \nu\lambda_k^2u_k + \sigma_k\dot{w}_k,
\end{equation}
where $\sigma_k=1$ for $|k|\leq4$ and $\sigma_k=0$ for $|k|>4$,
$\dot{w}_k$ is white noise, and $\lambda_k=k$.  In contrast to the KS
equation, which is deterministic and exhibits self-sustained chaos, the
viscous Burgers equation is dissipative: without forcing, all solutions
converge to the steady state $u\equiv0$ as $t\to\infty$.  Stationary
statistics of $u(x,t)$ thus reflect a balance between the forcing
$\zeta$ and dissipation through viscosity.  We note that the stochastic
Burgers equation has the so-called ``one force one solution'' (1F1S)
property~\cite{weinan2000invariant}: for a given realization of
$\zeta_t, t\geq0$, all initial conditions lead asymptotically to the
same (time dependent) solution.  Put another way, modulo transients,
solutions of \keqref{eq:stochastic-burgers} are determined by the
forcing.

In view of the 1F1S property, a natural question is: given a specific
realization of the forcing $\zeta_t$, can a reduced model correctly
predict the response of the system?  To test this, we compare a
fully-resolved, 128-mode truncation of \keqref{eq:stochastic-burgers}
with an under-resolved 9-mode truncation and a 9-mode reduced model
inferred from data.  Throughout, $\nu=0.05$.  (See
\cite{bunder2017resolution} for an alternate view of this problem.)

\heading{Data-driven reduced model.}  To generate data from the full model, we solve
\keqref{eq:stochastic-burgers} using a scheme of the form
\begin{equation}
  \label{eq:stochastic-numerical-scheme}
  u_k^{n+1} = G_k(u^n,\dt) + \sqrt{\dt}~\sigma_k~w_k^n, 
\end{equation}
where $G_k(u,\dt)$ is the result of applying ETDRK4 to the deterministic
part of \keqref{eq:stochastic-burgers}, $u^n_k = u_k(n\dt)$,
$u^n=(u^n_1,\cdots,u^n_K)$, and $w_k^n$ independent $N(0,1)$ random
variables.  Like the standard Euler-Maruyama scheme,
\keqref{eq:stochastic-numerical-scheme} has weak order 1, but is more
stable~\cite{KP99}.  We solve the full system with timestep $\dt =
0.00125$ and observe every 8th step, so the reduced model has timestep
$\delta=0.01$.  Except for minor differences, this has the form of
Eq.~(\ref{eq:stochastic-burgers-scheme}).

To account for the forcing, we modify \keqref{eq:recursion} to obtain
\begin{subequations}
  \label{eq:burgers-reduced-model}
  \begin{align}
    x_{n+1} =& ~y_n + \xi_{n+1}, \\
    y_{n} + a_{p-1}&y_{n-1} \terms a_0y_{n-p}\\
    =& ~\Psi_{n-p+r}\cdot b_r \terms \Psi_{n-p}\cdot b_0 +\label{eq:burgers-nonlinear-terms}\\
    \label{eq:skew-product-ma-term}
    & ~c_q\bar{w}_{n+q} \terms c_0\bar{w}_n.
  \end{align}
\end{subequations}
The $\bar{w}_n$ in the moving average~(\ref{eq:skew-product-ma-term})
are related to the forcing $w^n$ in
\keqref{eq:stochastic-numerical-scheme} by $\bar{w}^n =
(w^{8n}+\cdots+w^{8n+7})/\sqrt8$; this correlates the full model and the
reduced model during fitting.  The independent noise term $\xi_n$ is
inferred from the residuals as before, and permits one to quantify the
uncertainty in response prediction via ensemble forecasting.  As noted
in Sect.~\ref{sect:wiener-rds}, random dynamical systems like
\keqref{eq:stochastic-burgers} are encompassed within MZ theory, and
\keqref{eq:burgers-reduced-model} can be seen as a special case of the
Wiener projection.  As before, the orders $p$ and $r$ are selected by
trial-and-error.

We have also constructed reduced models of the form
(\ref{eq:recursion}), which do not correlate the reduced and full models
through shared forcing.  All else being equal, we found the performance
of \keqref{eq:burgers-reduced-model} to be strictly better in our tests
than \keqref{eq:recursion} because more information is retained.  We
report results obtained using \keqref{eq:burgers-reduced-model} with
$p=r=1$, leading to a $\sim50$-fold reduction in cost.

The exact form of the predictors $\Psi(\cdot)$ are given in
\ref{app:burgers}.  Interested readers are referred to
\cite{Lu20arxiv} for further investigation of this and other parametric
forms, consistency of estimators, and model selection.

\begin{figure}[tb!]
  \begin{center}
    \begin{tabular}{c@{\hskip -0.05in}c}
      \hspace{0.3in}Spacetime view, full model & Snapshots (legend below)\\
      \resizebox{!}{2.5in}{\includegraphics*[bb=0in 0.1in 1.8in 2.7in]{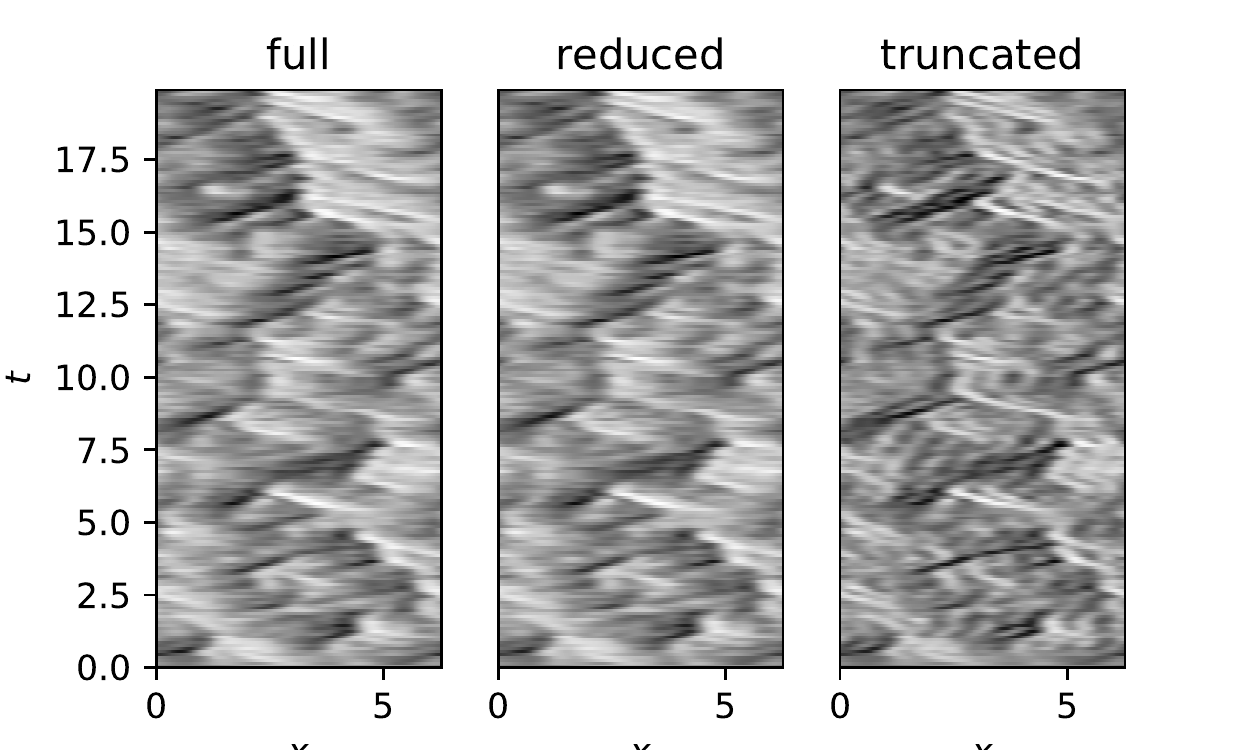}}&
      \resizebox{!}{2.5in}{\includegraphics*[bb=0.15in 0.2in 5in 3.7in]{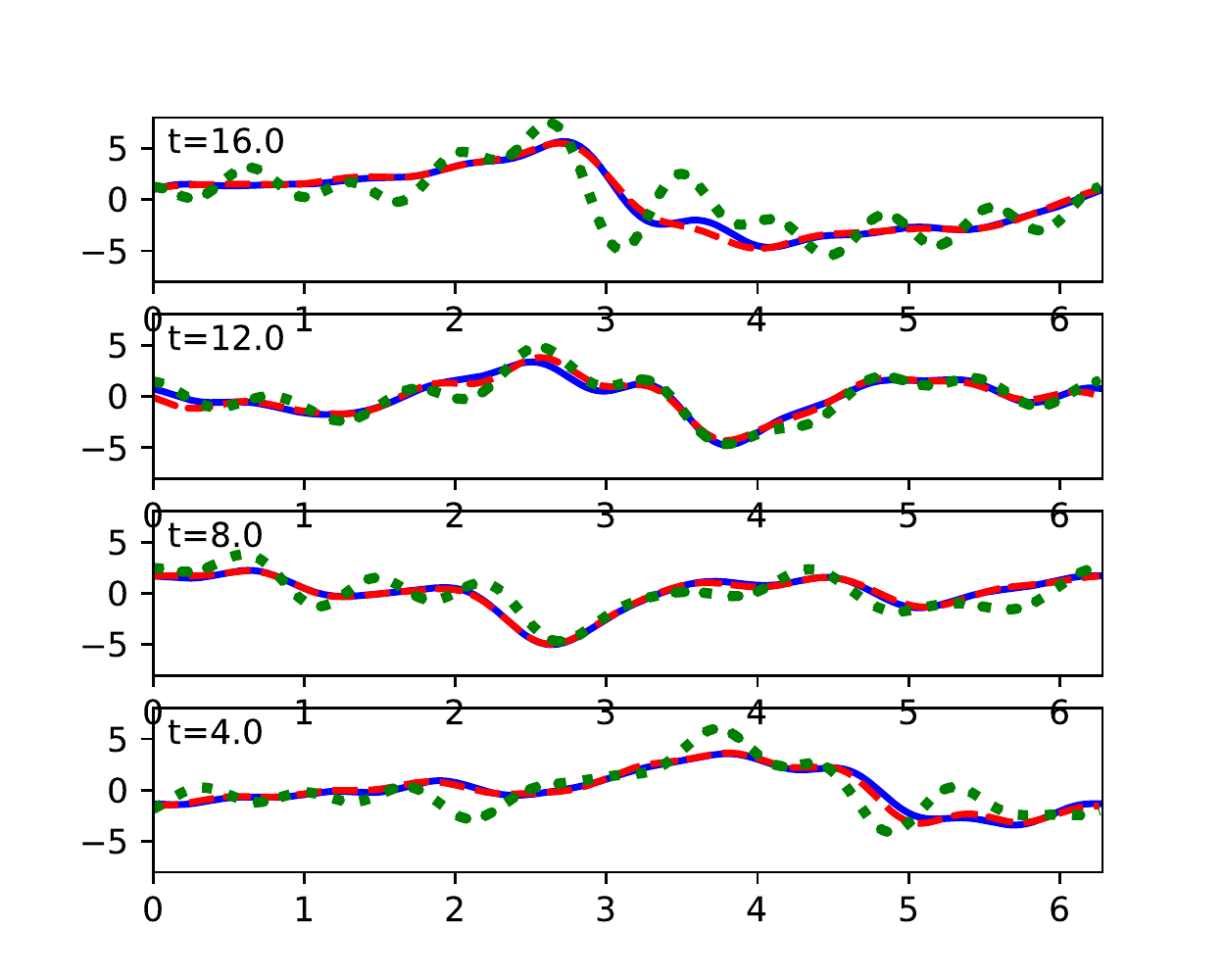}}\\[-2ex]
      \hspace{0.3in}$x$&$x$\\
    \end{tabular}\\
    (a) Burgers solutions\\[2ex]
    \begin{tabular}{c@{\hskip -0.2in}c}
      \begin{tabular}{r@{\hskip 0pt}c}
        \rotatebox{90}{\footnotesize\hspace{0.4in}$k=1$}&
        \resizebox{!}{1.3in}{\includegraphics{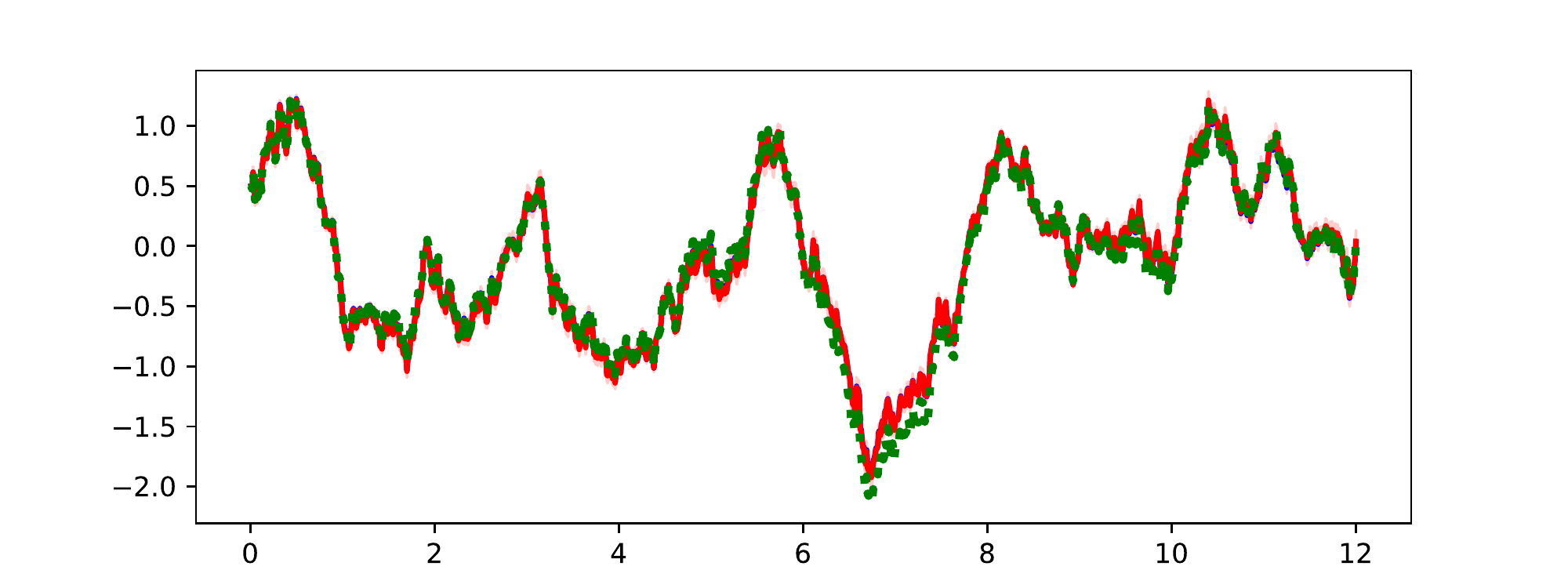}}\\[-1ex]
        \rotatebox{90}{\footnotesize\hspace{0.4in}$k=9$}&
        \resizebox{!}{1.3in}{\includegraphics{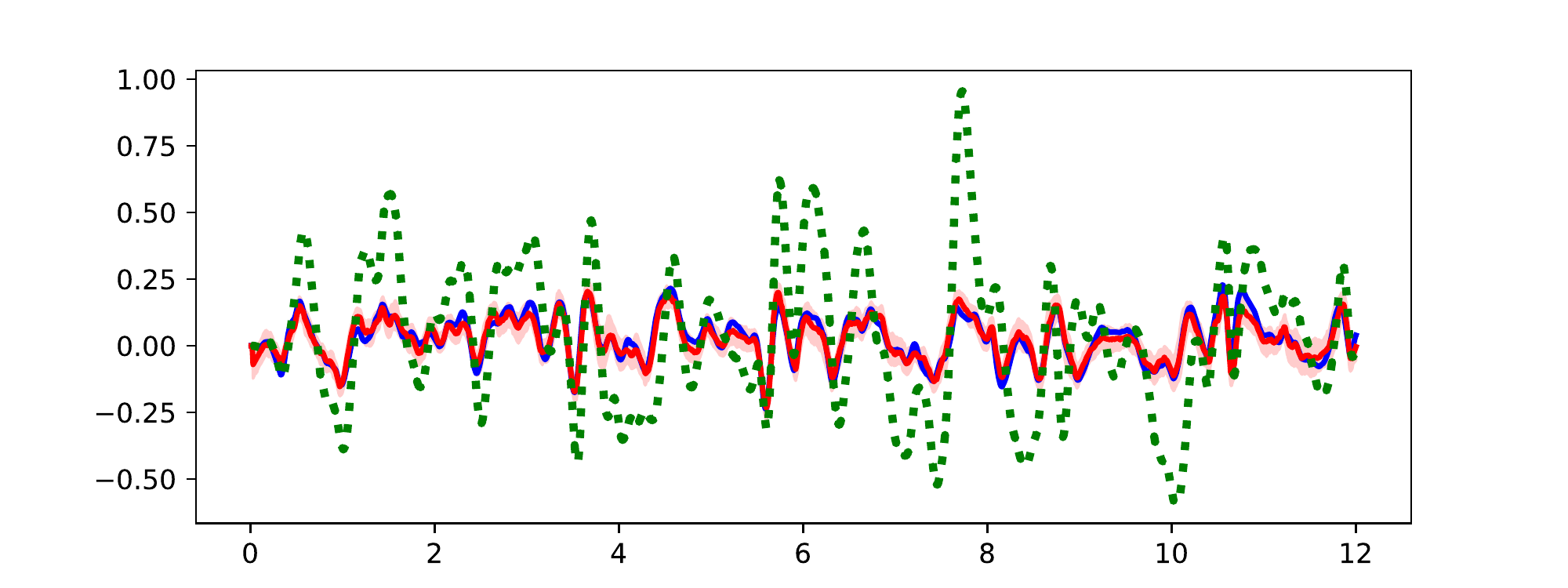}}\\[-1ex]
        &{\footnotesize $t$}\\
      \end{tabular} &

      \begin{tabular}{c}
        \resizebox{!}{0.6in}{\includegraphics*[bb=1.4in 1in 3in 2in]{./legend}}\\
        \resizebox{!}{2in}{\includegraphics{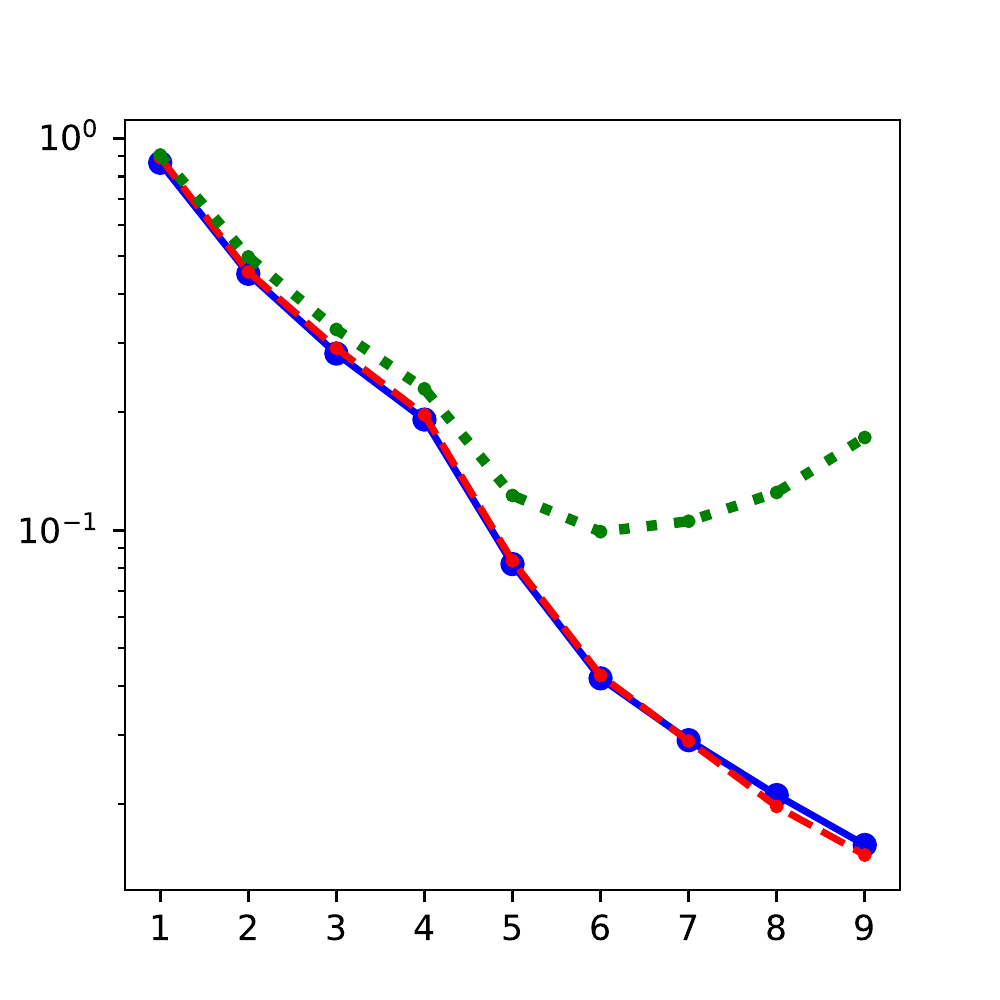}}\\[-2ex]
        {\footnotesize $k$}\\
      \end{tabular}\\
      (b) Trajectories of $Re(u_k(t))$ &  \hspace{6mm}
      (c) Energy $\braket{|u_k|^2}$\\
    \end{tabular}



  \end{center}
  
  \caption{Stochastic Burgers solutions.  Panel (a) shows results
    computed using the 128-mode truncation with $\dt=0.00125$ {\em
      (left)}, and snapshots of the full model, the 9-mode reduced model
    ($\dt=0.01$), and the 9-mode truncation ($\dt=0.00125$).  In (b), we
    plot two Fourier modes as functions of time, with 90\% confidence
    intervals for the reduced model.  Panel (c) shows the energy
    spectrum.}
  \label{fig:burgers-solutions}  
\end{figure}

\heading{Results.}  Fig.~\ref{fig:burgers-solutions}(a) shows sample
solutions.  The 1F1S property suggests that the low modes in the full,
reduced, and the 9-mode truncation models will all be strongly
correlated, as confirmed in the snapshots.  However, one also sees that
the 9-mode truncation exhibits significant deviations from the full
model, unlike the reduced model.  Panels (b) and (c) shows this behavior
in more details: because of the forcing, the low modes of all 3 models
stay close over time, but the 9-mode truncation shows relatively large
deviations from the full model.  As before,
Fig.~\ref{fig:burgers-solutions}(b) shows 90\% confidence intervals for
the reduced model, computed using an ensemble of 100 trajectories.  As
expected, the forced modes are tightly entrained to each other, whereas
the 9-mode truncation shows significant deviation in higher modes.
Because of the 1F1S property, the reduced model can be expected to
correctly forecast the response for as long as information about the
forcing is available.  As for the KS equation, the reduced model here
also reproduces long-time statistics; see
Fig.~\ref{fig:burgers-solutions}(c) for the energy spectrum and
\ref{app:burgers} for other statistics.

Finally, we note that while accurate response forecasting will clearly
become more difficult for larger observation intervals, the reduced
model can nevertheless capture long-time statistics for much larger
observation times.  Indeed, we have tested the reduced model for much
larger observation intervals, up to $0.1$ (see
Appendix~\ref{app:burgers}).

\begin{figure}[tb!]
  \begin{center}
    \begin{tabular}{c}
      \resizebox{!}{2in}{\includegraphics{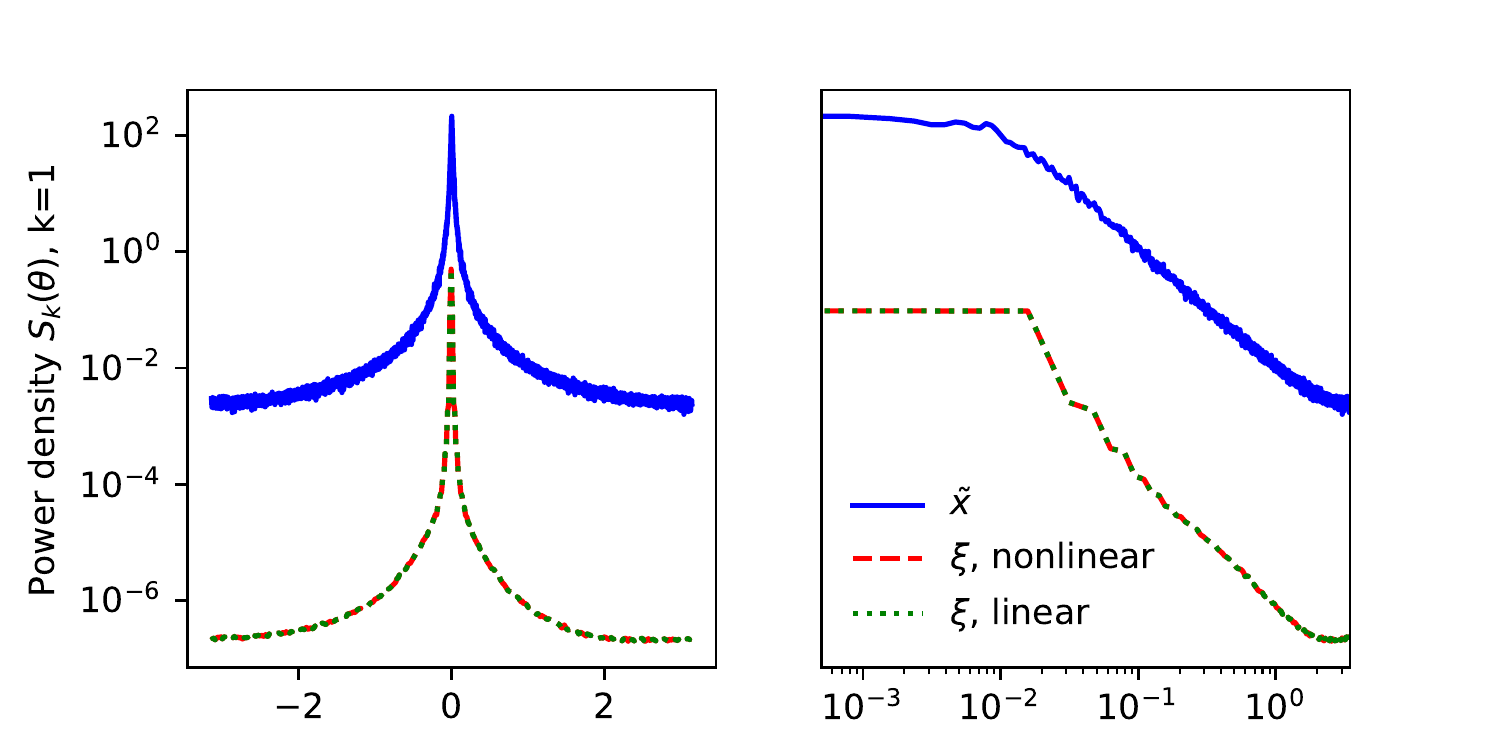}}\\[-1ex]
      $\theta$\hspace{1.6in}$\theta$\\
      \resizebox{!}{2in}{\includegraphics{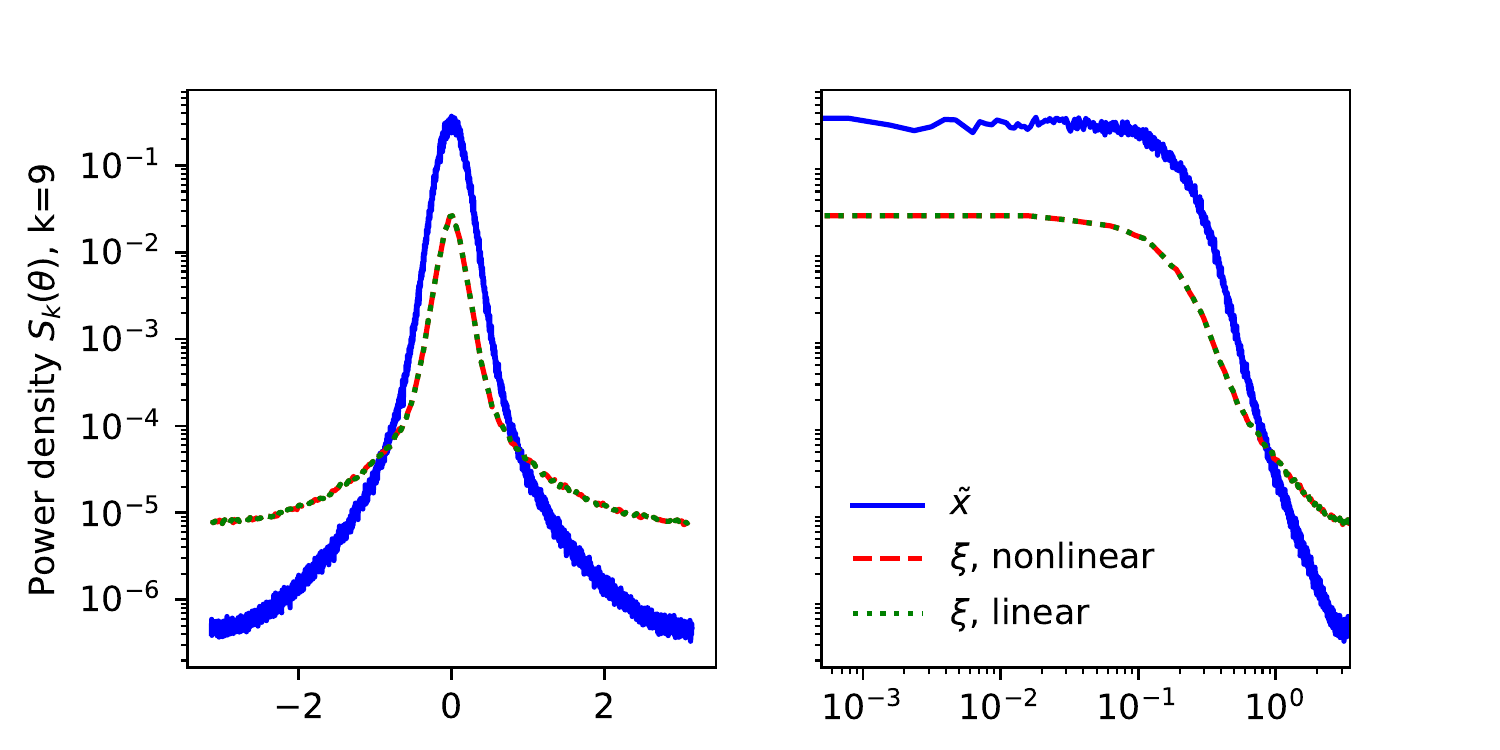}}\\[-1ex]
      $\theta$\hspace{1.6in}$\theta$\\
    \end{tabular}
  \end{center}
  \caption{Burgers power spectra.  Left panels show the spectral power
    density $S_{xx}(\theta)$ for the $k$th Fourier mode of the Burgers
    equation.  Right panels show the same power spectrum on a log-log
    scale.  Solid blue curves are the power spectrum of the Fourier mode
    $u_k$ from the full;model, dashed red curves that of the residuals
    $\xi$ from nonlinear regression; and dotted green curves the
    residual power spectrum from linear regression.  The green and red
    curves essentially coincide.}
  \label{fig:burgers-powerspectra}  
\end{figure}

\heading{Linear vs. nonlinear regression.}  In contrast to the KS
equation, linear and nonlinear regression produced essentially identical
results for the Burgers equation.  In particular, our tests show that
linear regression can produce marginal distributions and ACFs comparable
to the nonlinear regression model, and has nearly identical forecasting
skill; see Fig.~\ref{fig:burgers-nonlin-comparison} in
\ref{app:burgers}.

Fig.~\ref{fig:burgers-powerspectra} compares the spectral power
densities $S_{xx}(\theta)$ and $S_{\xi\xi}(\theta)$ for the relevant
variables $x_n$ and the residuals $\xi_n$, for the $k=3$ Fourier mode.
(The other modes are similar.)  Unlike the KS equation, here the linear
and nonlinear regression give essentially identical power spectra of the
noise.  The residual spectrum is not broader than the spectrum of the
Fourier mode itself, likely because the Fourier modes of the Burgers
equation are subjected to white noise forcing and therefore contain much
higher frequency content than their KS counterparts.  As in the KS
example, leaving out the noise terms entirely led to much worse results.

In view of the discussion in Sect.~\ref{sect:rat-narmax}, the remarkable
contrast between this and the KS equation may be due to the fact that
the Burgers equation is being driven by white noise.  The 1F1S property
implies that the dynamics is largely determined by the forcing, and
hence long-range temporal correlations play less of a role.

\section{Concluding Discussion}

Many issues surrounding this topic remain incompletely understood.  We
mention a few here:
\begin{itemize}

\item \emph{Nonparametric modeling.}  We have focused on parametric
  model reduction in this paper.  But in principle the observation
  functions $\Psi(x)$ can be inferred from data using nonparametric
  methods like delay coordinates, manifold learning, dynamic model
  decomposition, reservoir computing, and other machine learning
  techniques~\cite{freno2019machine,brunton2016discovering,jiang2020modeling,mukhin_principal_2015,berry2013time,ma2018model,pathakModelFreePrediction2018}.
  \keqref{eq:wiener-projection} still applies in these situations, and
  the Wiener projection formulation complements and extends existing
  strategies for data-driven modeling and model reduction by providing a
  systematic guide to incorporating memory and noise effects, in
  situations without sharp scale separation.  For example, one may infer
  $\Psi$ by a combination of delayed coordinates and manifold learning,
  or artificial neural network techniques.

\item \emph{Other rational approximations of $H(z)$.}  The simple
  rational approximation $H(z) = B(z)/A(z)$ is used here out of
  expedience.  Experience has shown that other rational approximations,
  e.g., those based on continued fractions, can sometimes yield
  effective approximations with relatively few undetermined
  parameters~\cite{zwanzig}.  These will be investigated in future work.

\item \emph{Structure-preserving reduced models.}  Most physical systems
  of interest are characterized by exact or approximate conservation
  laws and symmetries, and it is important for reduced models to
  preserve fundamental physical constraints such as these.
  Structure-preserving model reduction is an active area of study, and
  the framework described in this paper may provide a new perspective on
  this problem.

\item \emph{Numerical stability.} In a data-driven approach, one often
  finds that the estimated reduced model is numerically unstable.
  Heuristically, this is because (i) reduced models often coarse-grain
  in both time and space, and the relatively large time steps impose
  more stringent stability requirements; (ii) most loss functions used
  in data-driven model reduction reflect the accuracy of the
  approximation, and one runs the risk of overfitting data.  Indeed, our
  results have shown that the most accurate reduced models (i.e., those
  with the smallest residuals) are not always the best reduced model.  A
  general understanding of numerical stability in these models is
  currently lacking.  Because our nonlinear regression method always
  produces linearly stable models, understanding numerical stability
  will likely require tackling the strong nonlinearities inherent in
  these models.

\item \emph{Quantification of the accuracy of a reduced model.}
  Data-driven approaches have led to many model reduction strategies
  that can successfully reproduce key dynamical features such as energy
  spectrum and correlations.  The development of systematic approaches
  to quantify, analyze, and compare reduced models to full models remain
  incomplete.  It is our hope that the formalism developed in this paper
  will provide a new perspective on this fundamental problem.


\item \emph{Noise modeling.}  For both our examples, the residuals have
  small amplitude, and we have seen that additive noise models work
  relatively well.  We do not know if this approach will continue to be
  effective when the residuals have large amplitude, as occurs in, e.g.,
  molecular dynamics at finite temperatures.

\item \emph{Relationship to other data-driven modeling approaches.}  In
  recent years, a variety of data-driven modeling and model reduction
  techniques have been proposed, applicable in different dynamical
  regimes.  These include delay coordinate
  embedding\cite{broomhead1986extracting,berry2013time,chekroun2020efficient},
  manifold learning and kernel regression
  \cite{duan2014effective,jiang2020modeling}, dynamic mode decomposition
  (DMD)\cite{schmid2010dynamic,tu2014dynamic,kutz2016dynamic}, and many
  others.  The MZ framework should not be viewed as an alternative to
  these methods.  Rather, it is complementary in the sense that it
  provides a general scaffold into which different model reduction
  techniques can fit.  For example, for problems with low-dimensional
  attractors in high-dimensional phase spaces, delayed coordinates and
  extensions like DMD are natural.  But when the underlying assumptions
  (e.g., fast convergence to the low-dimensional attractor,
  deterministic dynamics) are only satisfied approximately, the MZ
  formalism may be useful for suggesting corrections.

\end{itemize}

\medskip

To conclude, we have shown the Wiener projection provides a framework
for data-driven modeling that is grounded in dynamical systems theory.
As such, we view it as a step towards bridging the gap between nonlinear
dynamics theory and data-driven model reduction.  Within this framework,
we give a heuristic derivation of a version of the NARMAX model widely
used in time series modeling and analysis, providing an interpretation
of NARMAX in terms of an underlying dynamical system and evidence that
it may be nearly optimal in the sense of least squares.  In addition to
giving a dynamical systems interpretation for NARMAX, this framework may provide a
starting point for systematic data-driven model reduction. 
Using the KS and stochastic Burgers equations, we have demonstrated the
flexibility and effectiveness of this view of model reduction for
deterministic chaotic and random dynamics.

\bigskip
\noindent
\textbf{Acknowledgments.}  We thank the Mathematics Group at Lawrence
Berkeley National Lab for its support of this work; to Alexandre Chorin
for many useful comments on the manuscript; and to Xiantao Li for
encouraging us to study the discrete Mori-Zwanzig formalism.  KL was
supported by NSF grant DMS-1821286, and FL was supported by NSF grant
DMS-1821211.

\bibliographystyle{elsarticle-num}
\bibliography{ref_Fei18_6,Model_reduction,KL,more}

\appendix

\section{The dual equation and Mori-Zwanzig closure}
\label{app:dual}

In this section, we give an alternate derivation of the MZ
equation~(\ref{eq:mz}) that makes use of a dual equation describing the
evolution of conditional probability distributions.  Though longer, it
gives some additional insights into the meaning of the MZ equation.

As before, let $F$ be a dynamical system with state space $\spaceX$.
Suppose an initial condition $X_0$ is drawn from the distribution
$\rho_0$.  Let $\rho_n$ denote the distribution of $X_n$; then
$\rho_{n+1} = L\rho_n$, where $L$ is the {\it transfer operator},
defined by
\begin{equation}
  \int\varphi\circ F~d\rho = \int\varphi~d(L\rho)
\end{equation}
for suitable test functions $\varphi$.  The above is equivalent to
\begin{equation}
  \int (M\varphi)~d\rho = \int\varphi~d(L\rho),
\end{equation}
i.e., the operator $L$ is the adjoint of the Koopman operator $M$, where
the adjoint of an operator $T$ acting on functions is the operator
$\ad{T}$ acting on distributions defined by $\int (T\varphi)~d\rho =
\int\varphi~d(\ad{T}\rho)$.  With this, and with $P$ and $Q$ as before,
we have
\begin{align*}
  \ad{Q}\rho_{n+1}
  &= \ad{Q} L\rho_n\\
  &= \ad{Q} L(\ad{P}+\ad{Q})\rho_n\\
  &= \ad{Q} L\ad{Q}\rho_n + \ad{Q} L\ad{P}\rho_n,
\end{align*}
using $\ad{P}+\ad{Q}=I$.  Solving the recurrence for $\ad{Q}\rho_n$
gives
\begin{equation}
  \label{ortho}
  \ad{Q}\rho_{n} = (\ad{Q} L)^n\ad{Q}\rho_0 + \sum_{k=0}^{n-1}(\ad{Q} L)^{n-k}\ad{P}\rho_k~.
\end{equation}
From this it follows that
\begin{align*}
  \rho_{n+1} &=  L\rho_n\\
  &=  L\ad{P}\rho_n +  L\ad{Q}\rho_n\\
  &=  L\ad{P}\rho_n +  L(\ad{Q} L)^n\ad{Q}\rho_0 +  L\sum_{k=0}^{n-1}(\ad{Q} L)^{n-k}\ad{P}\rho_k~;
\end{align*}
in the last line we just substituted Eq.~(\ref{ortho}).

The above is equivalent to the operator equation
\begin{equation}
   L^{n+1} =  L\ad{P} L^n
 +  L(\ad{Q} L)^n\ad{Q}
 +  L\sum_{k=0}^{n-1}(\ad{Q} L)^{n-k}\ad{P} L^k~.
\end{equation}
Taking adjoints, we get the Dyson formula
\begin{equation}
  \label{operator-mz}
  M^{n+1} = M^nPM + Q(MQ)^nM + \sum_{k=0}^{n-1}M^kP(MQ)^{n-k}M.
\end{equation}
From this, the MZ equation follows as before.

Suppose now $P$ is conditional expectation with respect to $\mu$.
Observe that for an observable $\varphi$ and a probability distribution
$\rho$, we have
\begin{align*}
  \int P\varphi~d\rho
  &= \iint\left[\int\varphi(x,y)~\mu_{Y|X}(dy|x)\right]~\rho(dx,dy')\\
  &= \int\left[\int\varphi(x,y)~\mu_{Y|X}(dy|x)\right]~\rho_X(dx)\\
  &= \int\varphi(x,y)\int~\mu_{Y|X}(dy|x)~\rho_X(dx)\\
  &= \int \varphi~d(\ad{P}\rho).
\end{align*}
So the dual $\ad{P}$ to the conditional expectation $P$ is
\begin{equation}
  (\ad{P}\rho)(dx,dy) = \rho_X(dx)\cdot\mu_{Y|X}(dy|x)~.
\end{equation}
That is, for a density $\rho~,$ $\ad{P}\rho$ is the product of the
$X$-marginal of $\rho$ and the conditional density $\mu_{Y|X}~.$ The
operator $\ad{P}$ preserves the $X$-marginals of densities, and
$\ad{P}\mu=\mu~.$ If one were to construct reduced models by keeping
only the Markov term in the MZ equation, this corresponds to the closure
assumption that the unresolved modes have statistics given by the
stationary distribution $\mu$ conditioned on the resolved modes.  This
is the discrete-time analog of the averaging principles for ODEs (see,
e.g., \cite{freidlin-wentzell}).

\section{Brief summary of $z$-transform and Wiener filters}
\label{app:wiener-and-z}

For the convenience of readers, this Appendix provides a brief
non-technical summary of some basic facts about $z$-transforms and
Wiener filtering.  See, e.g.,
\cite{wiener,kailath,yaglom1962,yaglom1987} for more details.

\heading{$z$-transforms and linear filtering.}  We first consider (real
or complex, scalar or vector) bi-infinite sequences $\cdots, x_{-1},
x_0, x_1, \cdots$ that are {\em causal}, i.e., $x_n=0$ for $n<0$.  For
simplicity, we assume $(x_n)\in\ell^1$ (though much of what we say below
holds as long as the $x_n$ decay sufficiently fast as $n\to\infty$).
For a causal sequence, its {\it $z$-transform} is the formal series
\begin{equation}
  \label{eq:z-transform-def}
  X(z) = \sum_{n\geq0}x_nz^{-n}.
\end{equation}
In the above expression, $z$ should be viewed as a complex variable,
though the series typically does not converge for all $z\in\C$.  The
$\ell^1$ assumption (which covers many examples in applications) means
the domain of convergence of $X(z)$ includes the unit circle and
$X(e^{-i\theta})$ is a Fourier series with $x_n$ as coefficients.  In
this case, the $z$-transform is invertible by
\begin{equation}
  \label{eq:cauchy}
  x_n = \frac1{2\pi}\int_0^{2\pi}e^{-in\theta}X(e^{-i\theta})~d\theta.
\end{equation}
More generally, the $z$ transform can be inverted by an appropriate
application of the Cauchy integral formula.

The $z$-transform is the analog of the Laplace transform for difference
equations.  Two key properties include:
\begin{enumerate}
\item {\em Shifts:} if $y_n = x_{n+1}$ for $n\geq0$, then $Y(z) =
  z(X(z)-x_0)$.
\item {\em Convolution:} if $w_n = (x\star y)_n =
  \sum_{k\geq0}x_ky_{n-k}$, then $W(z) = X(z)\cdot Y(z)$.
\end{enumerate}

In signal processing and time series analysis, the $z$-transform is
useful for representing the action of ``linear filters.''  That is,
suppose we have a signal $(x_n)$.  A linear filter is a linear
transformation mapping $(x_n)$ to $(y_n)$, with
\begin{equation}
  y_n = (x\star h)_n = \sum_{k\geq0}x_k\cdot h_{n-k}.
\end{equation}
The sequence $(h_n)$, which defines the filter, is known as its {\em
  impulse response}, so called because $h_n$ is the response of the
filter when $(x_n)$ is the unit impulse, i.e., $x_n=\delta_{n0}$,
$\delta_{mn}$ being the Kronecker delta function.  By the convolution
property, we then have $Y(z) = H(z)X(z)$.  $H(z)$ is the ``transfer
function'' of the linear filter.

One of the ways in which the $z$-transform is useful is that the
analytic properties of the transfer function encode the asymptotic
behavior of the impulse response.  For example, if the transfer function
$H(z)$ of a scalar filter were meromorphic and all its poles lie
strictly inside the unit disc, then Eq.~(\ref{eq:cauchy}) tells us $h_n$
is causal and decays exponentially as $n\to\infty$.  (If we only know
that the restriction of $H$ to the unit circle is continuous, then
$h_n\to0$ is implied by the Riemann-Lebesgue lemma.)  In the reverse
direction, if $(h_n)\in\ell^1$ (as we assume), then $H(z)$ cannot have
any poles outside the unit disc.

\heading{An application to NARMAX.}  In Sect.~\ref{sect:rat-narmax}, we
asserted the equivalence of Eqs.~(\ref{eq:recursion}) and
(\ref{eq:multistep}) modulo transients.  Here we show a derivation using
$z$-transforms; an alternative is to use the substitution
$y_n=x_{n+1}-\xi_{n+1}$ in Eq.~(\ref{eq:multistep}).  One of the
advantages of the $z$-transform method is that it provides an
operational calculus for keeping track of indices systematically.

First, take $z$-transforms of Eq.~(\ref{eq:recursion}), we get
\begin{subequations}
  \begin{align}
    \label{eq:recursion-z-a}
    z(X(z) - x_0) &= Y(z) + z(\Xi(z)-\xi_0) \\
    \label{eq:recursion-z-b}
    A(z)Y(z) + p_0(z) &= \Psi(z)\cdot B(z) + q_0(z)
  \end{align}
\end{subequations}
where $p_0(z)$ and $q_0(z)$ are polynomials whose coefficients are
functions of the initial conditions $x_0,\cdots,x_p$ and
$\Psi(x_0),\cdots,\Psi(x_q)$, with $\deg(p)\leq\deg(A)$ and
$\deg(q)\leq\deg(B)$, and $p_0\equiv q_0\equiv0$ if the
$x_0=\cdots=x_p=\Psi_0=\cdots=\Psi_q=0$.  Substituting
Eq.~(\ref{eq:recursion-z-b}) into (\ref{eq:recursion-z-a}) and
simplifying, we get
\begin{equation}
  \label{eq:recursion-z}
  zA(z)(X(z)-x_0) + p_0(z) = \Psi(z)\cdot B(z) + q_0(z) + zA(z)(\Xi(z)
  - \xi_0).
\end{equation}
For comparison, if we transform Eq.~(\ref{eq:multistep}), we get
\begin{equation}
  \label{eq:multistep-z}
  zA(z)X(z) + p_1(z) = \Psi(z)\cdot B(z) + q_1(z) + zA(z)\Xi(z).
\end{equation}
Comparing Eqs.~(\ref{eq:recursion-z}) and (\ref{eq:multistep-z}), we see
they are equivalent modulo terms involving initial conditions.  If all
the zeros of $A(z)$ lie inside the unit circle, then transients will
decay as $n\to\infty$, so modulo transients Eqs.~(\ref{eq:recursion-z})
and (\ref{eq:multistep-z}) are equivalent.  In particular, the
recursions are exactly equivalent for homogeneous initial conditions.

The above argument relies on the $z$-transform.  Because the recursions
are driven by the $\xi_n$, its validity hinges on what we assume about
$\xi_n$: if the $\xi_n$ were, e.g., white noise, then the $z$-transforms
are not well-defined, but if the $\xi_n$ decay sufficiently fast as
$n\to\infty$, then the $z$-transforms are valid.  Supposing now that
there is a sequence $\xi_n$ such that Eqs.~(\ref{eq:recursion-z}) and
(\ref{eq:multistep-z}) are not equivalent for homogeneous initial
conditions $x_0=\cdots=x_p=\Psi_0=\cdots=\Psi_q=0$.  Then there is a
least $N>0$ for which they would disagree.  But then if we set $\xi'_n =
\xi_n$ for $n\leq N+p$ and $\xi'_n =0$ for $n>N+p$, then (because the
recursion has order $p$) the two recursions would differ when driven by
$\xi'_n$.

\heading{Correlation functions and power spectra.}  The preceding
discussion of the $z$-transform only makes sense if the sequences
involved decay sufficiently fast as $n\to\infty$.  In our context, we
are interested in convolving such a sequence $(h_n)$ with stationary
stochastic processes.  The formal series~(\ref{eq:z-transform-def}) does
not make sense.

A standard approach is based on correlation functions.  Suppose $(X_k)$
and $Y_k$ are zero-mean stationary stochastic processes taking values in
$\R^d$.  We define the (matrix-valued) correlation function to be
\begin{equation}
  C_{xy}(k) = E\Big(X_k\cdot Y_0^*\Big)
\end{equation}
where ``$*$'' denotes the conjugate transpose.  The corresponding
\emph{power spectrum} is
\begin{equation}\label{eq:powerSpectrum}
  S_{xy}(z) = \sum_k z^{-k}C_{xy}(k).
\end{equation}
Note this generalizes the notion of spectrum introduced earlier, and we
are abusing notation slightly.  The spectrum introduced earlier is
$S_{xy}(e^{-i\theta})$.  We record some useful properties:
\begin{enumerate}

\item $C_{xx}(0)$ is hermitian positive-semidefinite.

\item $C_{xy}(k)^* = C_{yx}(-k),$ in particular $C_{xx}(k)^* =
  C_{xx}(-k).$

\item $S_{xy}(e^{-i\theta})^* = S_{yx}(e^{-i\theta}).$

\item $S_{xx}(e^{-i\theta})^*=S_{xx}(e^{-i\theta})$, i.e., the power
  spectrum is hermitian for all $\theta$.


\item If $Y=h\star X$, then
  \begin{equation}
    C_{yx}(n) = \sum_kh_{n-k}\cdot C_{xx}(k),
  \end{equation}
  or $C_{yx} = h\star C_{xx}$.

\item Taking $z$-transforms yields
  \begin{equation}
    S_{yx}(z) = H(z)\cdot S_{xx}(z).
  \end{equation}
  Note the above identities are valid for both scalar and matrix
  quantities.

\item Similarly,
  \begin{equation}
    C_{xy}(n) = \sum_kC_{xx}(n+k)\cdot h^*_k.
  \end{equation}
  The $z$-transform is now
  \begin{equation}
    S_{xy}(z) = S_{xx}(z)\cdot H^*(1/z)
  \end{equation}
  where $H^*$ is the $z$-transform of the sequence $h^*_n$.

\item Putting these relations together yields $C_{yy} = h\star
  C_{xx}\star h^*$, or equivalently
  \begin{equation}
    S_{yy}(z) = H(z)\cdot S_{xx}(z)\cdot H^*(1/z).
  \end{equation}
  On the unit circle, this simplifies to
  \begin{equation}
    S_{yy}(e^{-i\theta}) = H(e^{-i\theta})\cdot
    S_{xx}(e^{-i\theta})\cdot H^*(e^{i\theta}).
  \end{equation}
  In the scalar case, this reduces to $S_{yy}(e^{-i\theta}) =
  |H(e^{-i\theta})|^2S_{xx}(e^{-i\theta}).$

\end{enumerate}
These properties also form the basis for the random Fourier
representation of stationary stochastic processes in
Eq.~(\ref{eq:random-fourier}).

\heading{Wiener filtering.}  We now record some basic results of Wiener
filter theory for interested readers.  This material is not used
directly in the paper.

The Wiener filter is the linear filter $(h_n)$ that minimizes the MSE
\begin{equation}
  \E\Big|X_n - \sum_{k\geq0} \Psi_{n-k}\cdot h_{-k}\Big|^2.
\end{equation}
Equivalently, if we write
\begin{equation}
  X_n = \sum_k h_{n-k}\cdot \Psi_k + \xi_n
\end{equation}
this amounts of choosing $(h_n)$ to minimize the residuals
$\E|\xi_n|^2$.  One can show that the power spectrum satisfies
\begin{equation}
  \label{eq11}
  S_{\xi\xi}
  ~=~ \ub{S_{xx} - S_{x\psi}\cdot
  S_{\psi\psi}^{-1}\cdot S_{\psi
    x}}{(I)}
  ~+~\ub{(H\cdot S_{\psi\psi}-S_{x\psi})\cdot
  S_{\psi\psi}^{-1}\cdot(S_{\psi\psi}\cdot
  H^*-S_{\psi x})}{(II)}
\end{equation}
where $S_\cdot(\cdot)$ denotes power spectra as before, and $H(z)$ is
the $z$-transform of $(h_n)$.  Observe $S_{\xi\xi}(e^{-i\theta})\geq0$
for all $H$.  If we set
\begin{equation}
  \label{eq16}
  H(e^{-i\theta}) = S_{x\psi}(e^{-i\theta})\cdot
  S_{\psi\psi}^{-1}(e^{-i\theta})~,
\end{equation}
then (II) vanishes.  This means (I) is $\geq0$.  Since (II) is obviously
$\geq0$ as well, we see $Tr(S_{\xi\xi})$ is minimized by
Eq.~(\ref{eq16}).

Unfortunately, the linear filter $(h_n)$ defined by Eq.~(\ref{eq16}) may
not be {\em cauasal,} i.e., $h_n$ may be nonzero for $n<0$.  Such a
filter would use future values of $\Psi_m$ with $m>n$ to predict $X_n$.
How, then, do we find a causal filter, i.e., one with $h_n=0$ for $n<0$?
Let us first look at the special case where $S_{\psi\psi}(z)\equiv
I_{d\times d}~,$ i.e., $(\Psi_n)$ is ``white.''  Then the functional to
be minimized is
\begin{equation}
  Tr\left((H-S_{x\psi})\cdot(H^*-S_{\psi x})\right).
\end{equation}
By Plancherel's Theorem, the optimal \emph{causal} solution is given by
$H = [S_{x\psi}]_+$, where
\begin{equation}
  [S]_+(e^{-i\theta}) = \frac1{2\pi}\sum_{n=0}^\infty \int_0^{2\pi}
  e^{in(\theta-\theta')}S(e^{-i\theta'})~d\theta'.
\end{equation}
Summing over $n\geq0$ instead of $n\in\Z$ sets the impulse response
$s_n=0$ for $n<0$, thus making it causal.  The $[\cdot]_+$ operator
transforms a given function to the time domain, zero out all entries for
$n<0$, then transform back to frequency domain.

Now, in general $\Psi_n$ will not be white.  But, since
$S_{\psi\psi}(e^{-i\theta})\geq0$, there exist $C$ such that
$C(e^{-i\theta})\cdot C^*(e^{i\theta})=S_{\psi\psi}(e^{-i\theta})$.  So
if we take $W=C^{-1}$ (as a function on $S^1$), then $w\star\Psi$ will
be white.  A remarkable fact is that under very broad conditions, there
is a function $W(z)$ such that all its poles {\em and} zeros lie inside
the unit circle, and $W(e^{-i\theta}) = C(e^{-i\theta})$.  Such a $W$
defines a causal stable linear filter $(w_n)$ such that $w\star\Psi$ has
power spectrum
\begin{equation}
  W(e^{-i\theta})\cdot S_{\psi\psi}(e^{-i\theta})\cdot W(e^{-i\theta})^*
  \equiv I_{d\times d},
\end{equation}
i.e., $w\star\Psi$ is white.  (The filter $(w_n)$ is known as a
whitening filter.)  Using the whitening filter, one can check that
\begin{equation}
  H(z) = [S_{x\psi}(z)\cdot W^*(1/z)]_+\cdot W(z)
\end{equation}
is indeed the causal Wiener filter.

\section{Kuramoto-Sivashinsky equation}
\label{app:ks}

\subsection*{Nonlinear terms in the NARMAX model}

The Kuramoto-Sivashinsky example in
Sect.~\ref{sect:numAll}  uses the reduced
model from \cite{LLC17}.  For the convenience of readers, the full
{\it ansatz} is reproduced here:
\begin{subequations}
  \applabel{eq:ks-ansatz} 
  
  \begin{align}
    u_k^{n+1} =&~ u_k^n + R^{\dt}_k(u^n)~\dt +
    z_k^{n}~\dt\applabel{eq:ks-ansatz-a}, \\[2ex]
    z_k^{n+1} =&~ \Phi_k^n +  \xi_k^{n+1}, \applabel{eq:ks-ansatz-b}\\[2ex]
    \Phi^n_k =&~\sum_{j=0}^p a_{k,j}z_k^{n-j} +
    \sum_{j=0}^rb_{k,j}u^{n-j}_k + c_{k,(K+1)} R_k^{\dt}(u^n)\nonumber\\
    &+i\sum_{j=1}^K c_{k,j}\tilde{u}^n_{j+K} \bar{\tilde{u}^n_{j+K-k}}
    +\sum_{j=0}^q d_{k,j}\xi_k^{n-j},
    \applabel{eq:ks-ansatz-c}
  \end{align}  
  where
  \begin{equation}
    \applabel{eq:ks-ansatz-d}
    \tilde{u}^n_j = 
    \left\{
    \begin{array}{ll}
      u^n_j~, & 1\leq j\leq K\\[2ex]
      i\sum_{\ell=j-K}^K u_\ell^n u_{j-\ell}^n~, & K < j \leq 2K.
    \end{array}
    \right. 
  \end{equation}
\end{subequations}
The nonlinear terms in Eqs.~(\appref{eq:ks-ansatz-c}) and
(\appref{eq:ks-ansatz-d}) are suggested by inertial manifold theory.
See \cite{LLC17} for details.

We compare the above {\it ansatz} to the model used in this study, of the form~\eqref{eq:recursion} 
with predictors in \eqref{eq:Psi_KSE}. 
It is straightforward to show that the {\it ansatz} in
Eq.~(\appref{eq:ks-ansatz}) is equivalent to a model of the form in \meqref{eq:multistep}: 
\begin{equation}
  \applabel{eq:multistep}
  u_{n+p'+1} + a_{p'-1}u_{n+p'}+\cdots+a_0u_{n+1} = \Psi'_{n+q'}\cdot
  b_{q'} + \cdots+\Psi_n\cdot b_0 + \xi'_{n+1},
\end{equation}
for some choice of orders $p',q',$ coefficients $a_i,b_i$, functions $\{\Psi'_n\}$ and noise
$\xi'_n$. 
In addition to the different approaches estimating the parameters $(a,b)$, the models are different in the following aspect: 
\begin{enumerate}
\item Here we model the noise by a Gaussian process using
  power spectrum from the residual $\tilde{\xi}_n$, whereas
  Eq.~(\appref{eq:ks-ansatz}) models the noise by a moving average process.

\item as suggested by the Wiener projection formalism,
the model \eqref{eq:recursion} in this study contains time-delayed copies of all nonlinear
  terms, whereas Eq.~(\appref{eq:ks-ansatz}) does not.

\end{enumerate}

\subsection*{Detailed Numerical results}

\begin{figure}
  \begin{center}
    \begin{tabular}{r@{\hskip 0pt}c@{\hskip -6pt}c}
      &{\small Trajectories} & {\small Marginals}\\
      \rotatebox{90}{\footnotesize\hspace{0.5in}Mode $1$}&
      \resizebox{!}{1.5in}{\includegraphics{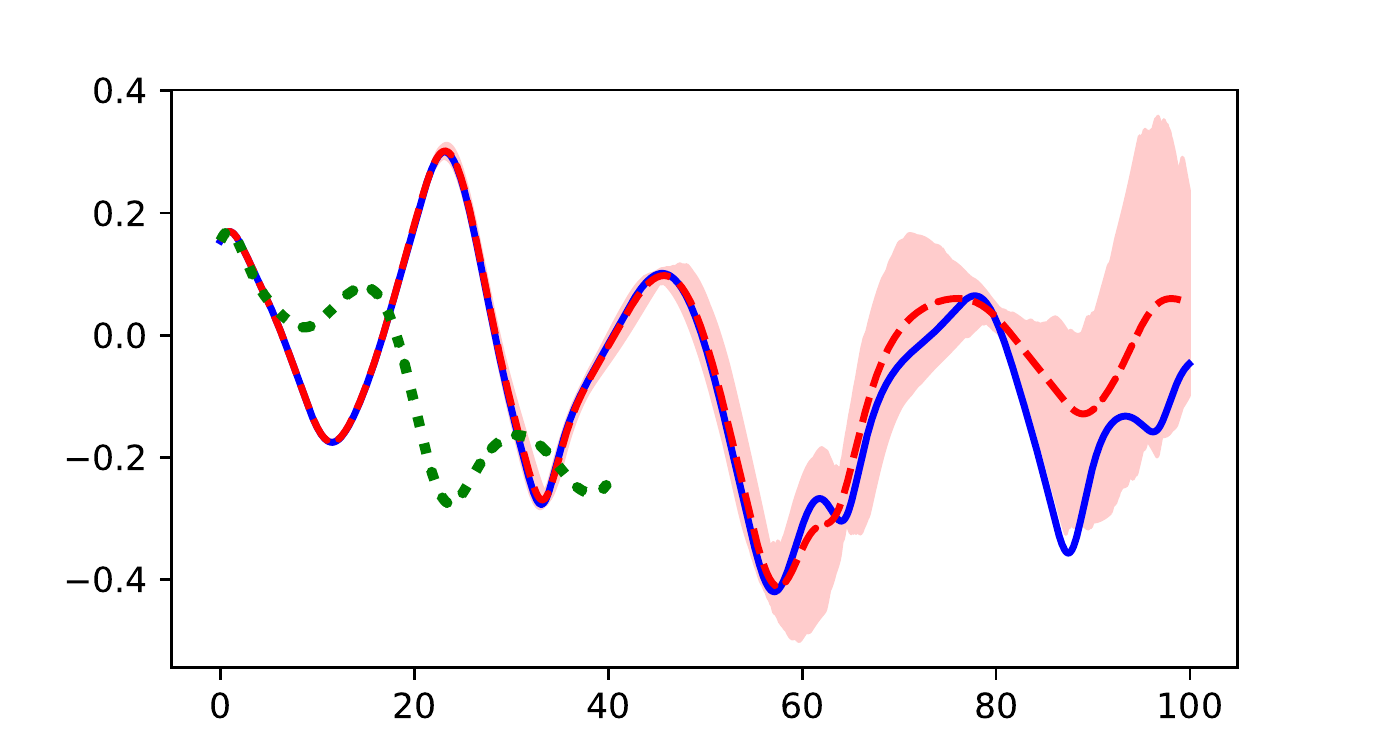}}&
      \resizebox{!}{1.5in}{\includegraphics{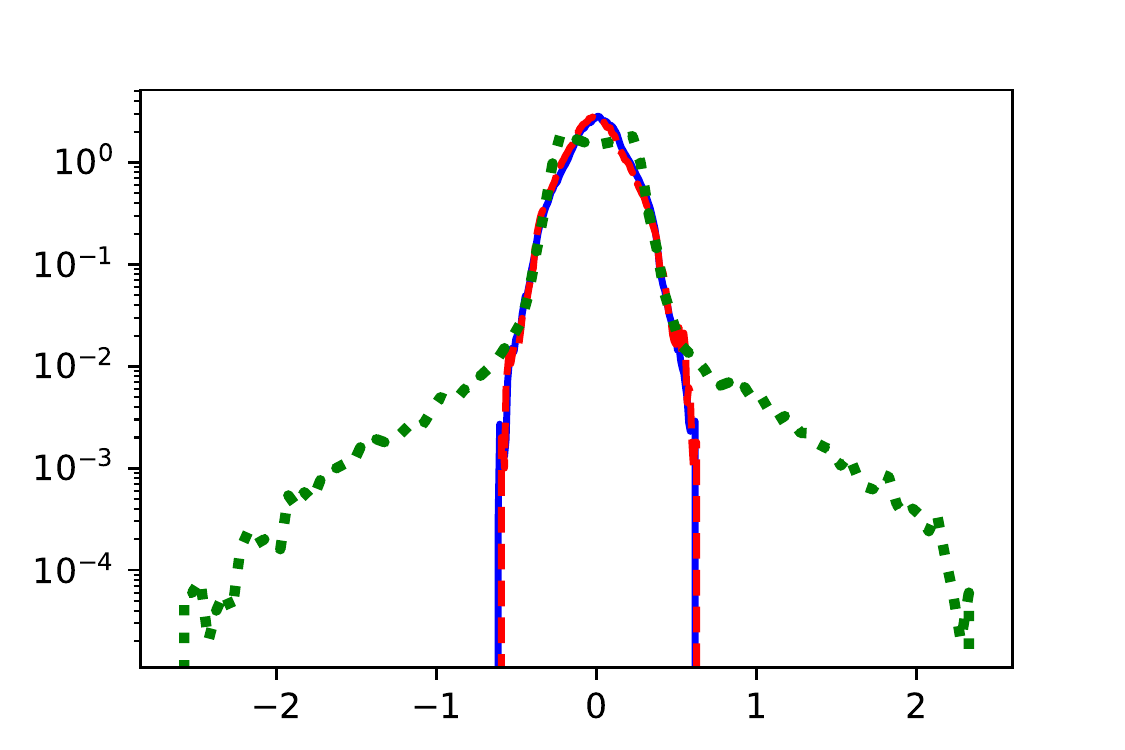}}\\[-13pt]
      \rotatebox{90}{\footnotesize\hspace{0.5in}Mode $2$}&
      \resizebox{!}{1.5in}{\includegraphics{./ks108-pqr331-traj-mode2}}&
      \resizebox{!}{1.5in}{\includegraphics{./ks108-pqr331-marginal-mode2}}\\[-13pt]
      \rotatebox{90}{\footnotesize\hspace{0.5in}Mode $3$}&
      \resizebox{!}{1.5in}{\includegraphics{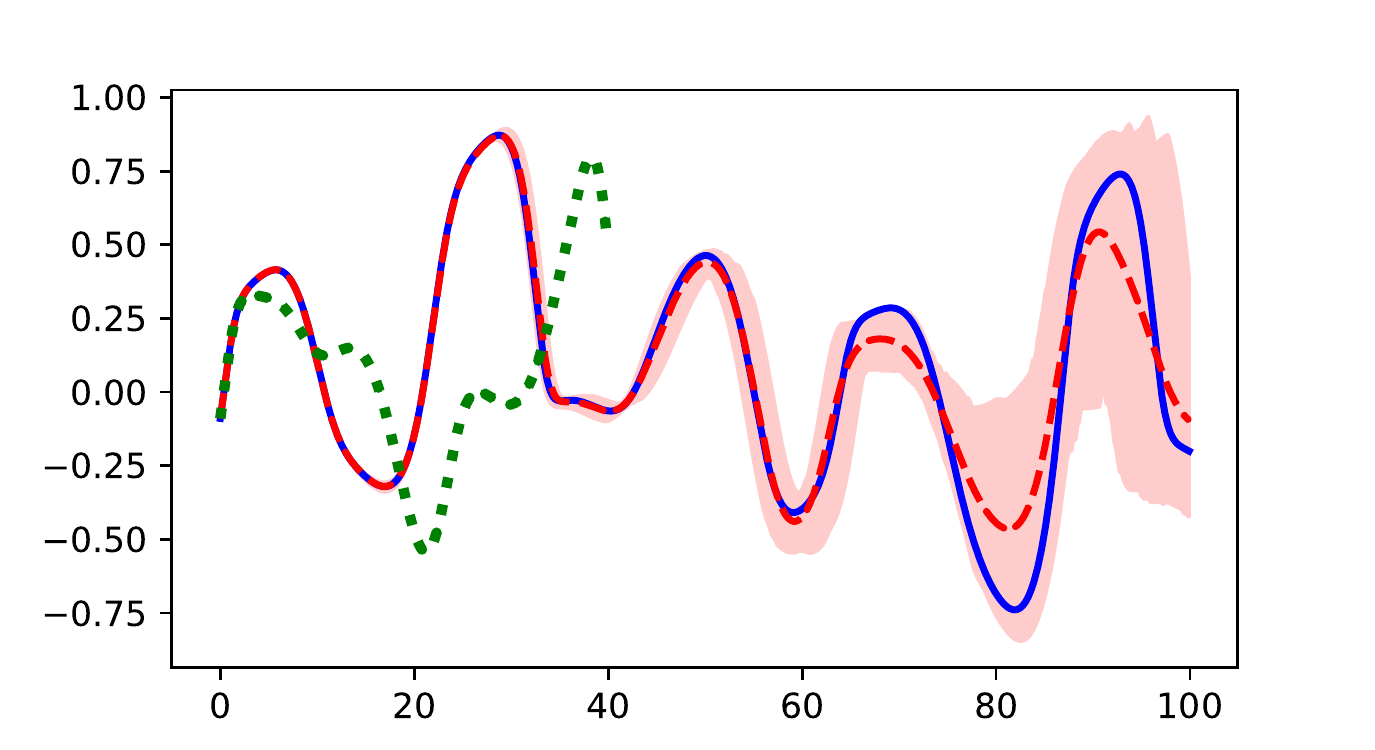}}&
      \resizebox{!}{1.5in}{\includegraphics{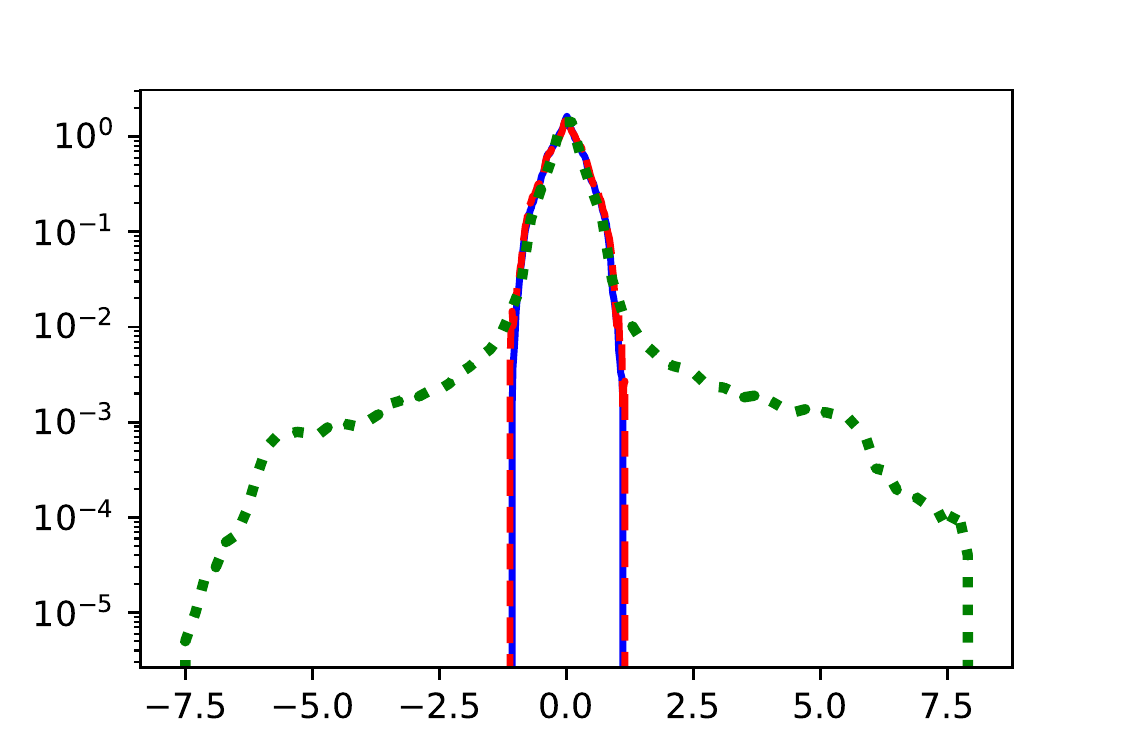}}\\[-13pt]
      \rotatebox{90}{\footnotesize\hspace{0.5in}Mode $4$}&
      \resizebox{!}{1.5in}{\includegraphics{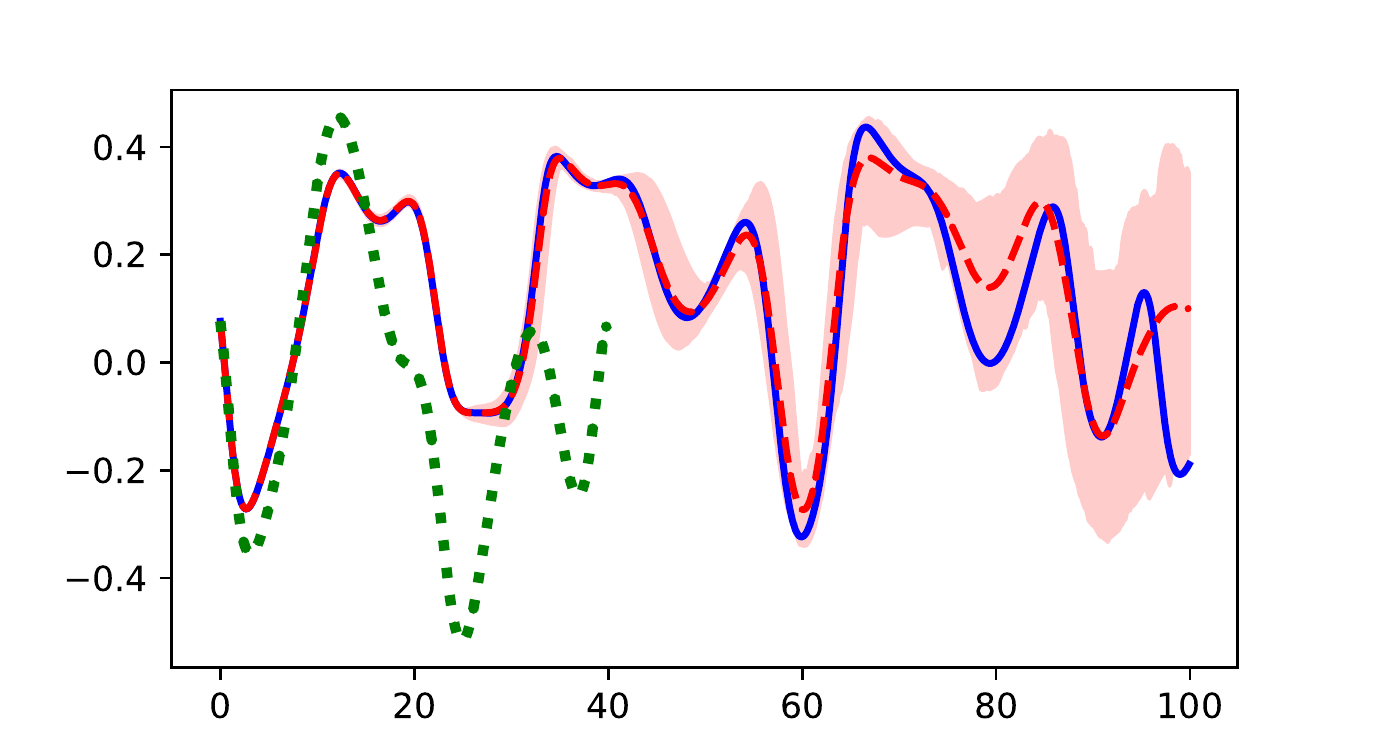}}&
      \resizebox{!}{1.5in}{\includegraphics{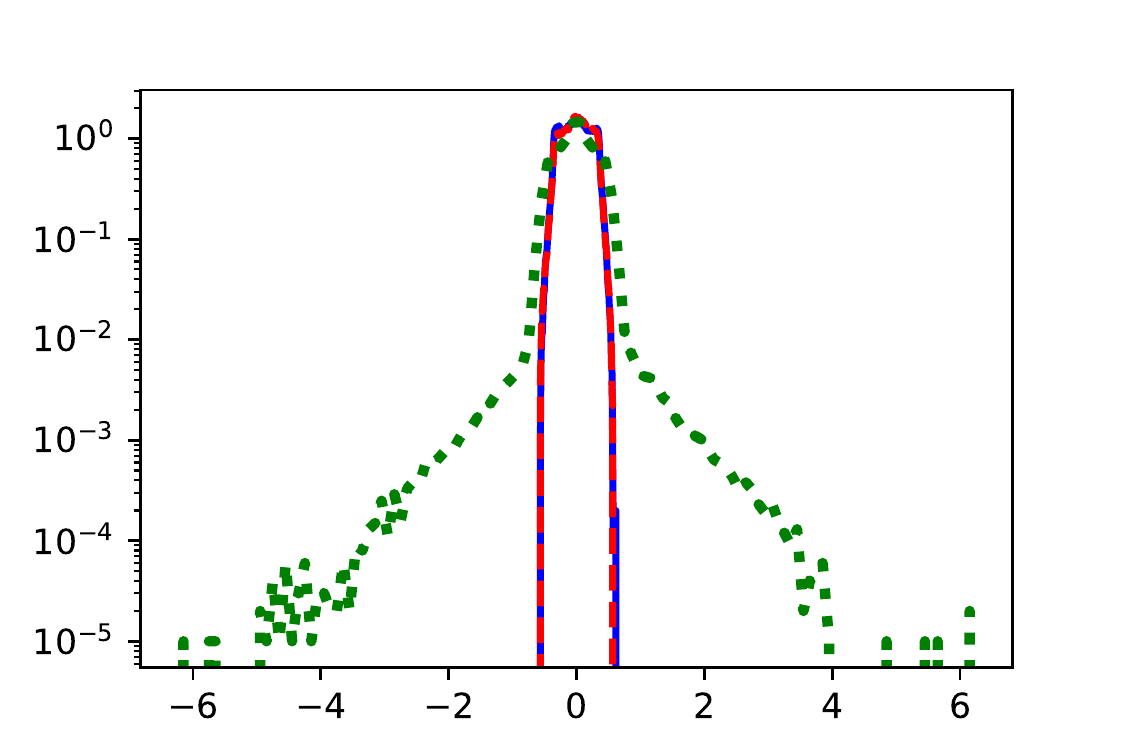}}\\[-13pt]
      \rotatebox{90}{\footnotesize\hspace{0.5in}Mode $5$}&
      \resizebox{!}{1.5in}{\includegraphics{./ks108-pqr331-traj-mode5}}&
      \resizebox{!}{1.5in}{\includegraphics{./ks108-pqr331-marginal-mode5}}\\[-6pt]
      &{\footnotesize Time} & {\footnotesize $Re(u_k)$}\\
    \end{tabular}
  \end{center}
  \caption{Comparison of finite-time forecasts and marginal
    distributions.  In all panels, solid blue line is the full model
    (108-mode truncation), dashed red line is the 5-mode reduced model,
    and dotted green line the 5-mode truncation.  \emph{Left:}
    trajectories starting from the same initial conditions.  For the
    reduced model, we show the 5th percentile, mean, and 95th
    percentile, computed with an ensemble of size 100.  The truncated
    model was terminated at $t=40$ to reduce clutter.  \emph{Right:}
    marginal densities.}
    
  \applabel{fig:ks-traj-pdf}
\end{figure}

\begin{figure}
  \begin{center}
    \begin{tabular}{r@{\hskip 0pt}c@{\hskip -6pt}c}
      &{\small ACF} & {\small CCF}\\
      \rotatebox{90}{\footnotesize\hspace{0.5in}Mode $1$}&
      \resizebox{!}{1.5in}{\includegraphics{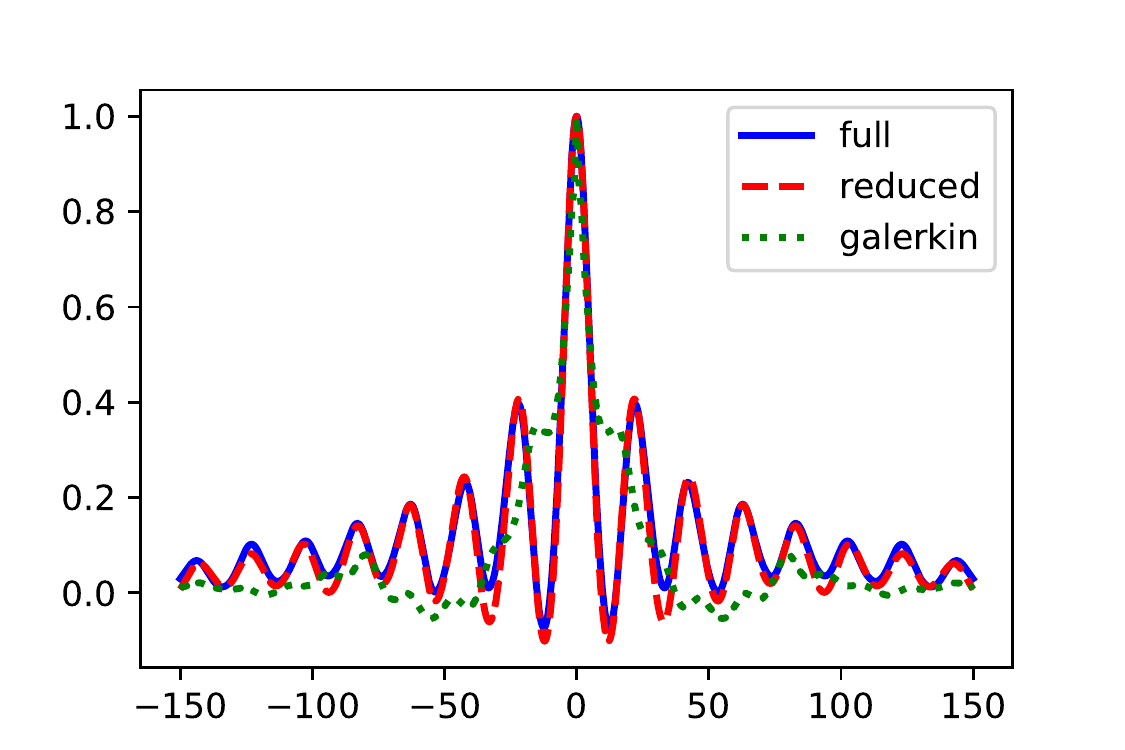}}&
      \resizebox{!}{1.5in}{\includegraphics{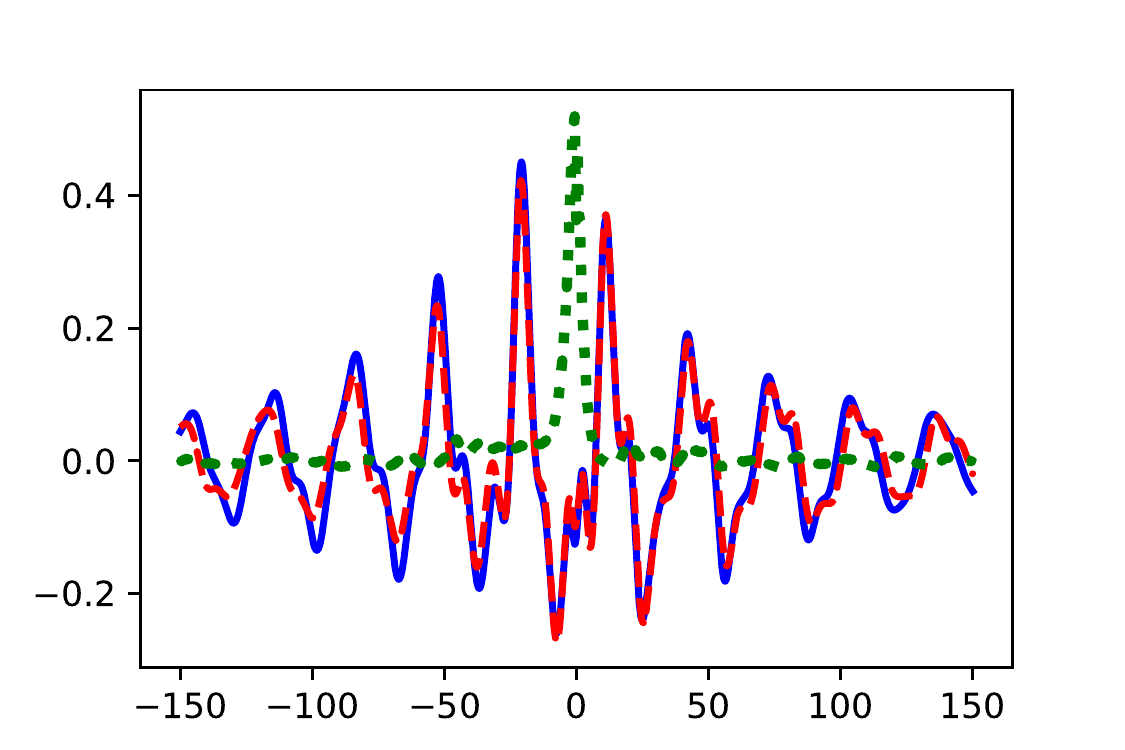}}\\[-13pt]
      \rotatebox{90}{\footnotesize\hspace{0.5in}Mode $2$}&
      \resizebox{!}{1.5in}{\includegraphics{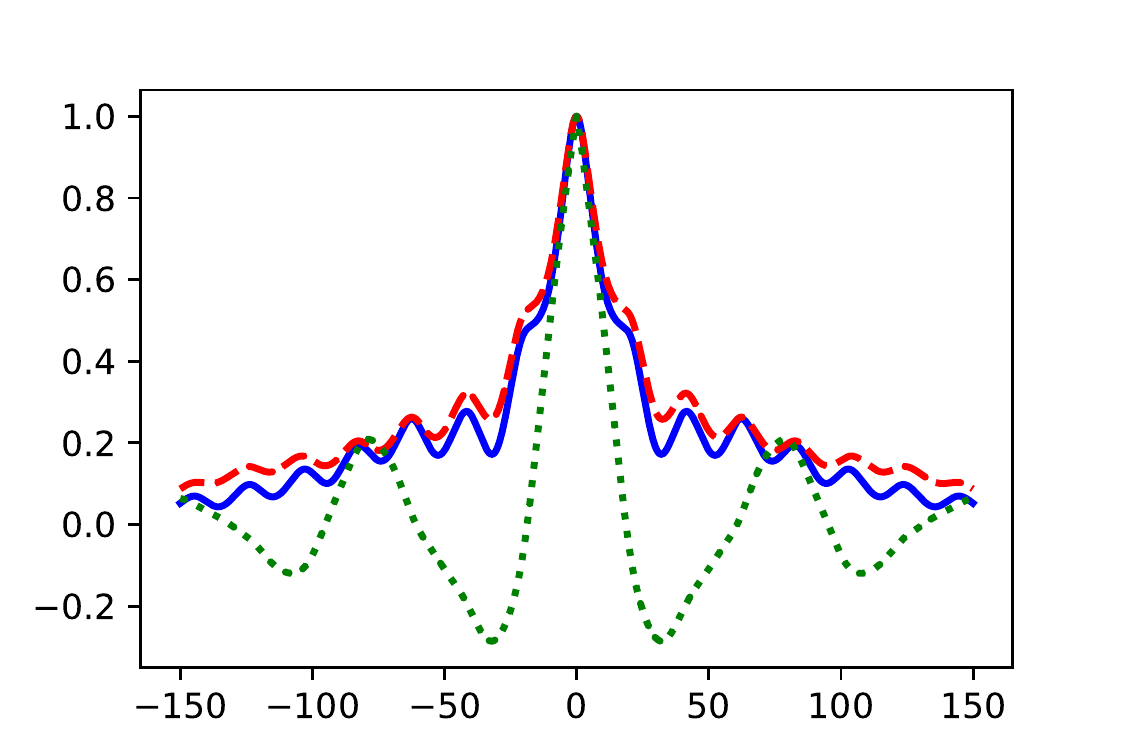}}&
      \resizebox{!}{1.5in}{\includegraphics{./ks108-pqr331-ccf-mode42}}\\[-13pt]
      \rotatebox{90}{\footnotesize\hspace{0.5in}Mode $3$}&
      \resizebox{!}{1.5in}{\includegraphics{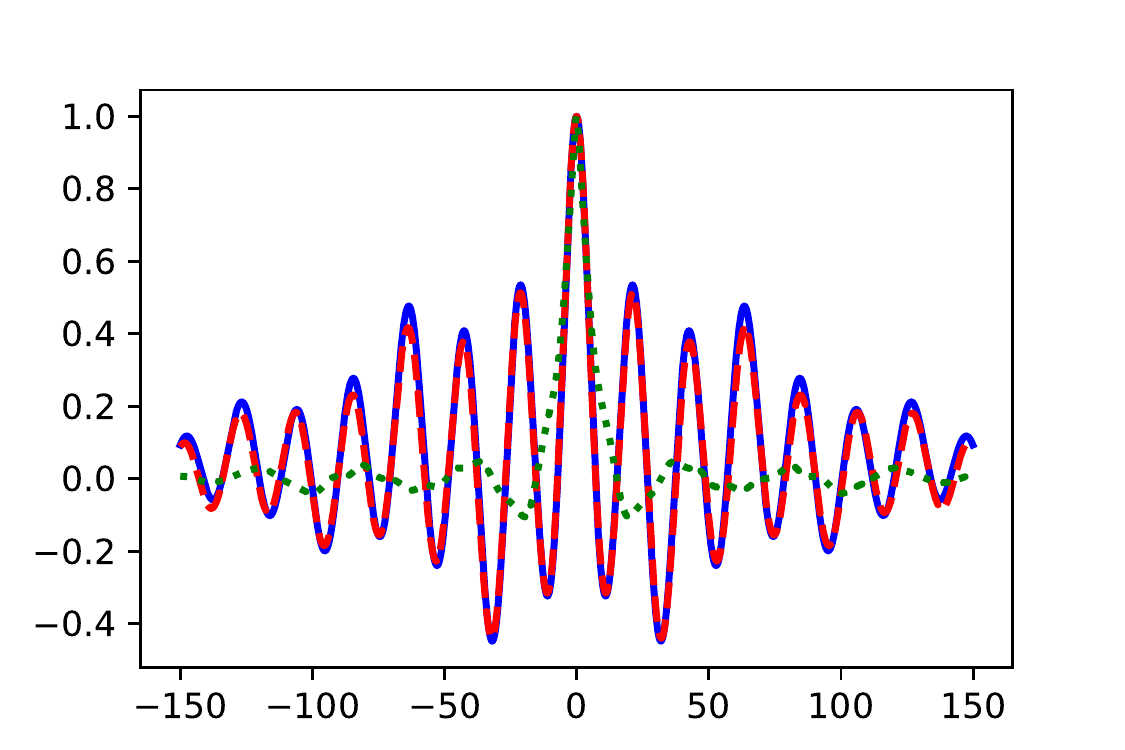}}&
      \resizebox{!}{1.5in}{\includegraphics{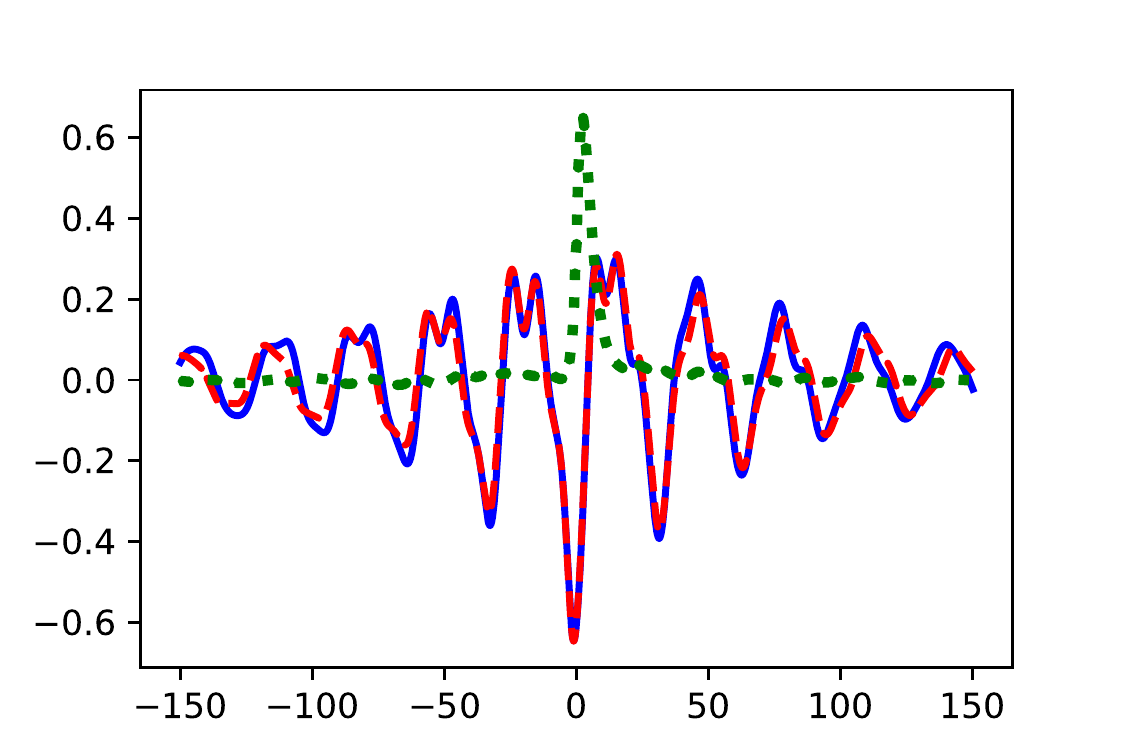}}\\[-13pt]
      \rotatebox{90}{\footnotesize\hspace{0.5in}Mode $4$}&
      \resizebox{!}{1.5in}{\includegraphics{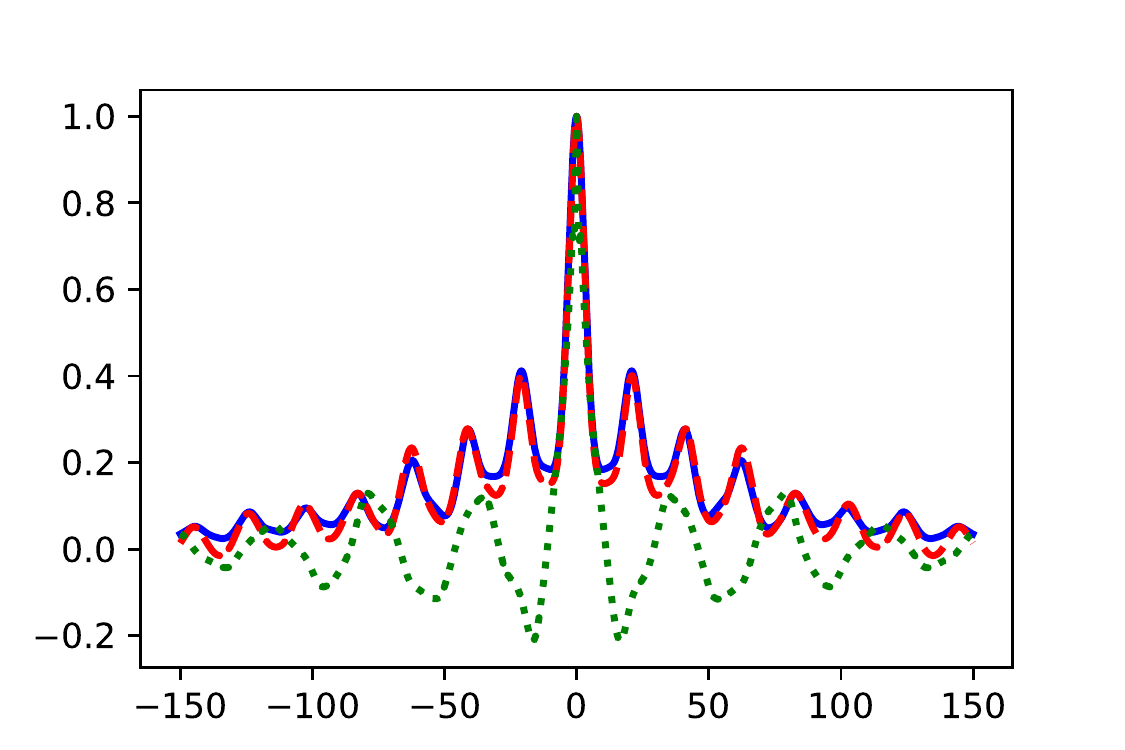}}&
      \resizebox{!}{1.5in}{\includegraphics{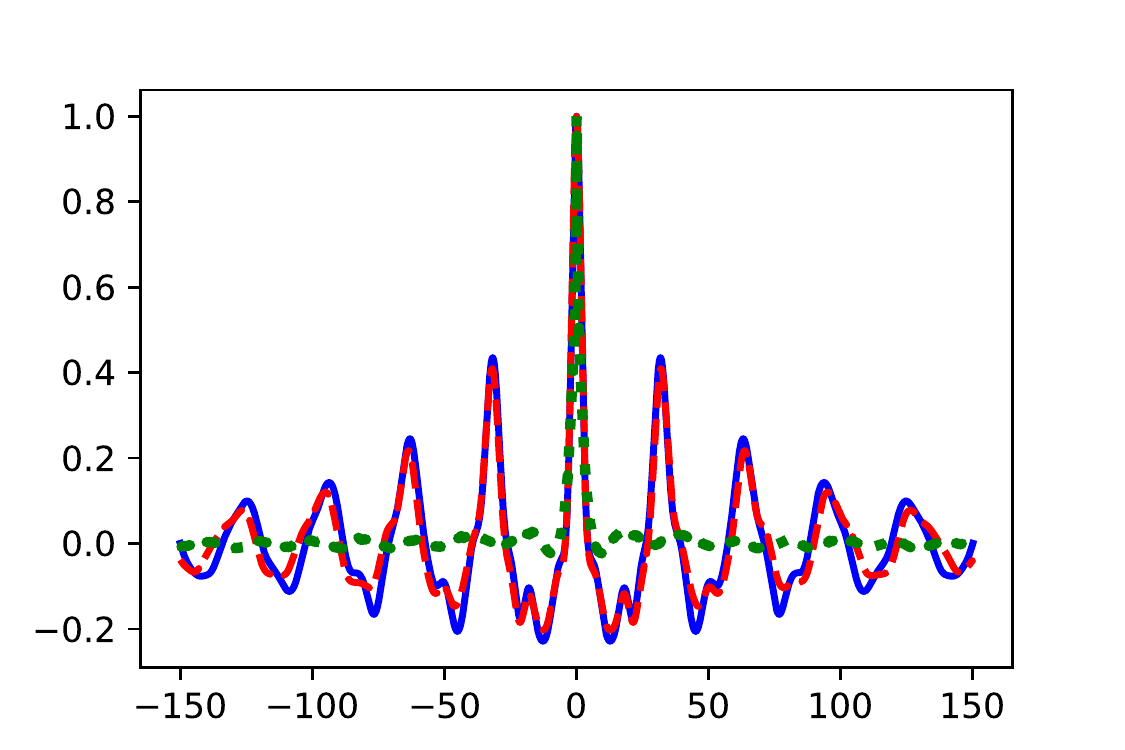}}\\[-13pt]
      \rotatebox{90}{\footnotesize\hspace{0.5in}Mode $5$}&
      \resizebox{!}{1.5in}{\includegraphics{./ks108-pqr331-acf-mode5}}&
      \resizebox{!}{1.5in}{\includegraphics{./ks108-pqr331-ccf-mode45}}\\[-6pt]
      &{\footnotesize Time lag} & {\footnotesize Time lag}\\
    \end{tabular}
  \end{center}
  \caption{Comparison of autocovariance functions (ACFs) and energy
    cross correlation functions (CCFs).  In all panels, solid blue line
    is the full model (108-mode truncation), dashed red line is the
    5-mode reduced model, and dotted green line the 5-mode truncation.
    \emph{Left:} Autocovariance functions for $Re(u_k(t))$ for
    $k=1,\cdots,5$.  \emph{Right:} Cross correlations between
    $|u_4(t)|^2$ and $|u_k(t)|^2$ for $k=1,\cdots,5$.}
    
  \applabel{fig:ks-stats}
\end{figure}

Figs.~\appref{fig:ks-traj-pdf} and \appref{fig:ks-stats} are full versions of
the numerical results shown in Sect.~\ref{sect:ks}.

\begin{figure}
  \begin{center}
    \begin{tabular}{cc}
      RMSE & ANCR \\[-0.5ex]
      \resizebox{!}{2in}{\includegraphics*[bb=4pt 4pt 3.5in 2.8in]{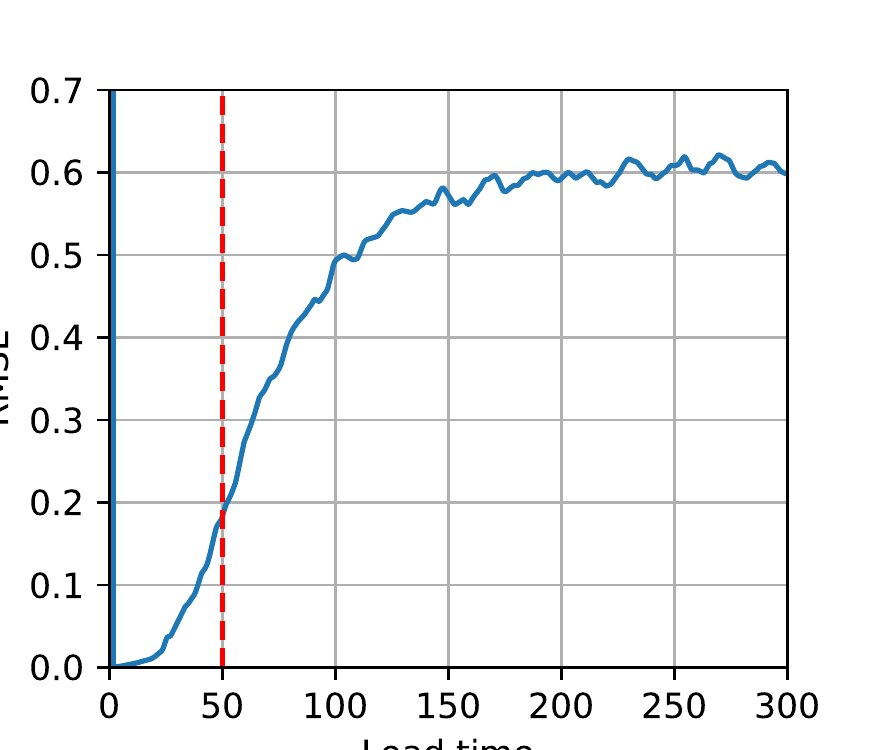}}&
      \resizebox{!}{2in}{\includegraphics*[bb=4pt 4pt 3.5in 2.8in]{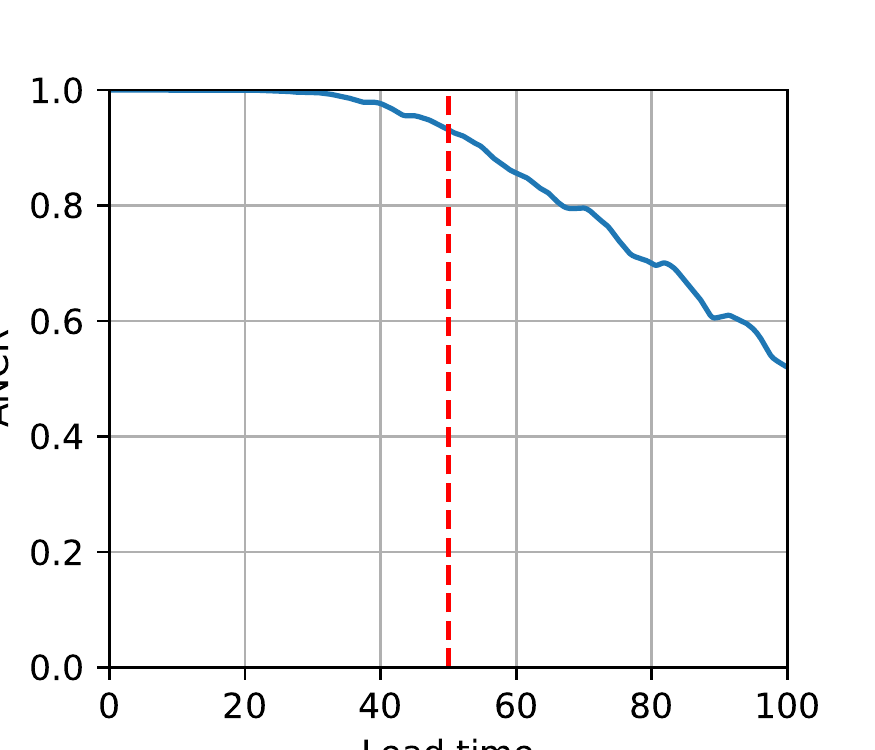}}\\
      Lead time & Lead time \\
    \end{tabular}
  \end{center}
  \caption{Forecasting skill as function of lead time of the reduced
    model for the KS equation.  \emph{Left:} root mean squared error
    (RMSE).  \emph{Right:} anomaly correlation (ANCR).  See text for
    details.}
  \applabel{fig:rmse}
\end{figure}

To further quantify finite-time forecasts as a function of the ``lead
time'' (i.e., time since initial observation), we introduce two standard
measures of forecasting ``skill,'' the root mean squared error and the
anomaly correlation.  Both are based on ensemble forecasts in the
following way: let $v(t_n)$ denote the time series data for the full
model, and take $N_{0}$ short pieces, i.e., $\left\{ \left(
v(t_{n}),n=n_{i},n_{i}+1,\dots ,n_{i}+T\right) \right\} _{i=1}^{N_{0}}$
with $n_{i+1}=n_{i}+T_{lag}/\dt$, where $T=$ $T_{lag}/\dt$ is the length
of each piece and $T_{lag}$ is the time gap between two adjacent pieces.
For each short piece $\left( v(t_{n}),n=n_{i},\dots ,n_{i}+T\right) $,
we generate $N_{ens}$ trajectories of length $T$ from the reduced model,
starting all ensemble members from the same initial segment $\left(
v(t_{n_{i}}),v(t_{n_{i}+1}),\dots ,v(t_{n_{i}+m})\right) $, where
$m=2p+1$, and denote the sample trajectories by $\left(
u^{n}(i,j),n=1,\dots ,T\right) $ for $i=1,\dots ,N_{0}$ and $j=1,\dots
,N_{ens}$.

Again, we do not introduce artificial perturbations into the initial
conditions, because the exact initial conditions are known, and by
initializing from data, we preserve the memory of the system so as to
generate better ensemble trajectories.

The \emph{root mean squared error} is
\begin{equation}
  \mathrm{RMSE}(\tau _{n}):=\left(
  \frac{1}{N_{0}}\sum_{i=1}^{N_{0}}\left\vert \operatorname{Re}
  v(t_{n_{i}+n})-\operatorname{Re}\bar{u}^{n}(i)\right\vert
  ^{2}\right)^{1/2},
\end{equation} 
where $\tau _{n}=n\dt$, $\bar{u}^{n}(i)=
\frac{1}{N_{ens}}\sum_{j=1}^{N_{ens}}u^{n}(i,j)$, and the \emph{anomaly
  correlation} (see, e.g., \cite{CVE08}) is
\begin{equation}
  \mathrm{ANCR}(\tau _{n}):=\frac{1}{N_{0}}\sum_{i=1}^{N_{0}}\frac{\mathbf{a}
    ^{v,i}(n)\cdot \mathbf{a}^{u,i}(n)}{\sqrt{|\mathbf{a}^{v,i}(n)|^{2}\left
      \vert \mathbf{a}^{u,i}(n)\right\vert ^{2}}},
\end{equation}
where
$\mathbf{a}^{v,i}(n)=\operatorname{Re}v(t_{n_{i}+n})-\operatorname{Re}\left\langle
v\right\rangle $ and
$\mathbf{a}^{u,i}(n)=\operatorname{Re}\bar{u}^{n}(i)-\operatorname{Re}
\left\langle v\right\rangle $ are the anomalies in data and the ensemble
mean.  Here $\mathbf{a\cdot b=}\sum_{k=1}^{K}a_{k}b_{k}$, $\left\vert
\mathbf{ a}\right\vert ^{2}=\mathbf{a\cdot a}$, and $\left\langle
v\right\rangle $ is the time average of the long trajectory of $v$.
Both statistics measure the accuracy of the mean ensemble prediction:
the RMSE measures, in an average sense, the difference between the mean
ensemble trajectory, and the ANCR shows the average correlation between
the mean ensemble trajectory and the true data trajectory.
$\rm{RMSE}=0$ and $\rm{ANCR}=1$ would correspond to a perfect
prediction, and small RMSEs and large (close to 1) ANCRs are desired.

For our reduced model, we computed the RMSE and ANCR using ensembles of
$N_{ens}=100$ trajectories with independent initial conditions.
Fig.~\appref{fig:rmse} (left) shows the RMSE and ANCR for a range of
lead times.  As expected, the RMSE increases with lead time, and
consistent with Fig.~\ref{fig:ks-solutions}(a), it is relatively small
compared to its apparent asymptotic value (about 0.6) for lead times
$<50$.  The ANCR in Fig.~\appref{fig:rmse} (right) corroborates this.
The two figures are comparable to Fig.~5 of \cite{LLC17} and show very
similar trends.

\heading{Role of the noise terms $\xi_n$.}  We experimented with running
the reduced model with $\xi_n\equiv0$, i.e., without any noise term.
This does not appreciably change the ACF or marginal distributions, nor
the forecasting skill of the reduced model.  However, the kind of
ensemble prediction and uncertainty quantification illustrated in
Fig.~\appref{fig:ks-traj-pdf} cannot be carried out without noise terms
calibrated to the reduced model.

\section{Stochastic Burgers equation}
\label{app:burgers}

The nonlinear terms $\{\Psi_{n-j}\}$ in~\meqref{eq:burgers-nonlinear-terms} are defined by
\begin{equation*}%
  \Psi^a_{n-j} = u^{n-j}~,~~
  \Psi^b_{n-j} = R^{\dt}(u^{n-j})~,~~\mbox{and}~~
  \Psi^c_{n-j,k} =  \sum_{\substack{ |k-l|\leq K, K< |l| \leq 2K \\ \text{ or }  |l|\leq K, K< |k-l| \leq 2K} }\wu^{n-1}_l \wu^{n-j}_{k-l} \text{ for} ~~k=1,\cdots,K,
\end{equation*}
where the terms $\{\wu\}$ are defined as
\begin{align}  \applabel{eq:Burgers-ansatz-d}
  \widetilde{u}^{n-j}_{k }= 
  \left\{
  \begin{array}{ll}
    u^{n-j}_k~, & 1\leq k\leq K; \\[1ex]
    \frac{i\lambda_k}{2}e^{-\nu \lambda_k^2 j\delta} \sum_{\substack{ | l |\leq K,  \\ |k-l|\leq K} } \wv^{n-j}_{k-l} \wv^{n-j}_{l}, & K < k \leq 2K.
  \end{array}
  \right.
\end{align}
These terms resemble those in Eq.~(\appref{eq:ks-ansatz-d}) as they are also
introduced to represent the high modes by the low modes. But there is a
major difference: they represent the high modes as a functional of the
history of the low modes, rather than a function of the current state of
the low modes. This is due to the lack of an inertial manifold for the
Burgers equation, unlike the KSE.  These terms are derived from an
Riemann sum approximation of the integral equation for the high modes,
with suitable linear parametrization of the quadratic interaction.  A
detailed derivation of the ansatz is presented in a forthcoming paper.

Figs.~\appref{fig:burgers-traj}--~\appref{fig:burgers-ccf} show numerical
results for the stochastic Burgers equation.

\begin{figure}
  \begin{center}
    \begin{tabular}{r@{\hskip 0pt}c@{\hskip -6pt}r@{\hskip 0pt}c@{\hskip 0pt}r@{\hskip 0pt}c}
      \rotatebox{90}{\footnotesize\hspace{0.5in}Mode $1$}&
      \resizebox{2.3in}{!}{\includegraphics{./burgers-pqr110-response1}}&
      \rotatebox{90}{\footnotesize\hspace{0.5in}Mode $2$}&
      \resizebox{2.3in}{!}{\includegraphics{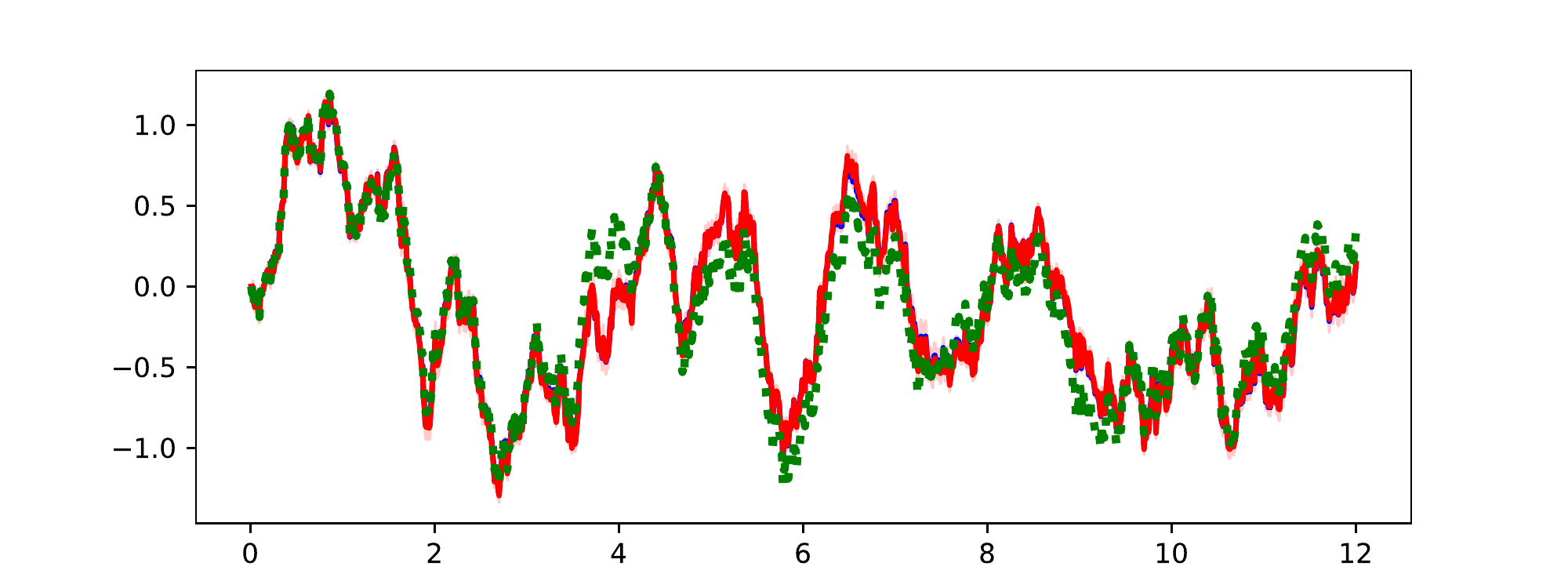}}&
      \rotatebox{90}{\footnotesize\hspace{0.5in}Mode $3$}&
      \resizebox{2.3in}{!}{\includegraphics{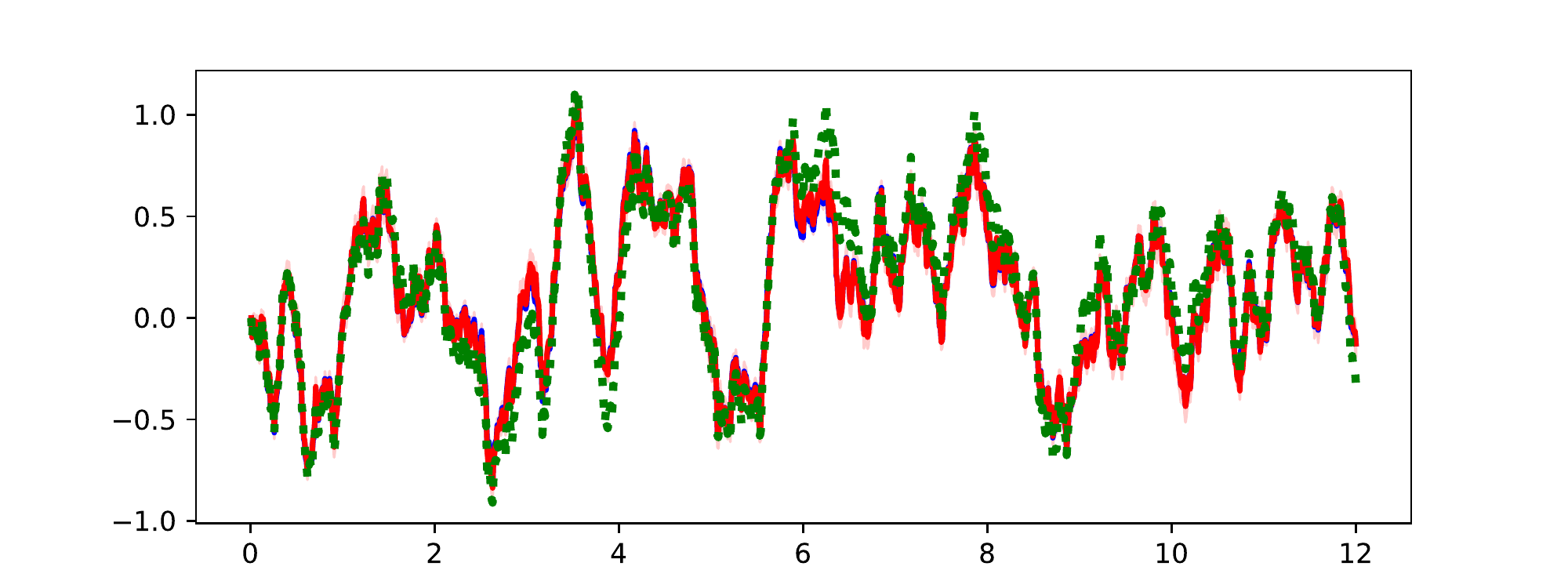}}\\[-6pt]
      \rotatebox{90}{\footnotesize\hspace{0.5in}Mode $4$}&
      \resizebox{2.3in}{!}{\includegraphics{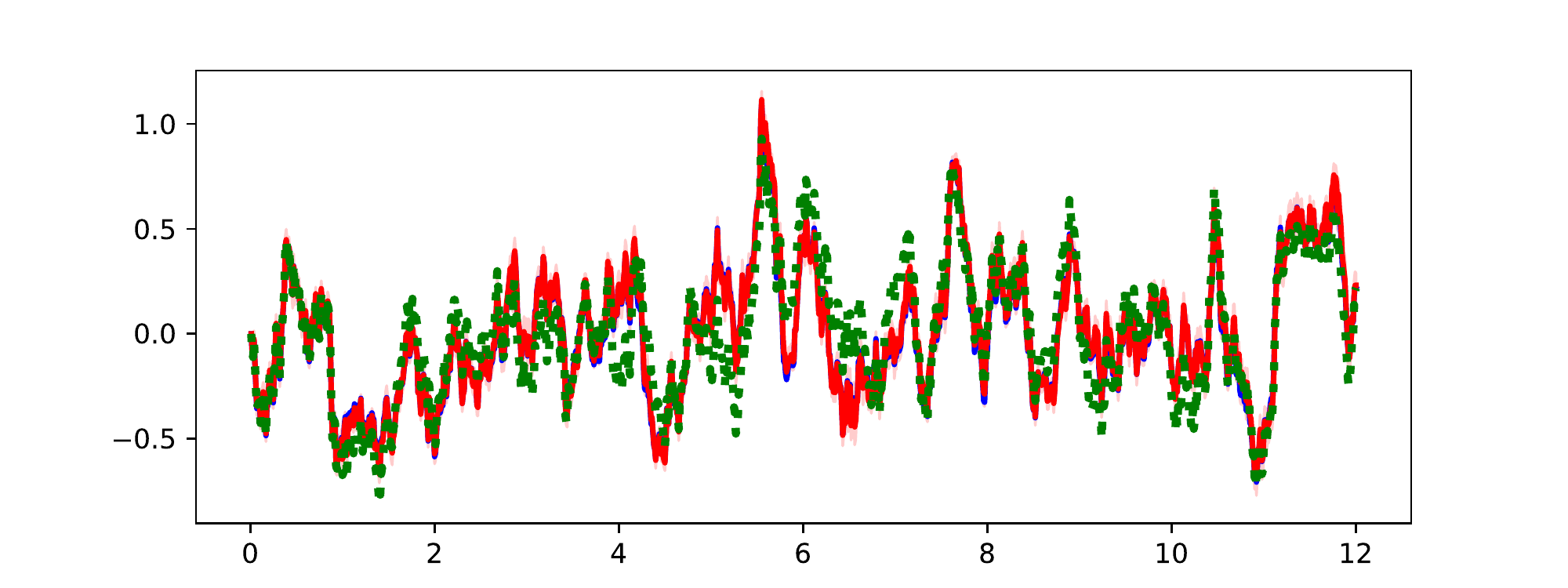}}&
      \rotatebox{90}{\footnotesize\hspace{0.5in}Mode $5$}&
      \resizebox{2.3in}{!}{\includegraphics{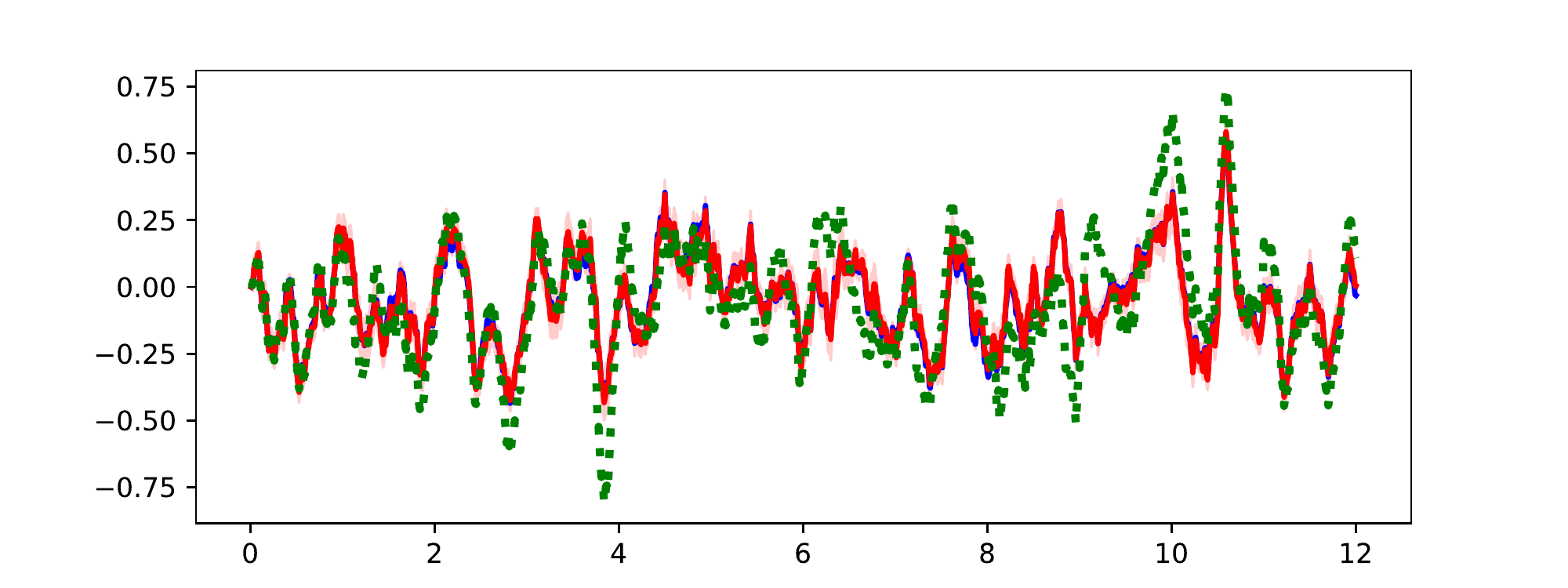}}&
      \rotatebox{90}{\footnotesize\hspace{0.5in}Mode $6$}&
      \resizebox{2.3in}{!}{\includegraphics{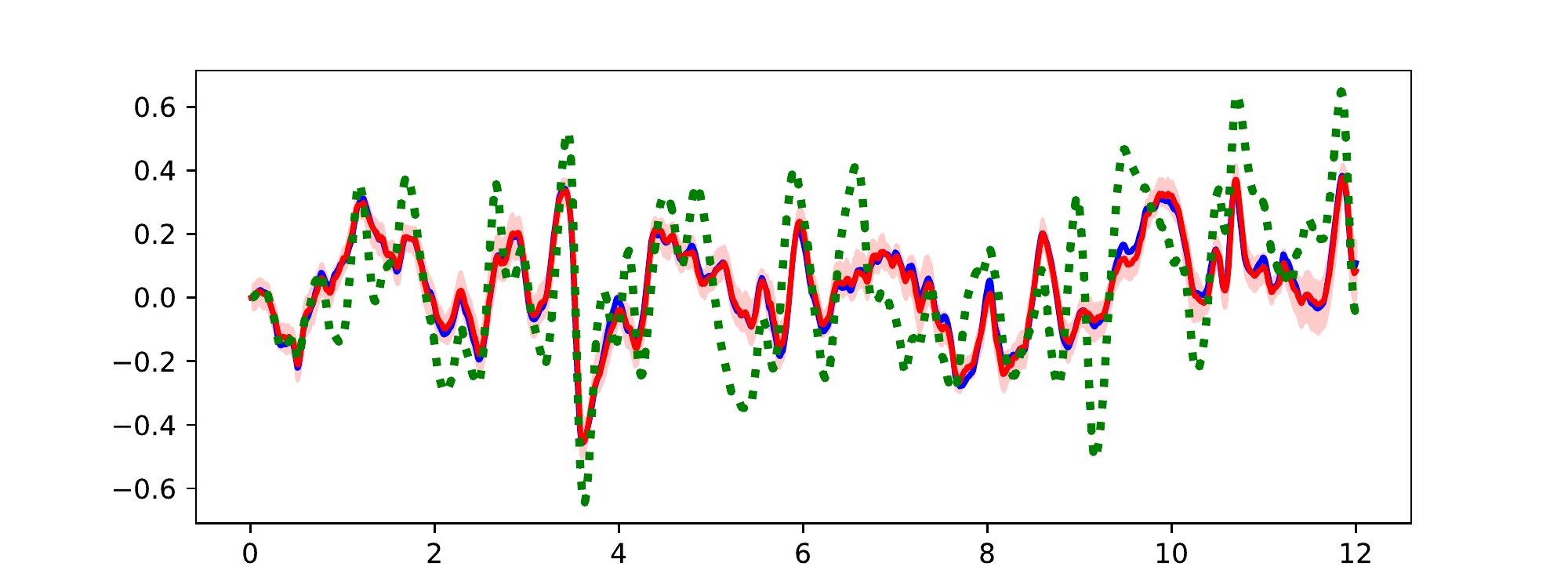}}\\[-6pt]
      \rotatebox{90}{\footnotesize\hspace{0.5in}Mode $7$}&
      \resizebox{2.3in}{!}{\includegraphics{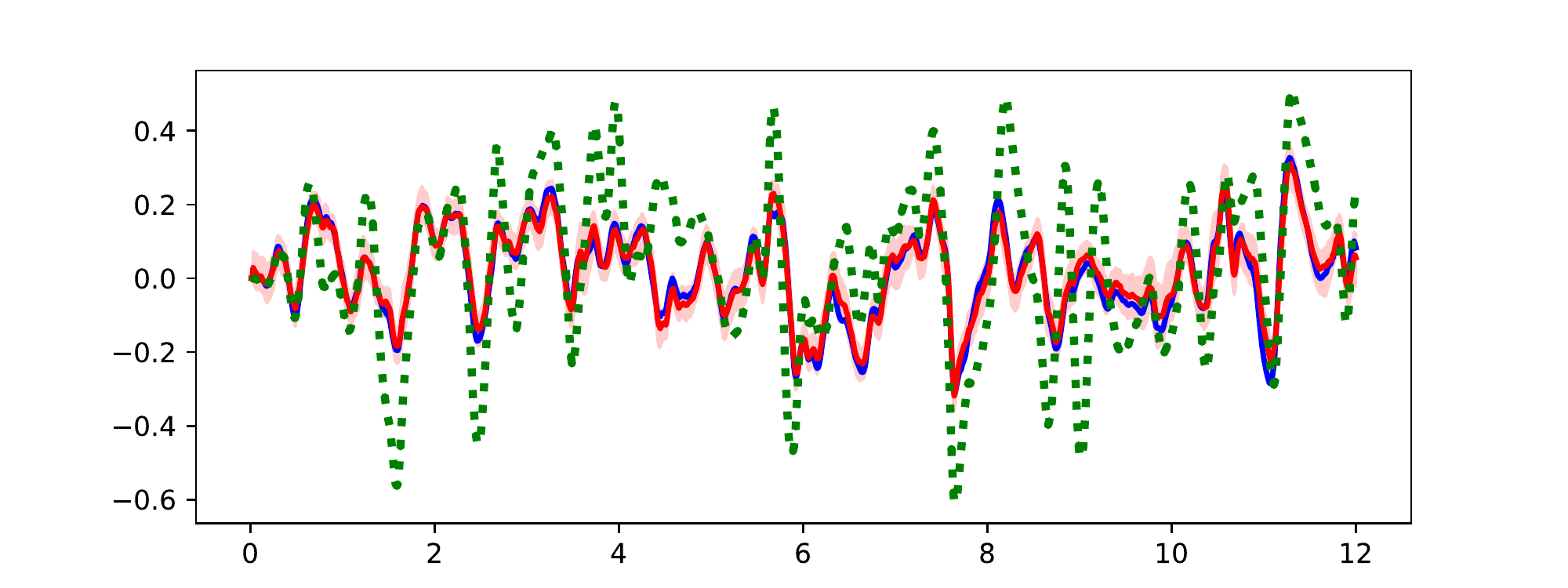}}&
      \rotatebox{90}{\footnotesize\hspace{0.5in}Mode $8$}&
      \resizebox{2.3in}{!}{\includegraphics{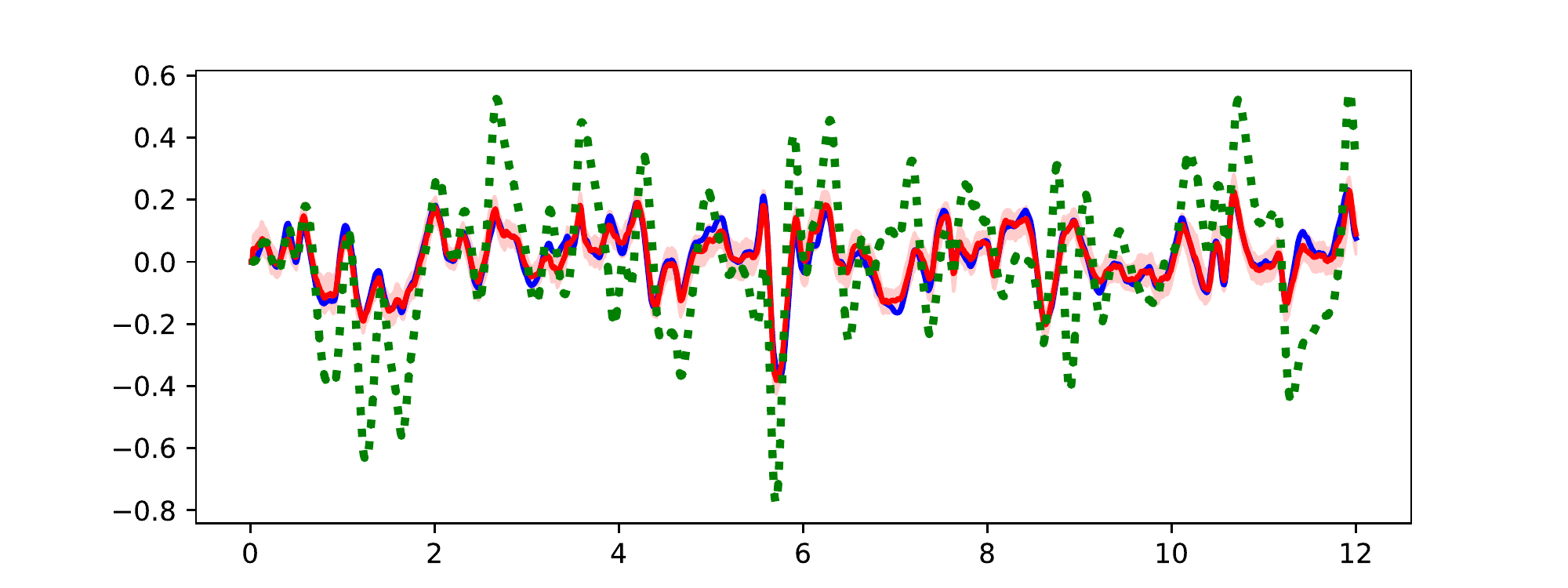}}&
      \rotatebox{90}{\footnotesize\hspace{0.5in}Mode $9$}&
      \resizebox{2.3in}{!}{\includegraphics{./burgers-pqr110-response9}}\\
      &{\footnotesize Time} && {\footnotesize Time} && {\footnotesize Time}\\
    \end{tabular}
  \end{center}
  \caption{Response forecasting for the stochastic Burgers equation.
    For $k=1,\cdots,9$, we plot $Re(u_k(t))$ as functions of $t$.  In
    all panels, solid blue line is the full model (128-mode truncation),
    dashed red line is the 9-mode reduced model, and dotted green line
    the 9-mode Galerkin truncation.  Initial transients ($t<8$) are not
    shown.}
    
  \applabel{fig:burgers-traj}
\end{figure}

\begin{figure}
  \begin{center}
    \begin{tabular}{r@{\hskip 0pt}c@{\hskip -6pt}r@{\hskip 0pt}c@{\hskip 0pt}r@{\hskip 0pt}c}
      \rotatebox{90}{\footnotesize\hspace{0.5in}Mode $1$}&
      \resizebox{2.3in}{!}{\includegraphics{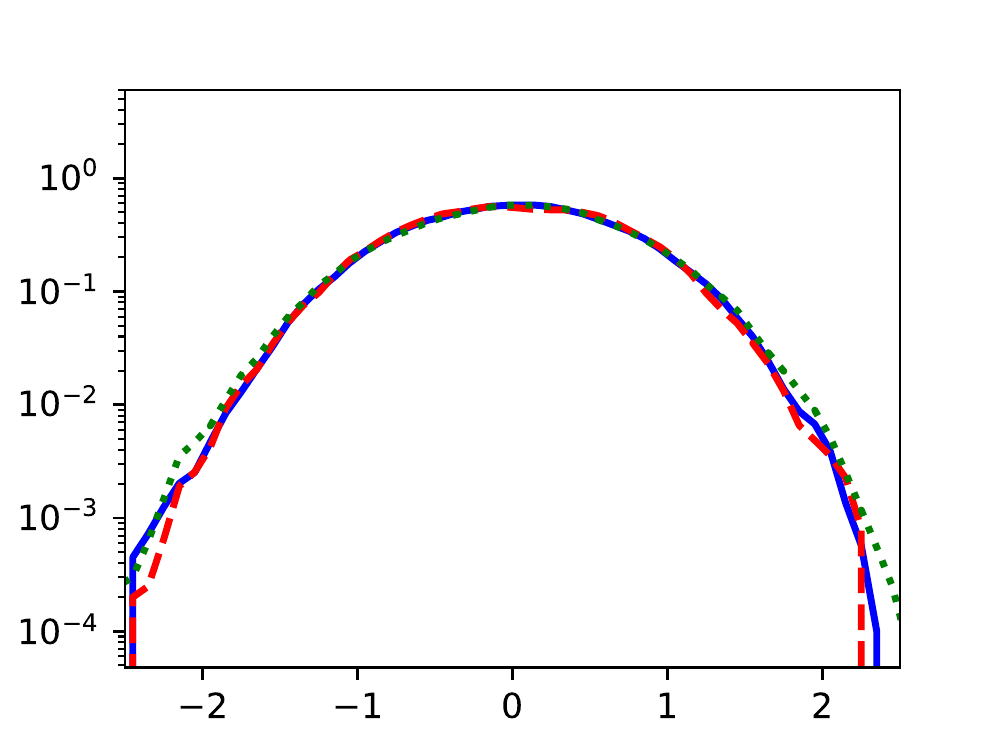}}&
      \rotatebox{90}{\footnotesize\hspace{0.5in}Mode $2$}&
      \resizebox{2.3in}{!}{\includegraphics{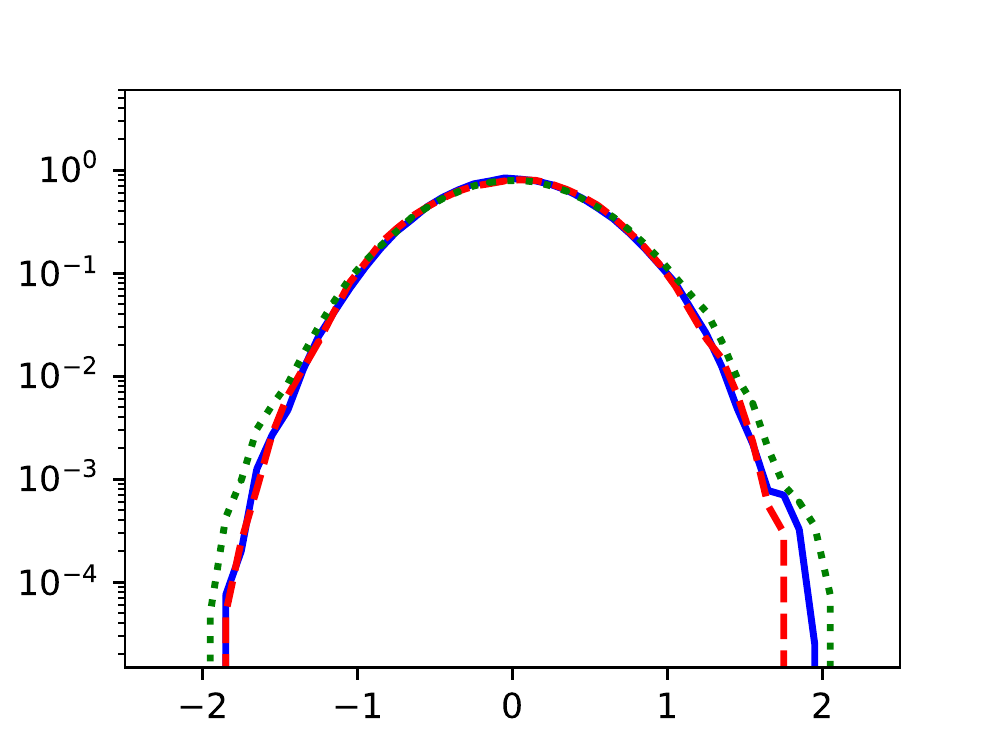}}&
      \rotatebox{90}{\footnotesize\hspace{0.5in}Mode $3$}&
      \resizebox{2.3in}{!}{\includegraphics{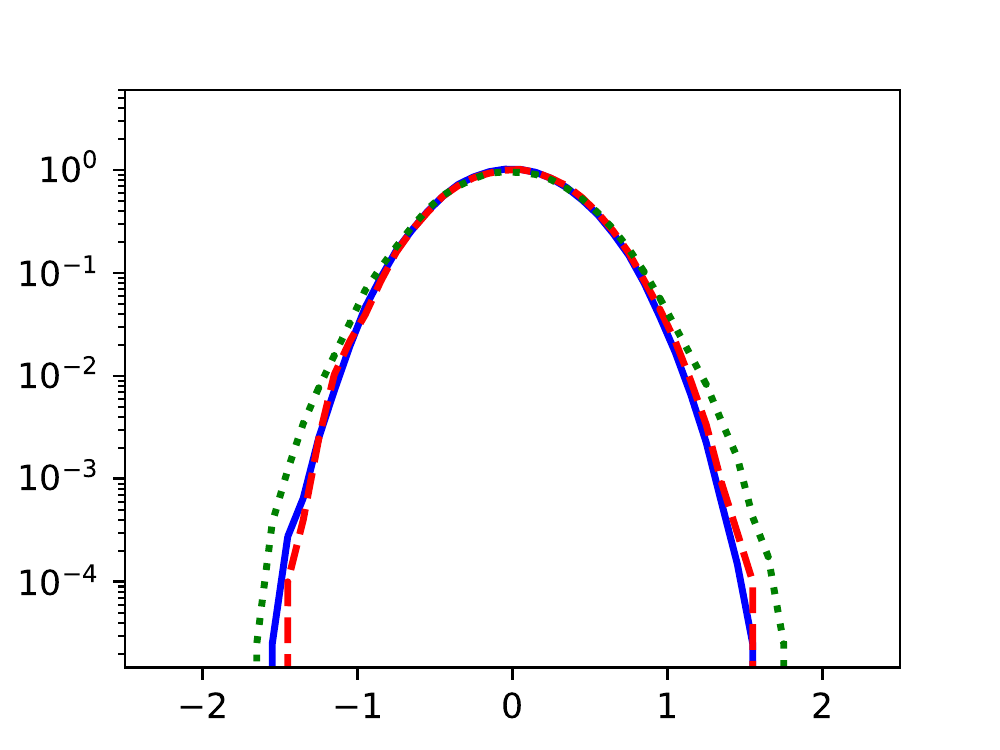}}\\[-6pt]
      \rotatebox{90}{\footnotesize\hspace{0.5in}Mode $4$}&
      \resizebox{2.3in}{!}{\includegraphics{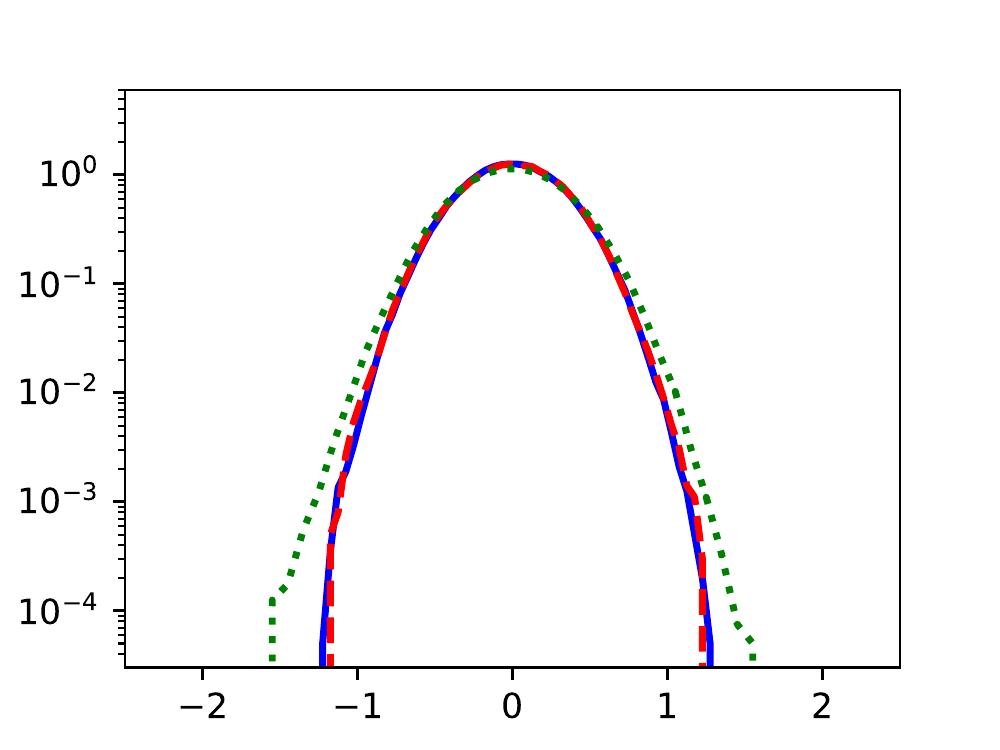}}&
      \rotatebox{90}{\footnotesize\hspace{0.5in}Mode $5$}&
      \resizebox{2.3in}{!}{\includegraphics{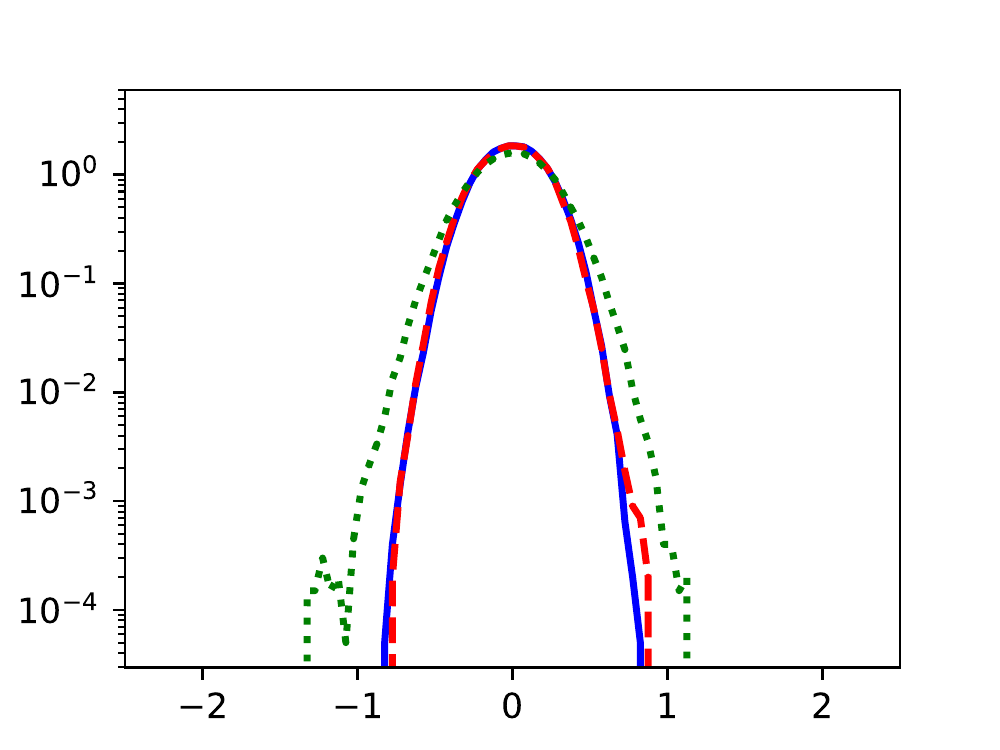}}&
      \rotatebox{90}{\footnotesize\hspace{0.5in}Mode $6$}&
      \resizebox{2.3in}{!}{\includegraphics{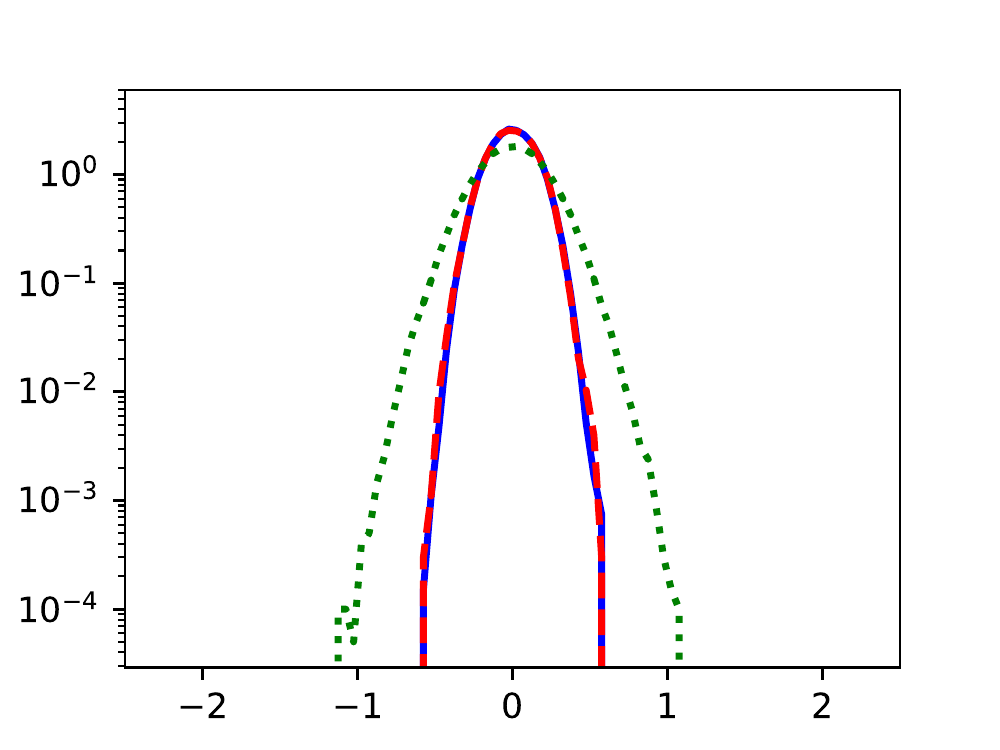}}\\[-6pt]
      \rotatebox{90}{\footnotesize\hspace{0.5in}Mode $7$}&
      \resizebox{2.3in}{!}{\includegraphics{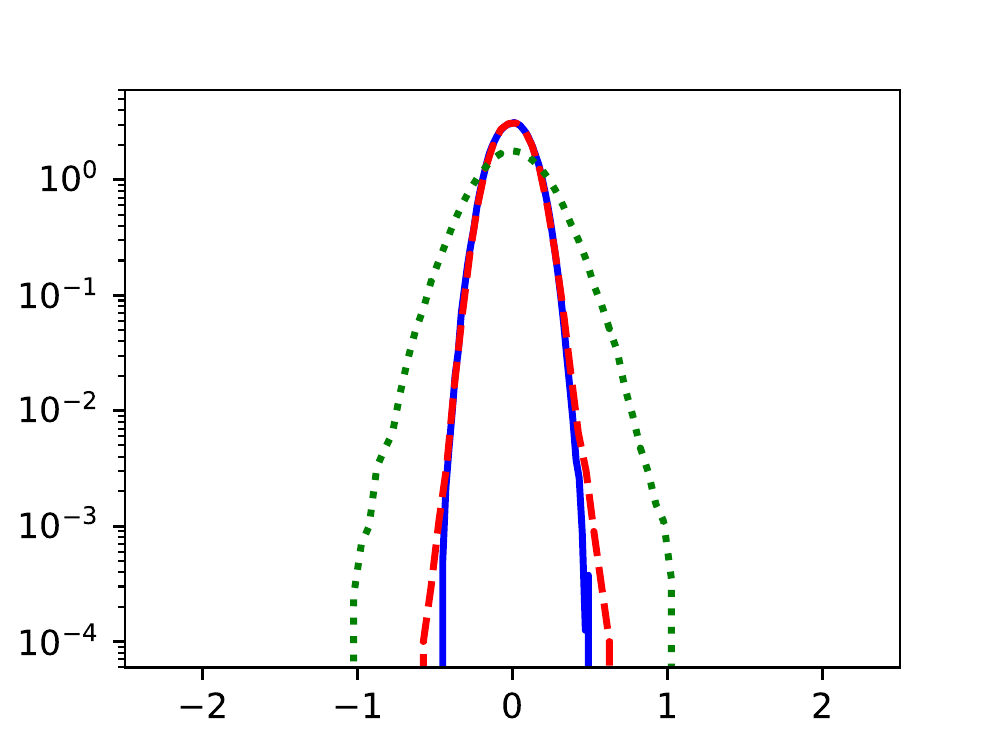}}&
      \rotatebox{90}{\footnotesize\hspace{0.5in}Mode $8$}&
      \resizebox{2.3in}{!}{\includegraphics{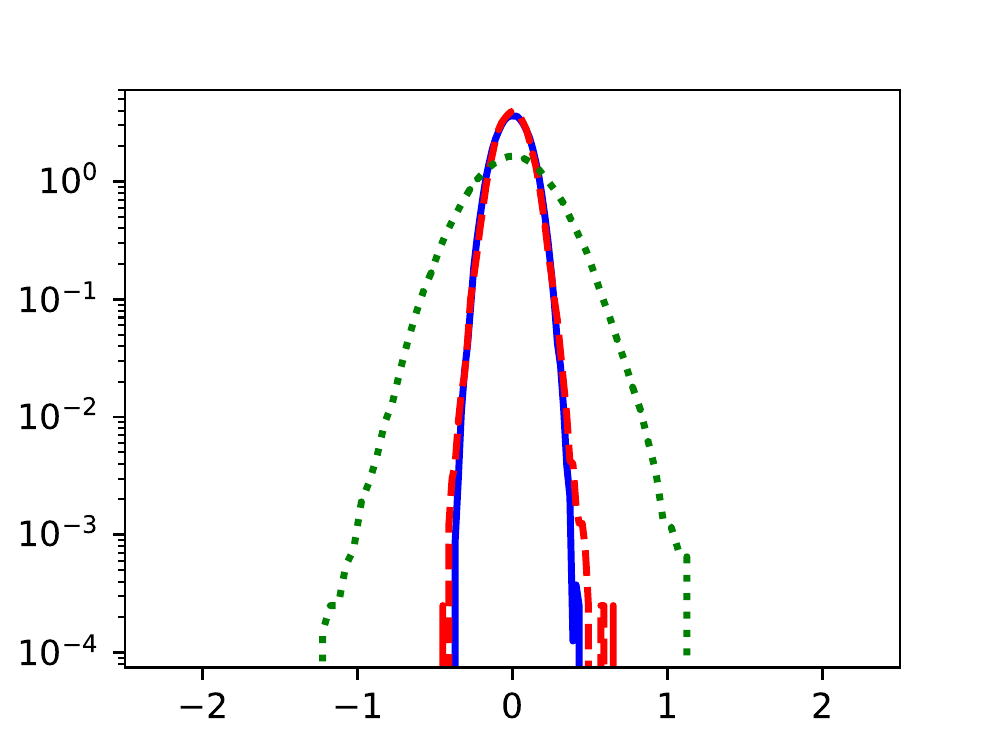}}&
      \rotatebox{90}{\footnotesize\hspace{0.5in}Mode $9$}&
      \resizebox{2.3in}{!}{\includegraphics{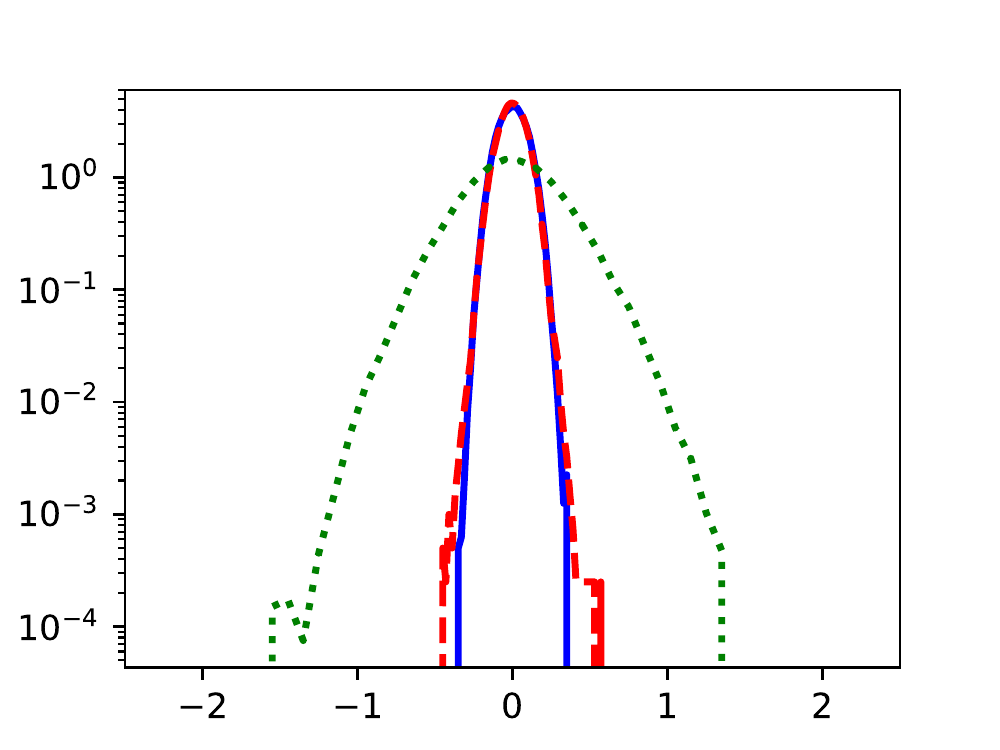}}\\
      &{\small $Re(u_k)$} && {\small $Re(u_k)$} && {\small $Re(u_k)$}\\
    \end{tabular}
  \end{center}
  \caption{Marginal densities for the stochastic Burgers equation.  We
    plot estimated densities for $Re(u_k)$ for $k=1,\cdots,9$.  In all
    panels, solid blue line is the full model (128-mode truncation),
    dashed red line is the 9-mode reduced model, and dotted green line
    the 9-mode Galerkin truncation.}
    
  \applabel{fig:burgers-pdf}
\end{figure}

\begin{figure}
  \begin{center}
    \begin{tabular}{r@{\hskip 0pt}c@{\hskip -6pt}r@{\hskip 0pt}c@{\hskip 0pt}r@{\hskip 0pt}c}
      \rotatebox{90}{\footnotesize\hspace{0.5in}Mode $1$}&
      \resizebox{2.3in}{!}{\includegraphics{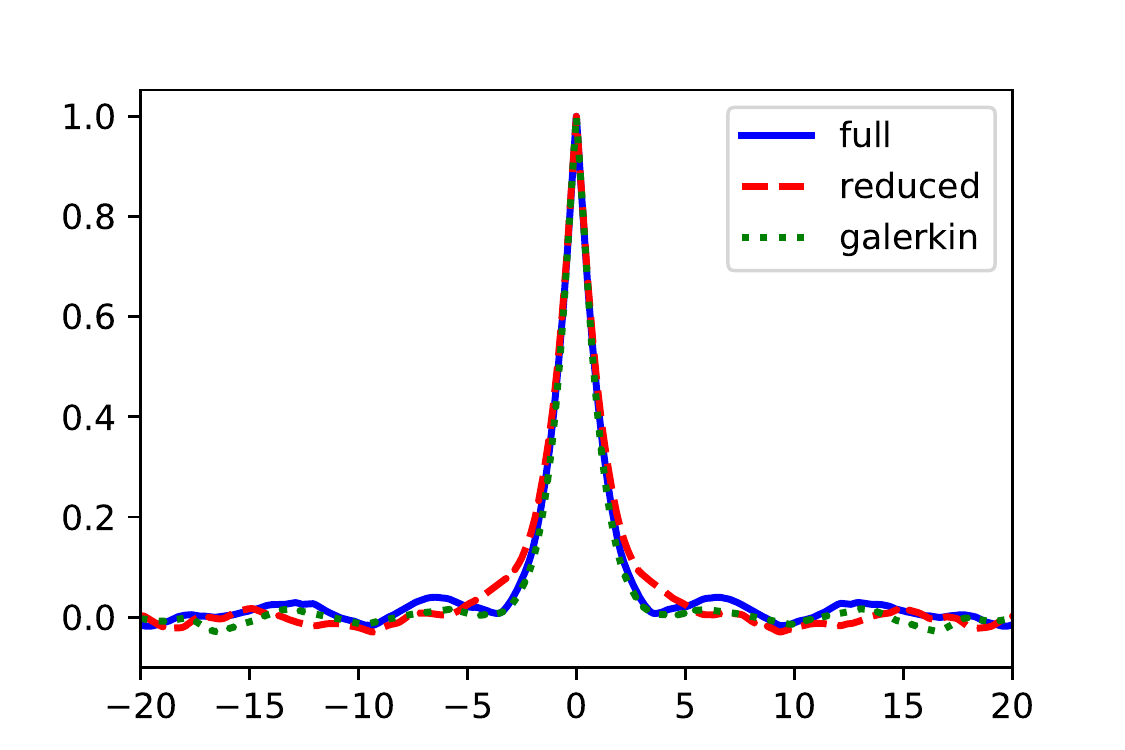}}&
      \rotatebox{90}{\footnotesize\hspace{0.5in}Mode $2$}&
      \resizebox{2.3in}{!}{\includegraphics{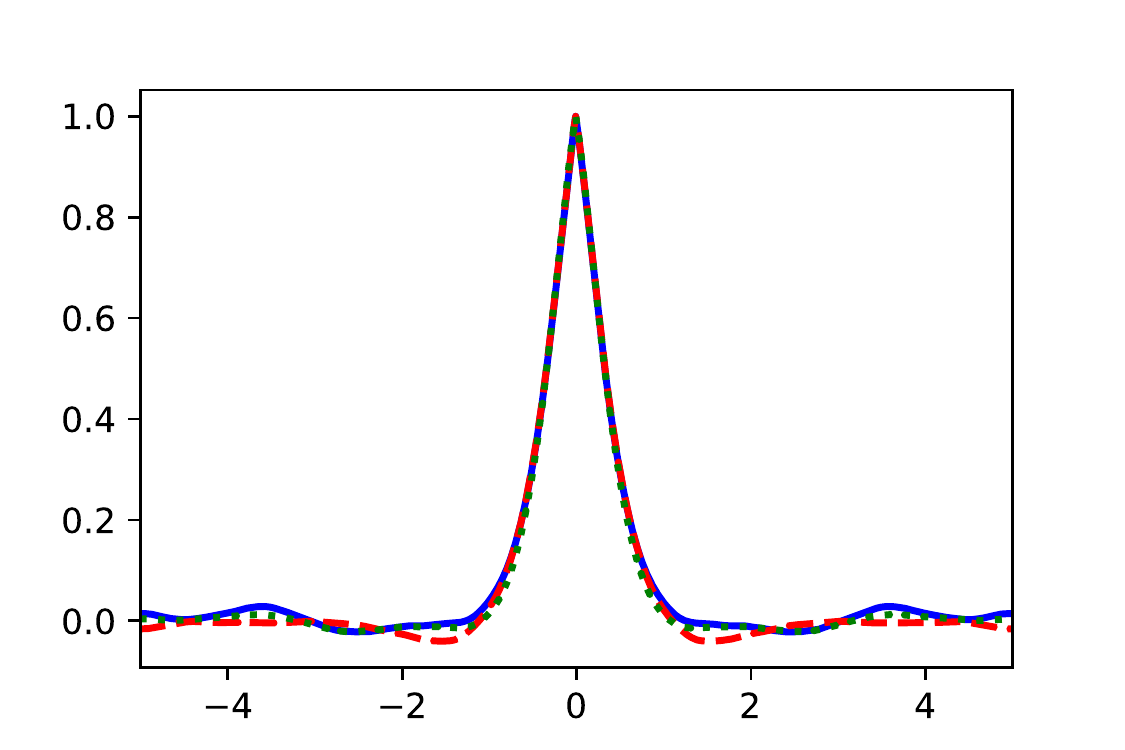}}&
      \rotatebox{90}{\footnotesize\hspace{0.5in}Mode $3$}&
      \resizebox{2.3in}{!}{\includegraphics{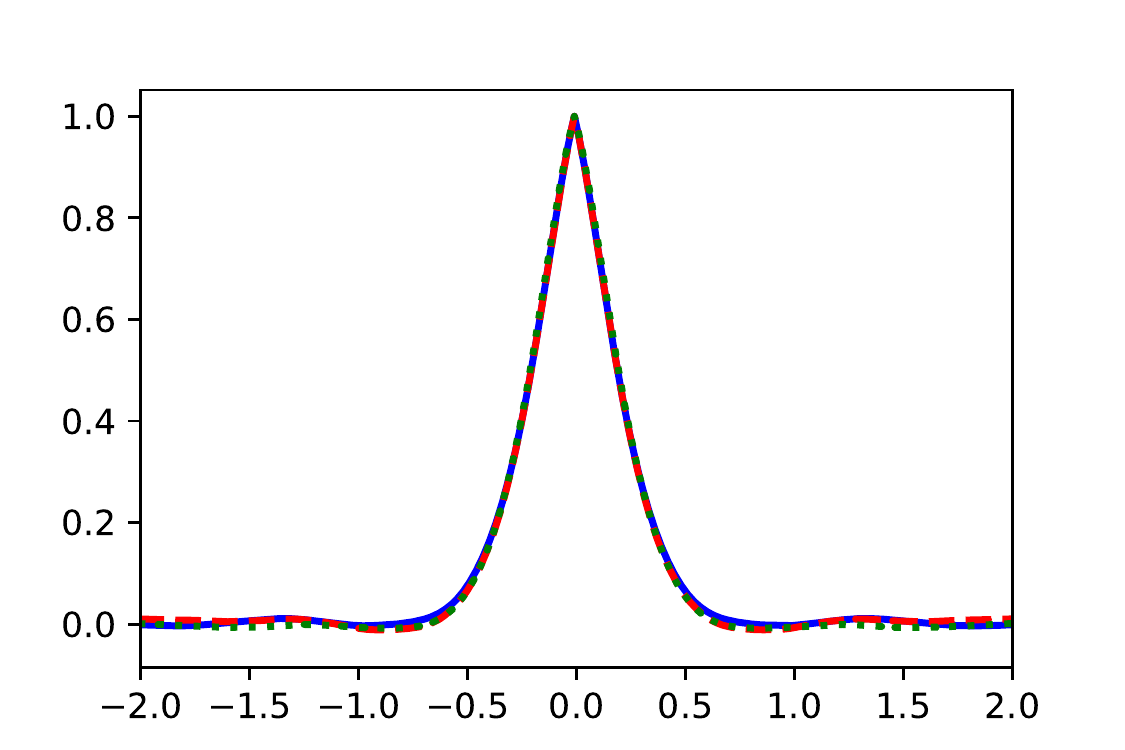}}\\[-6pt]
      \rotatebox{90}{\footnotesize\hspace{0.5in}Mode $4$}&
      \resizebox{2.3in}{!}{\includegraphics{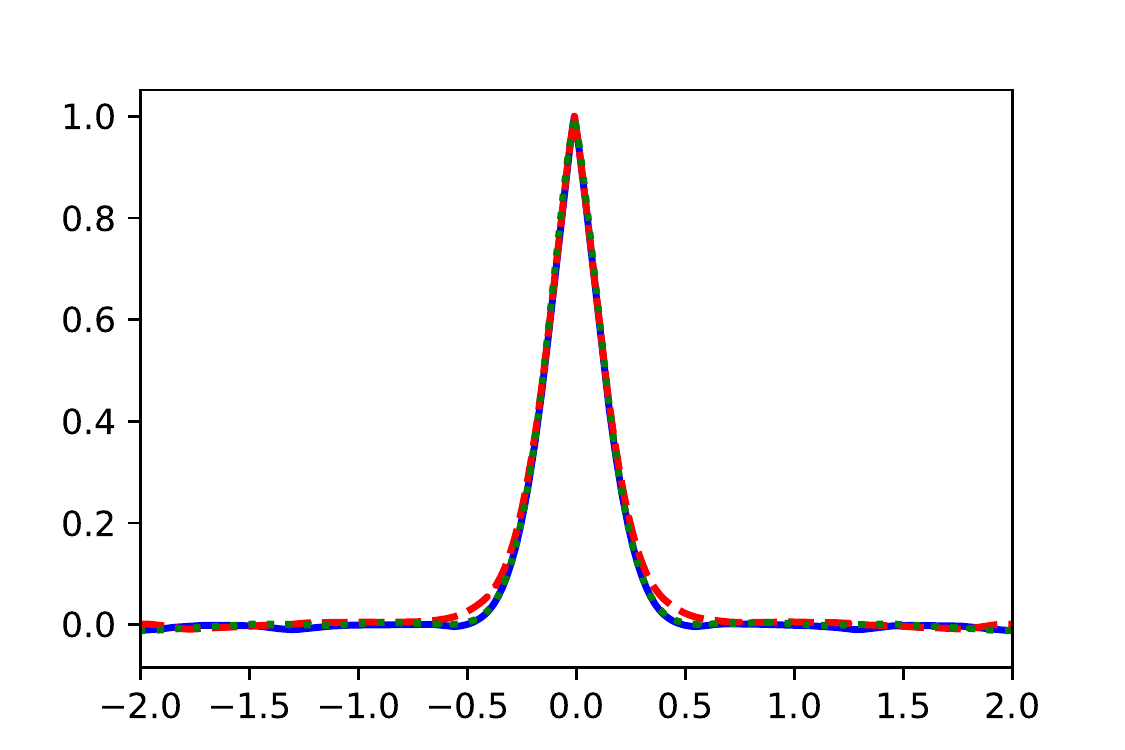}}&
      \rotatebox{90}{\footnotesize\hspace{0.5in}Mode $5$}&
      \resizebox{2.3in}{!}{\includegraphics{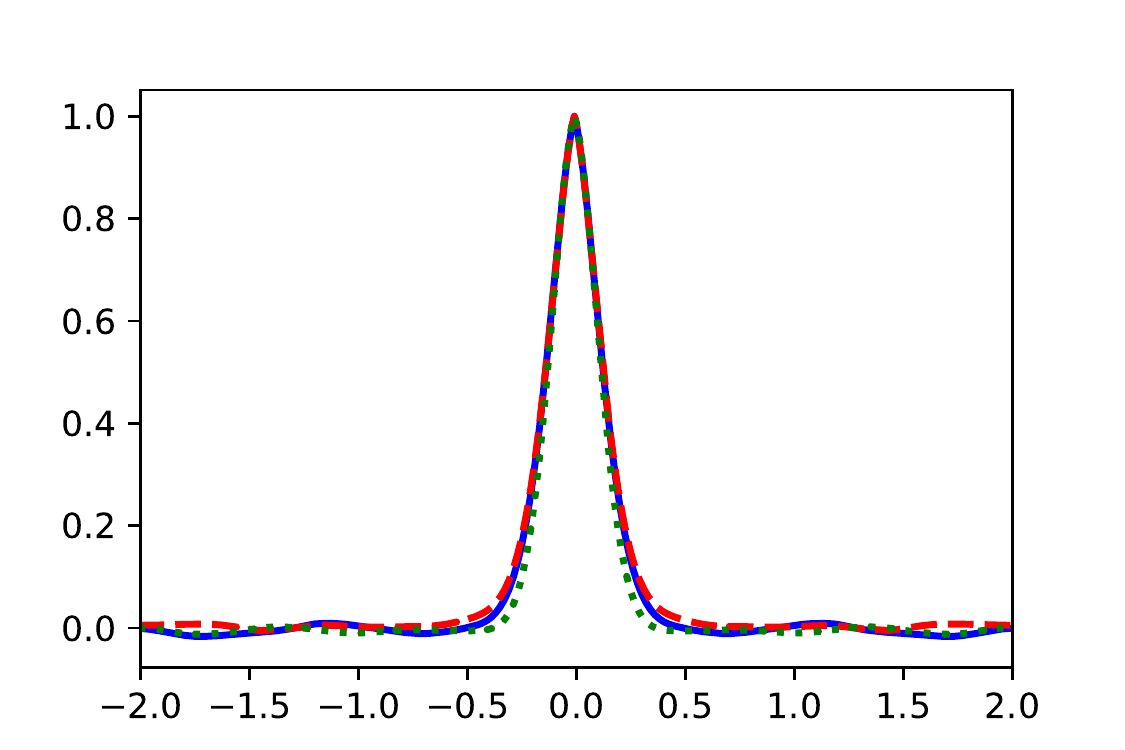}}&
      \rotatebox{90}{\footnotesize\hspace{0.5in}Mode $6$}&
      \resizebox{2.3in}{!}{\includegraphics{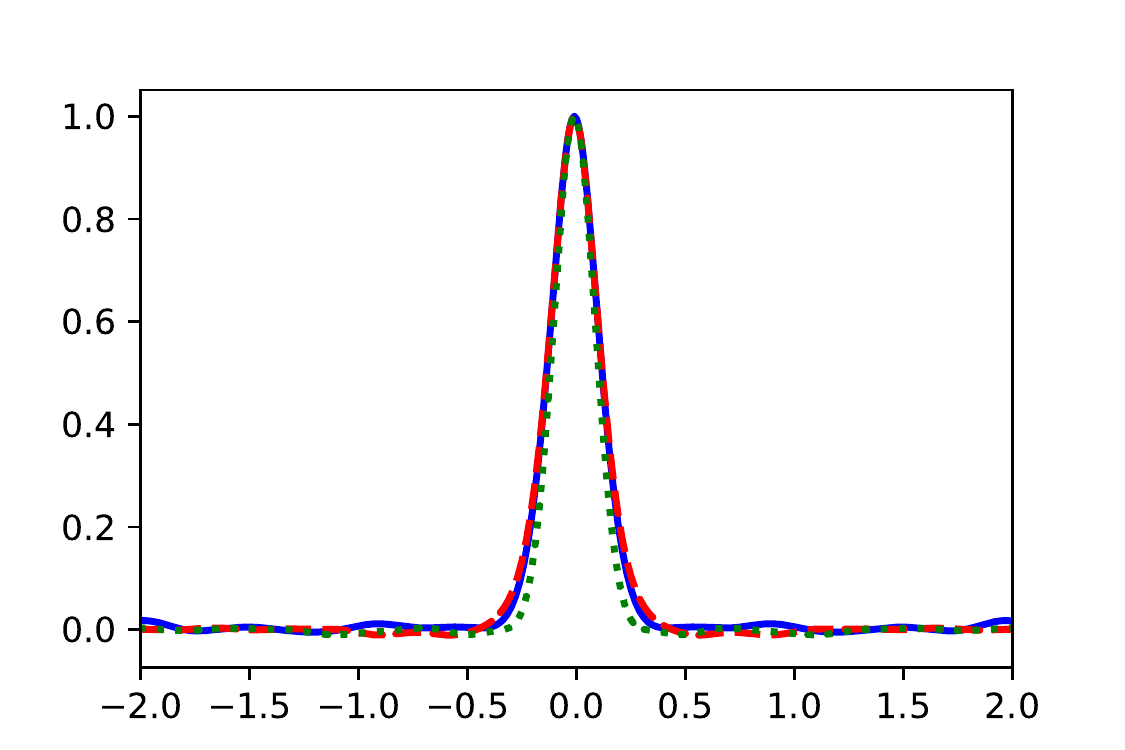}}\\[-6pt]
      \rotatebox{90}{\footnotesize\hspace{0.5in}Mode $7$}&
      \resizebox{2.3in}{!}{\includegraphics{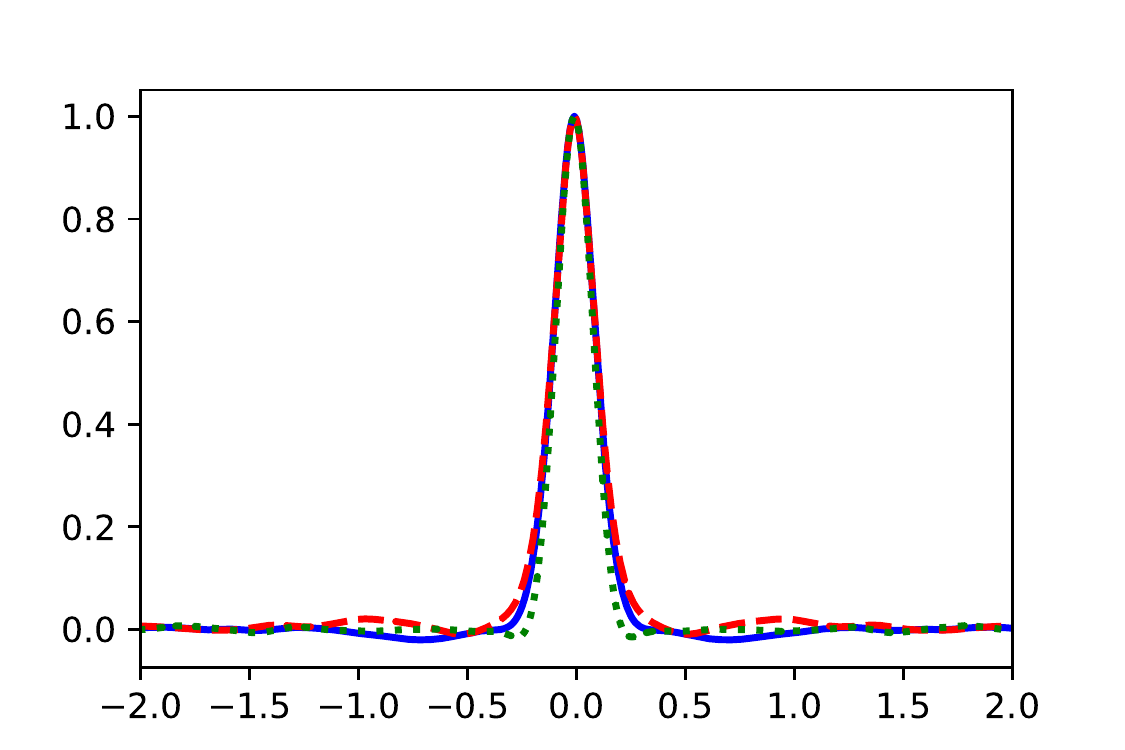}}&
      \rotatebox{90}{\footnotesize\hspace{0.5in}Mode $8$}&
      \resizebox{2.3in}{!}{\includegraphics{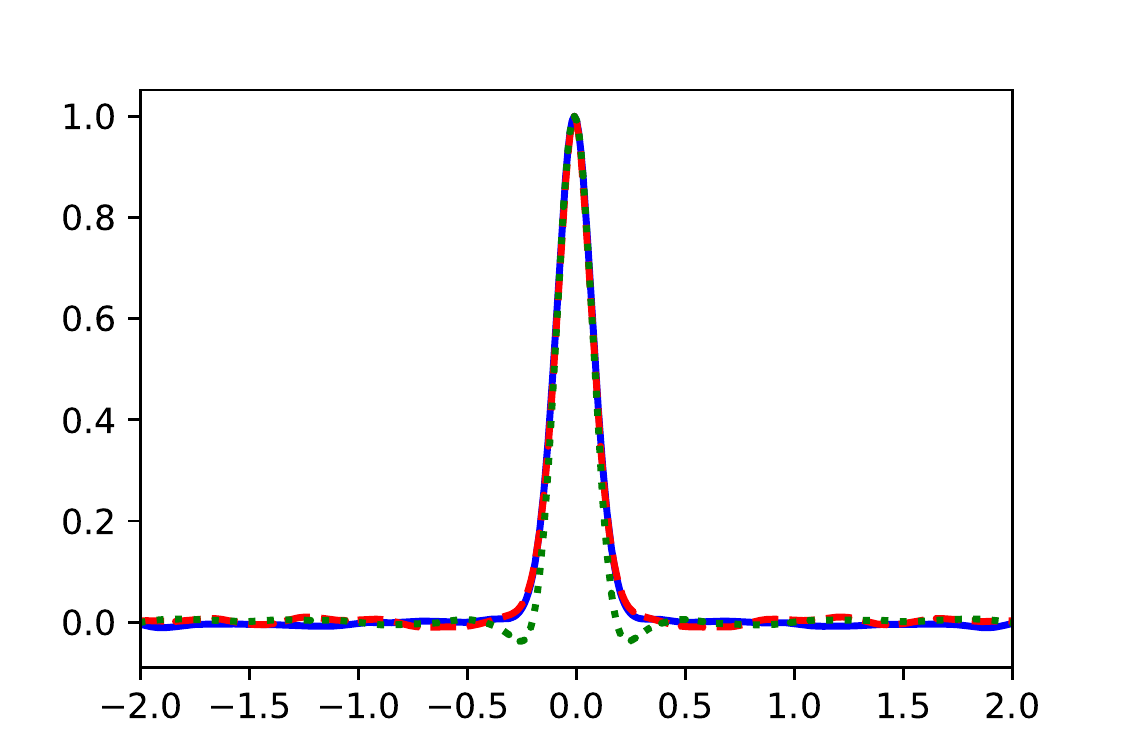}}&
      \rotatebox{90}{\footnotesize\hspace{0.5in}Mode $9$}&
      \resizebox{2.3in}{!}{\includegraphics{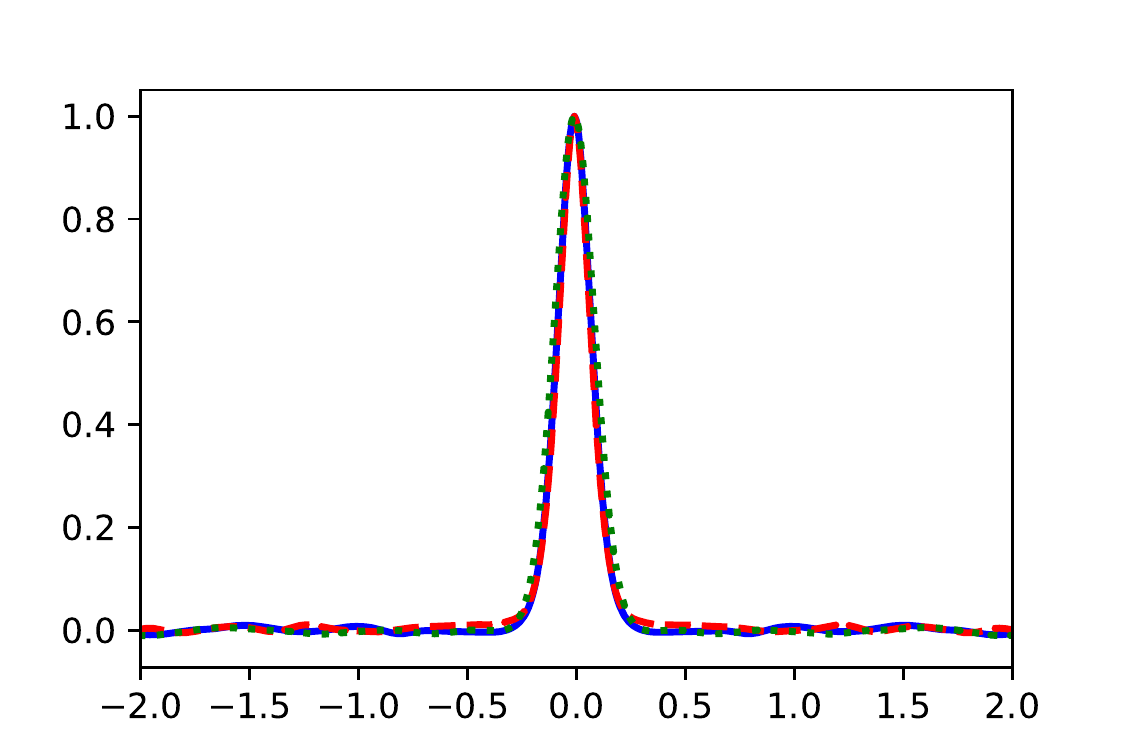}}\\
      &{\footnotesize Time lag} && {\footnotesize Time lag} &&
                {\footnotesize Time lag}\\
    \end{tabular}
  \end{center}
  \caption{Autocovariance functions for the stochastic Burgers
    equation.  We plot autocovariance functions for $Re(u_k)$ for
    $k=1,\cdots,9$.  In all panels, solid blue line is the full model
    (128-mode truncation), dashed red line is the 9-mode reduced model,
    and dotted green line the 9-mode Galerkin truncation.}
    
  \applabel{fig:burgers-acf}
\end{figure}

\begin{figure}
  \begin{center}
    \begin{tabular}{r@{\hskip 0pt}c@{\hskip -6pt}r@{\hskip 0pt}c@{\hskip 0pt}r@{\hskip 0pt}c}
      \rotatebox{90}{\footnotesize\hspace{0.5in}Mode $1$}&
      \resizebox{2.3in}{!}{\includegraphics{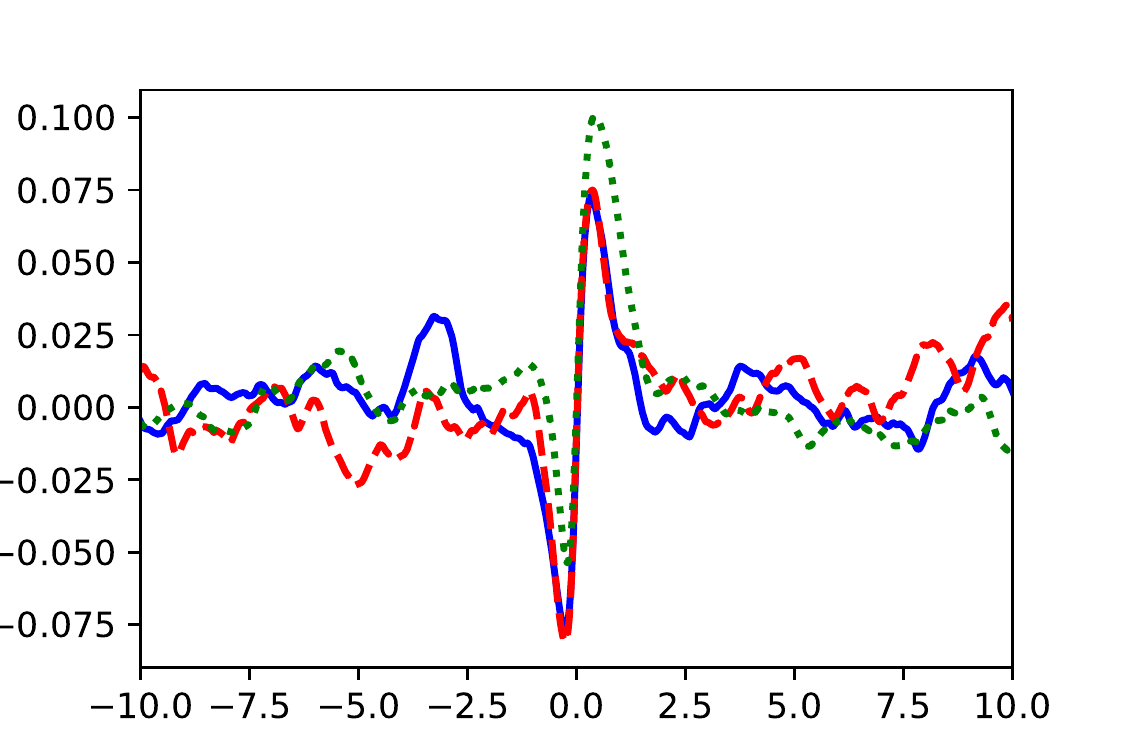}}&
      \rotatebox{90}{\footnotesize\hspace{0.5in}Mode $2$}&
      \resizebox{2.3in}{!}{\includegraphics{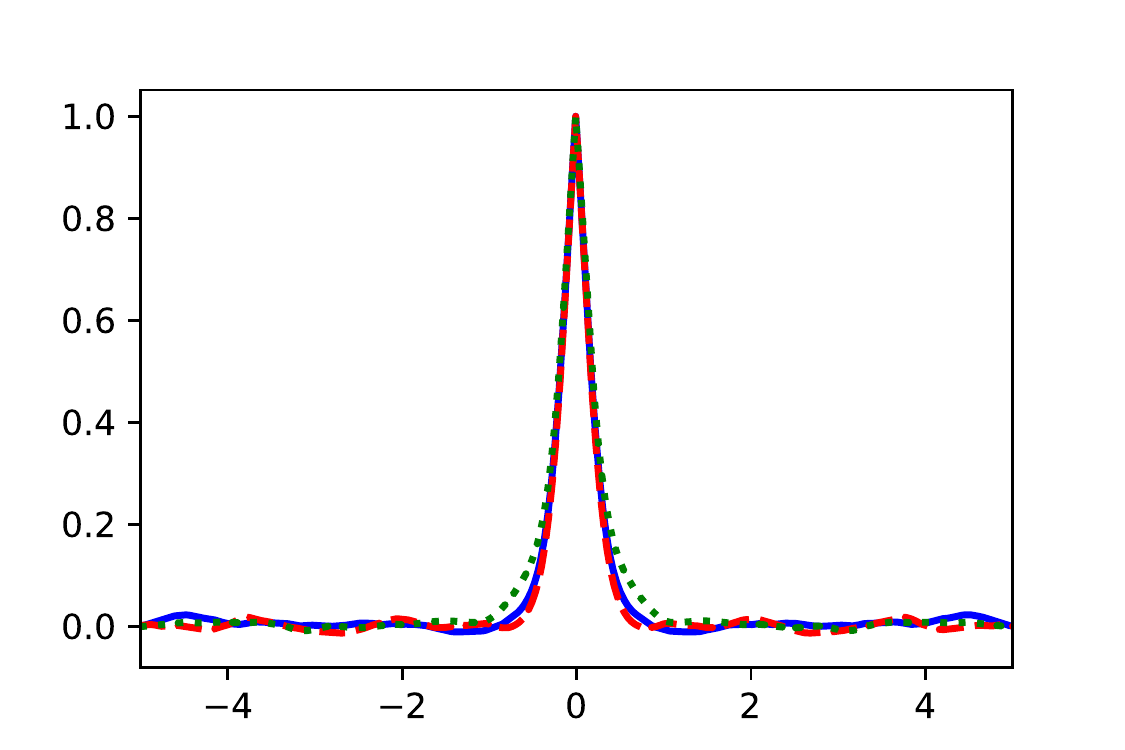}}&
      \rotatebox{90}{\footnotesize\hspace{0.5in}Mode $3$}&
      \resizebox{2.3in}{!}{\includegraphics{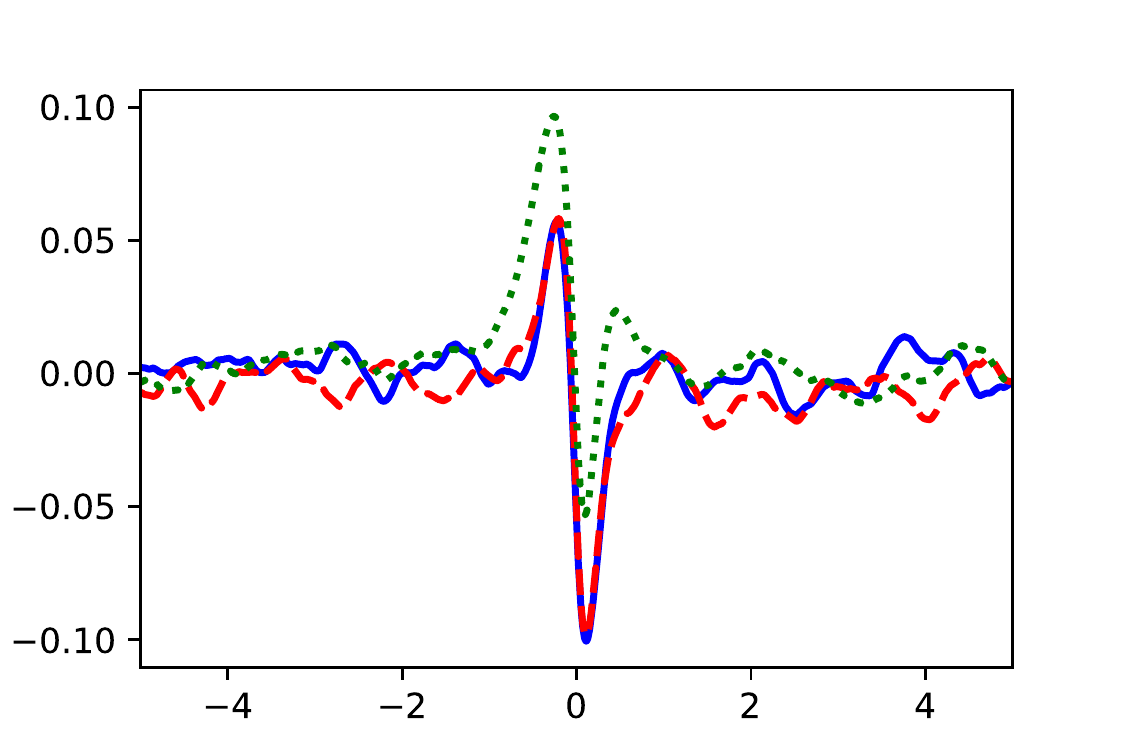}}\\[-6pt]
      \rotatebox{90}{\footnotesize\hspace{0.5in}Mode $4$}&
      \resizebox{2.3in}{!}{\includegraphics{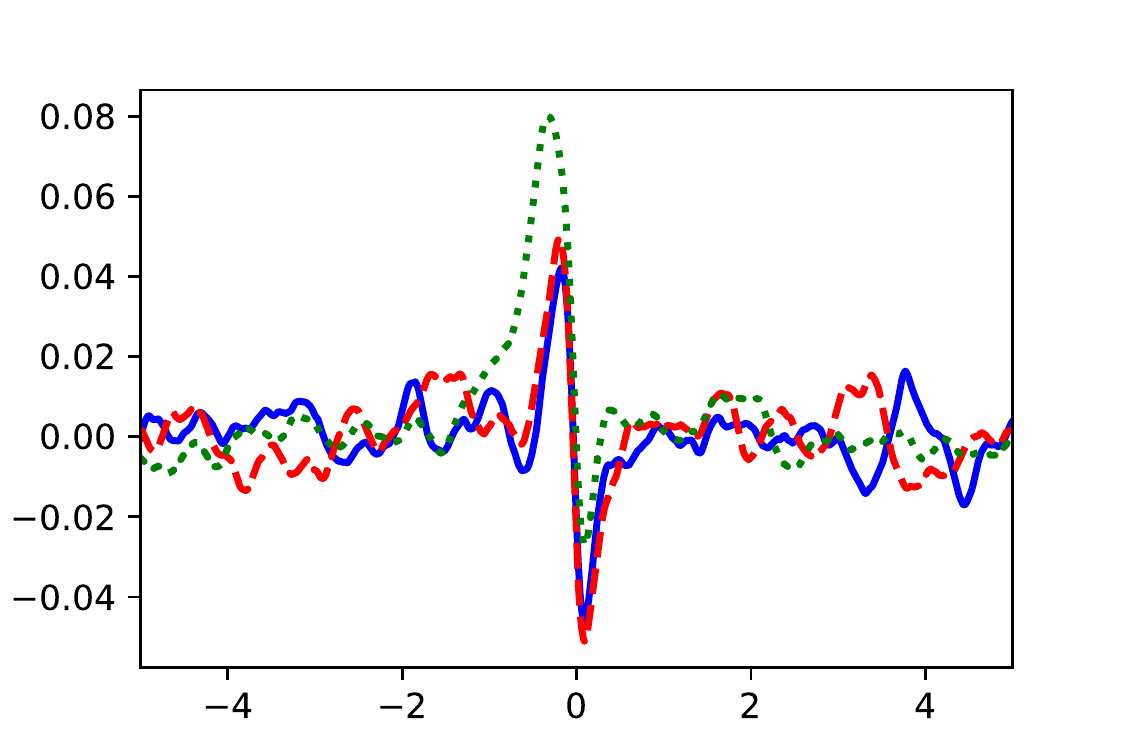}}&
      \rotatebox{90}{\footnotesize\hspace{0.5in}Mode $5$}&
      \resizebox{2.3in}{!}{\includegraphics{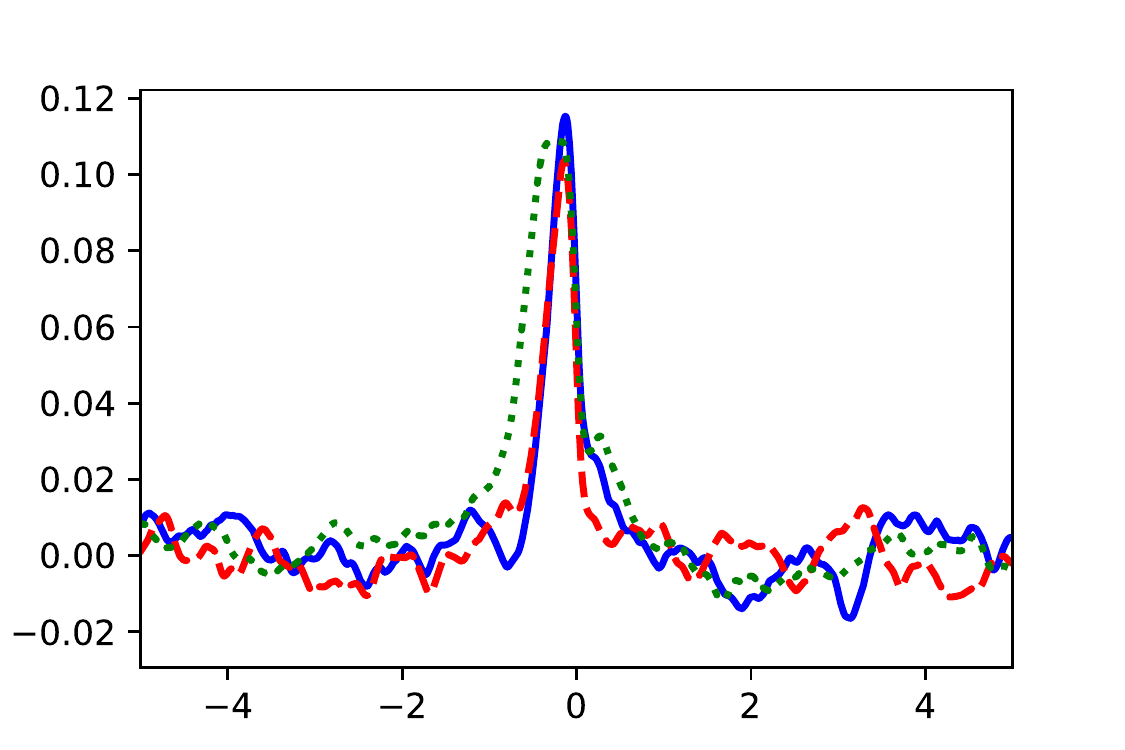}}&
      \rotatebox{90}{\footnotesize\hspace{0.5in}Mode $6$}&
      \resizebox{2.3in}{!}{\includegraphics{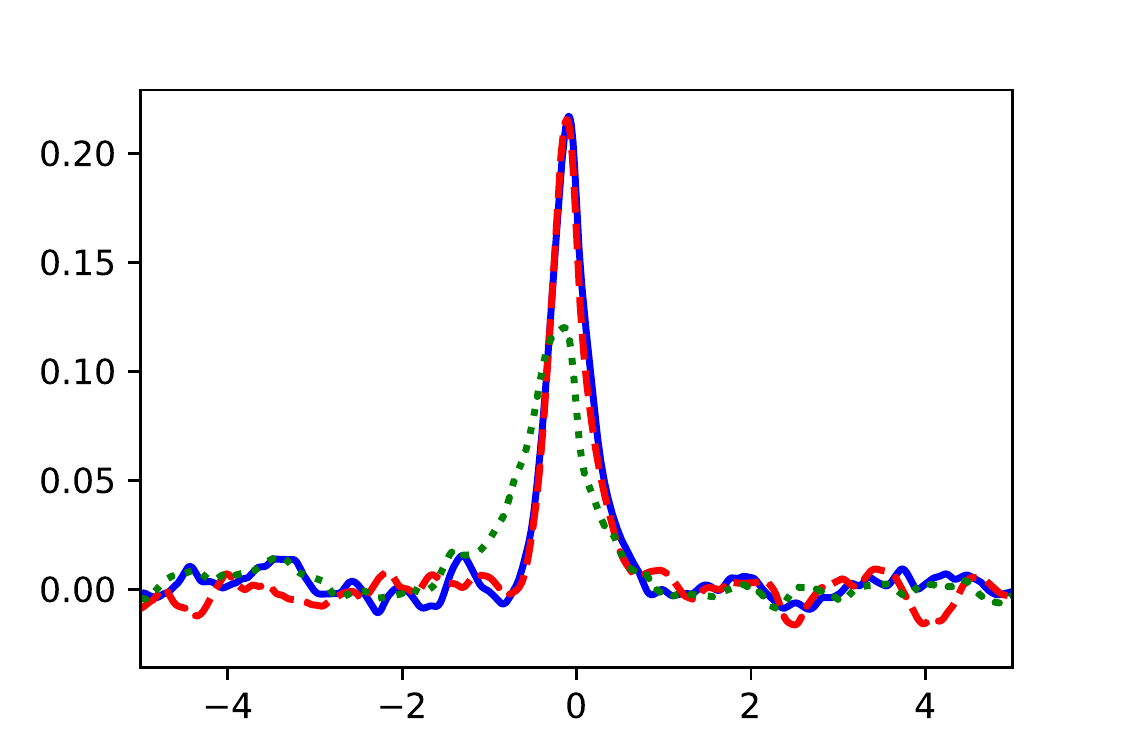}}\\[-6pt]
      \rotatebox{90}{\footnotesize\hspace{0.5in}Mode $7$}&
      \resizebox{2.3in}{!}{\includegraphics{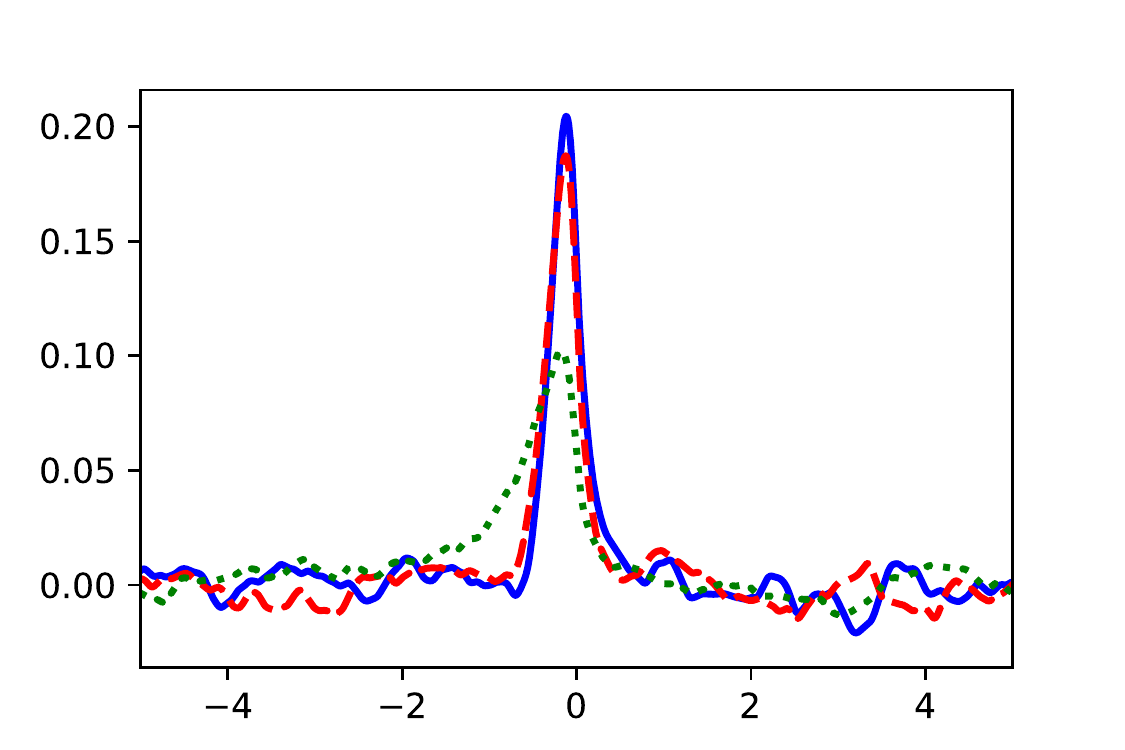}}&
      \rotatebox{90}{\footnotesize\hspace{0.5in}Mode $8$}&
      \resizebox{2.3in}{!}{\includegraphics{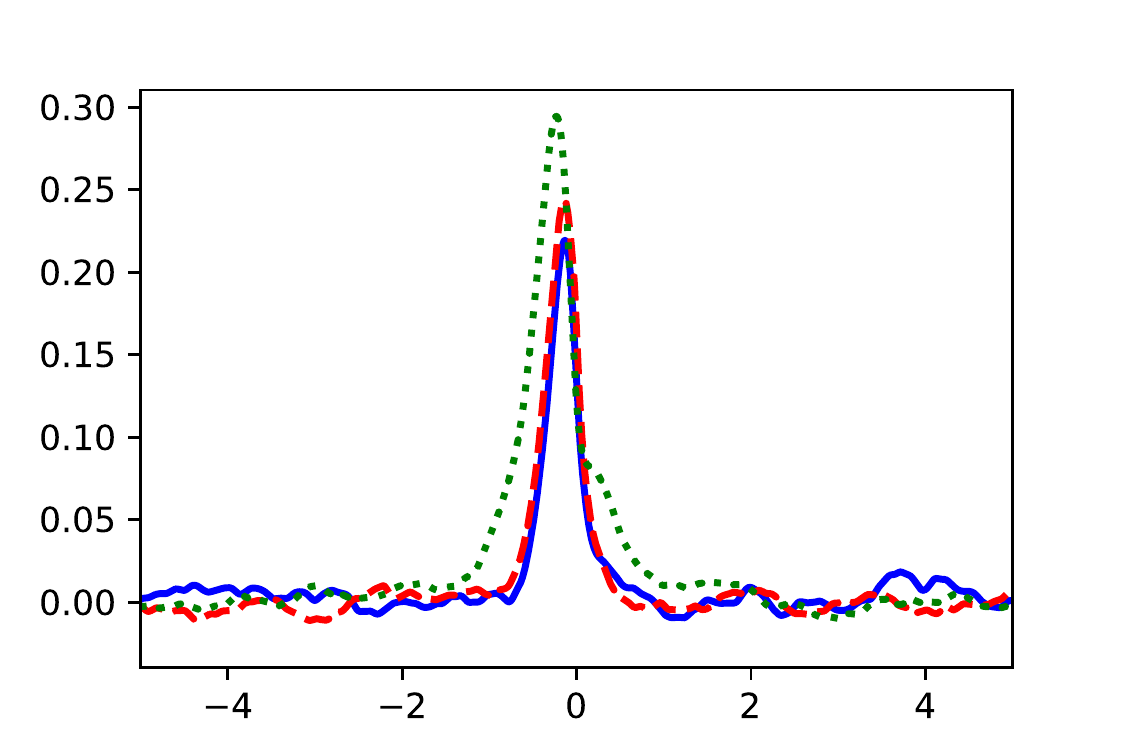}}&
      \rotatebox{90}{\footnotesize\hspace{0.5in}Mode $9$}&
      \resizebox{2.3in}{!}{\includegraphics{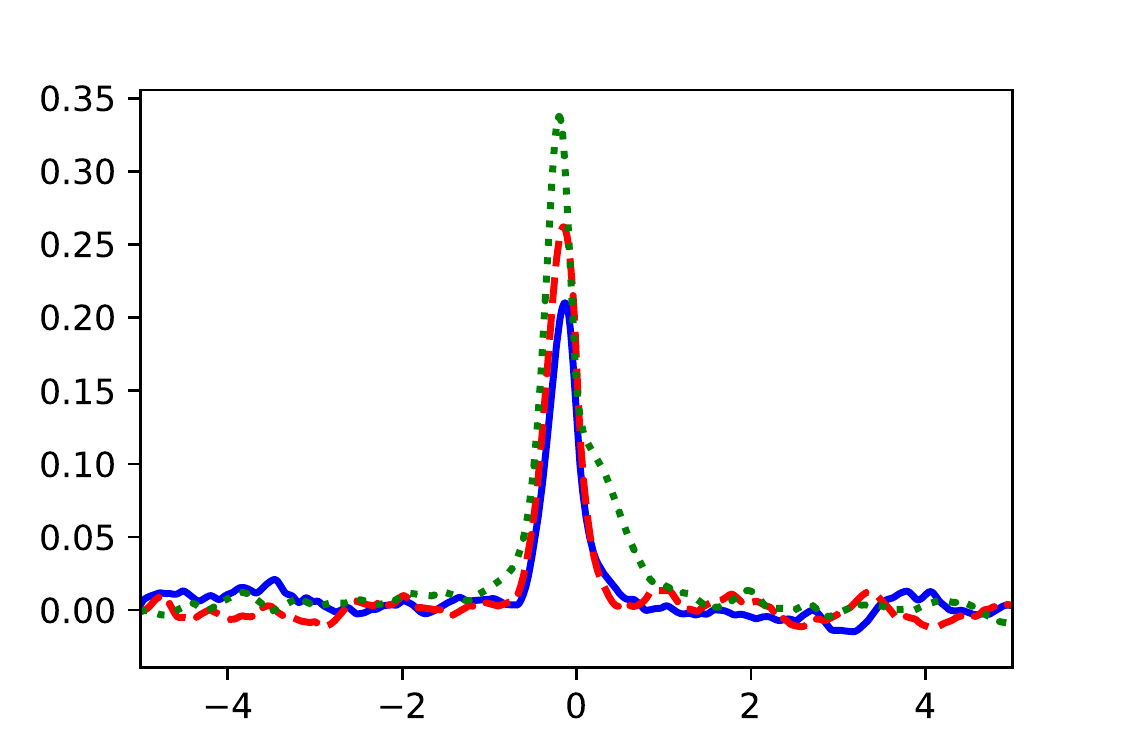}x2}\\
      &{\footnotesize Time lag} && {\footnotesize Time lag} &&
                {\footnotesize Time lag}\\
    \end{tabular}
  \end{center}
  \caption{Energy cross-correlation functions for the stochastic Burgers
    equation.  We plot cross correlation functions for $|u_2|^2$ and
    $|u_k|^2$ for $k=1,\cdots,9$.  In all panels, solid blue line is the
    full model (128-mode truncation), dashed red line is the 9-mode
    reduced model, and dotted green line the 9-mode Galerkin
    truncation.}
    
  \applabel{fig:burgers-ccf}
\end{figure}

\begin{figure}[tb!]
  \begin{center}
    \begin{tabular}{c@{\hskip -0.05in}c}
      \hspace{0.3in}Spacetime view, full model & Snapshots (legend below)\\
      \resizebox{!}{2.5in}{\includegraphics*[bb=0in 0.1in 1.8in 2.7in]{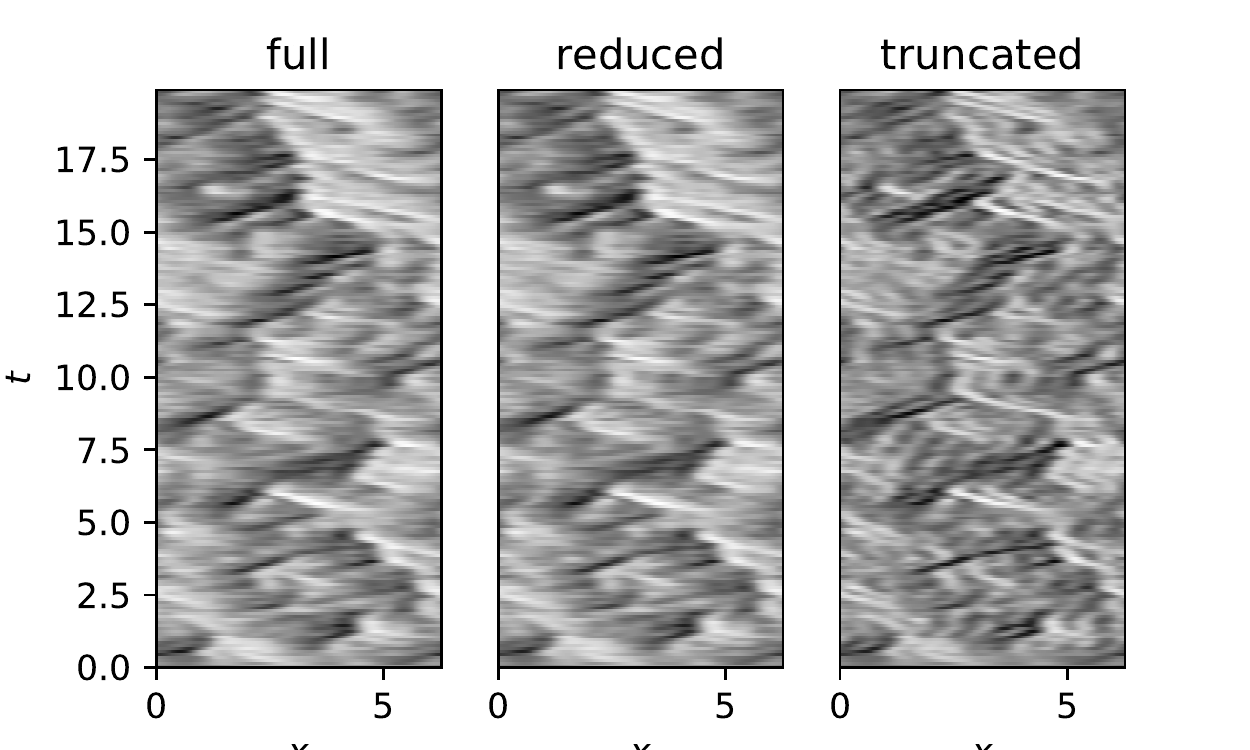}}&
      \resizebox{!}{2.5in}{\includegraphics*[bb=0.15in 0.2in 5in 3.7in]{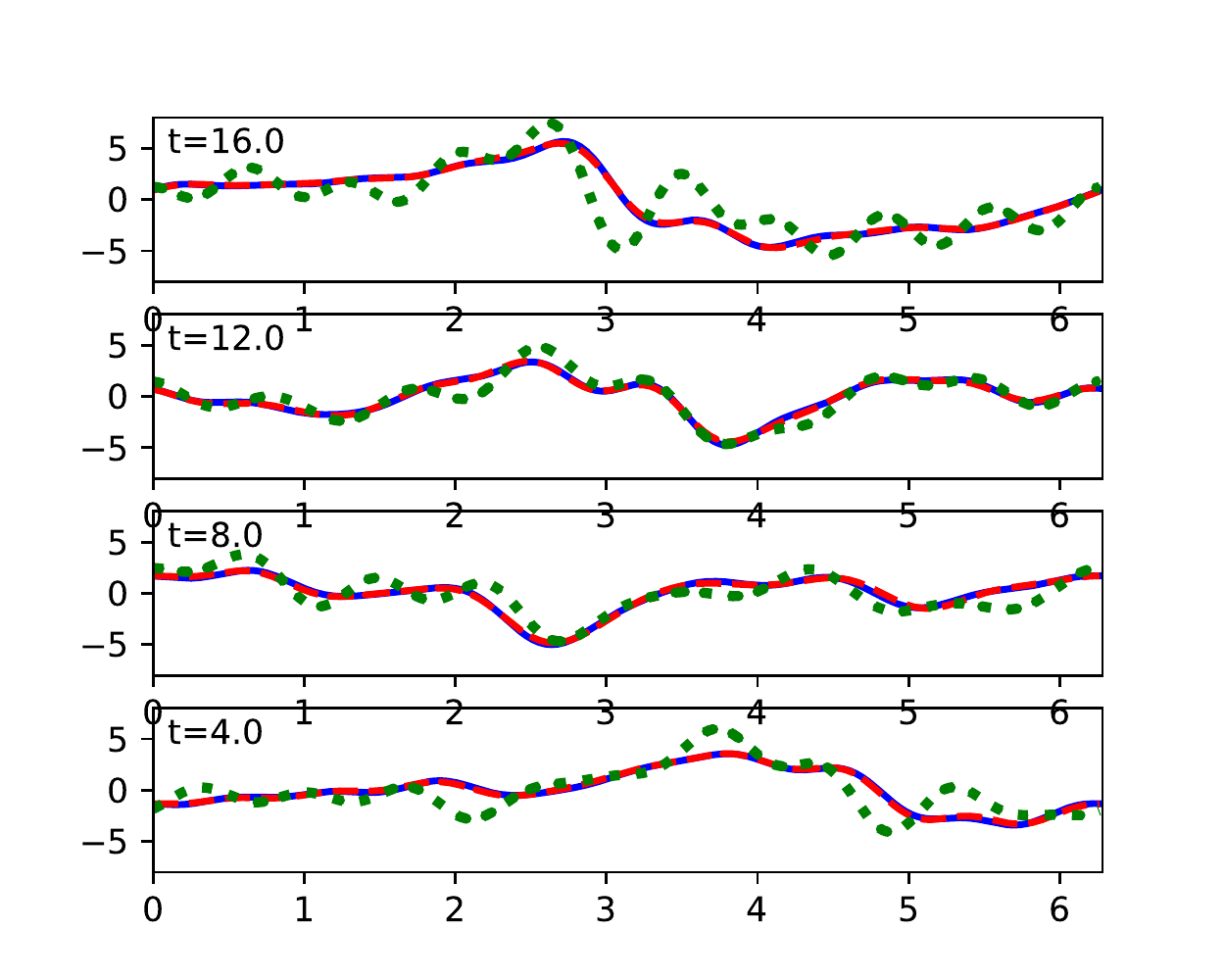}}\\[-2ex]
      \hspace{0.3in}$x$&$x$\\
    \end{tabular}\\
    (a) Burgers solutions\\[2ex]
    \begin{tabular}{c@{\hskip -0.2in}c}
      \begin{tabular}{r@{\hskip 0pt}c}
        \rotatebox{90}{\footnotesize\hspace{0.4in}$k=1$}&
        \resizebox{!}{1.3in}{\includegraphics{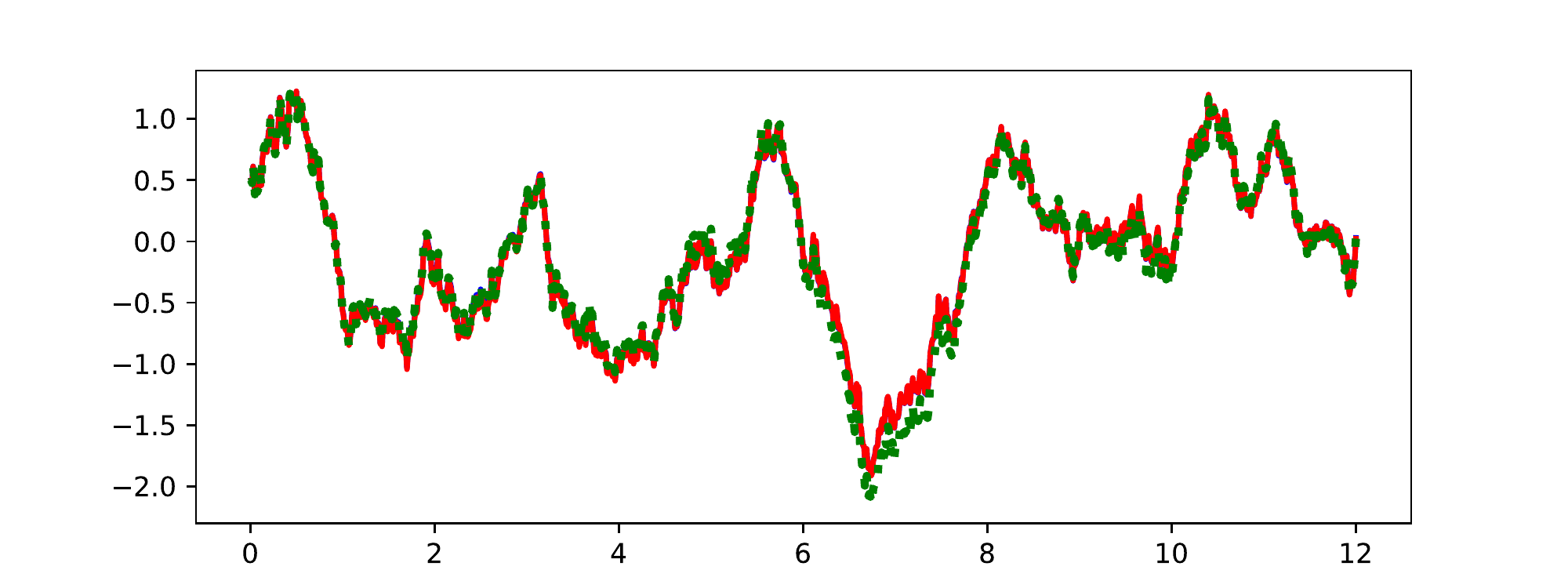}}\\[-1ex]
        \rotatebox{90}{\footnotesize\hspace{0.4in}$k=9$}&
        \resizebox{!}{1.3in}{\includegraphics{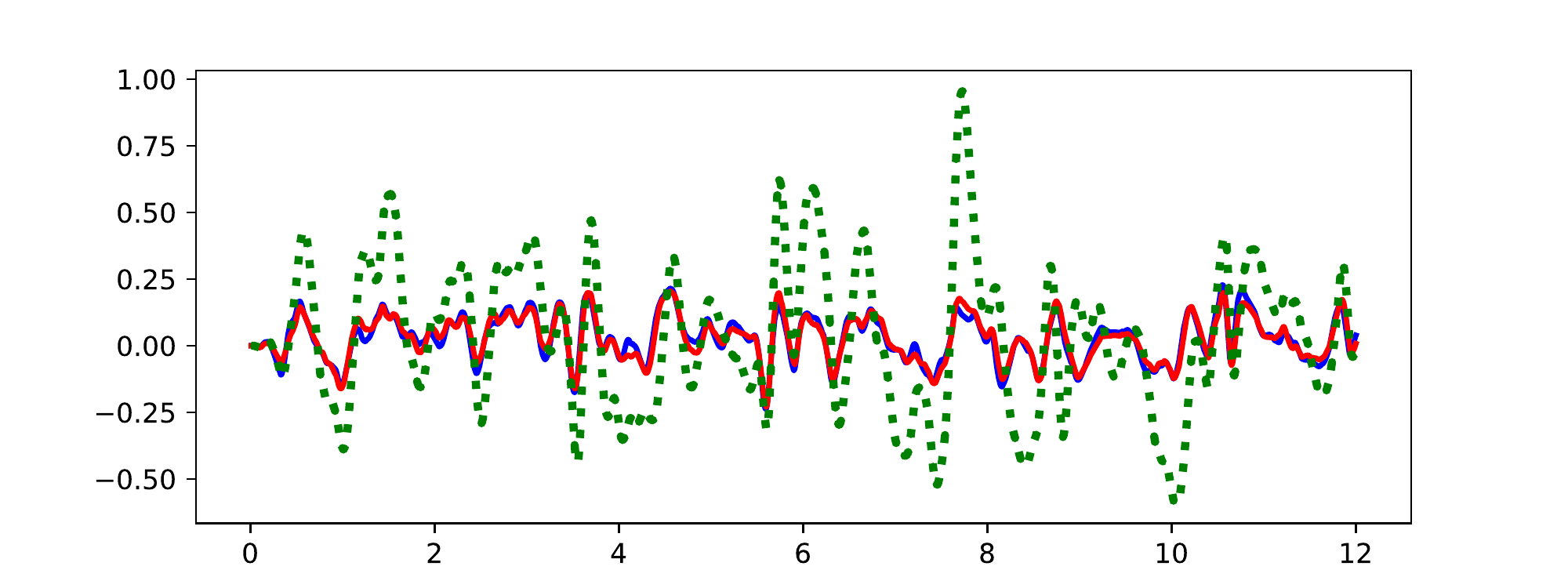}}\\[-1ex]
        &{\footnotesize $t$}\\
      \end{tabular} &

      \begin{tabular}{c}
        \resizebox{!}{0.6in}{\includegraphics*[bb=1.4in 1in 3in 2in]{./legend}}\\
        \resizebox{!}{2in}{\includegraphics{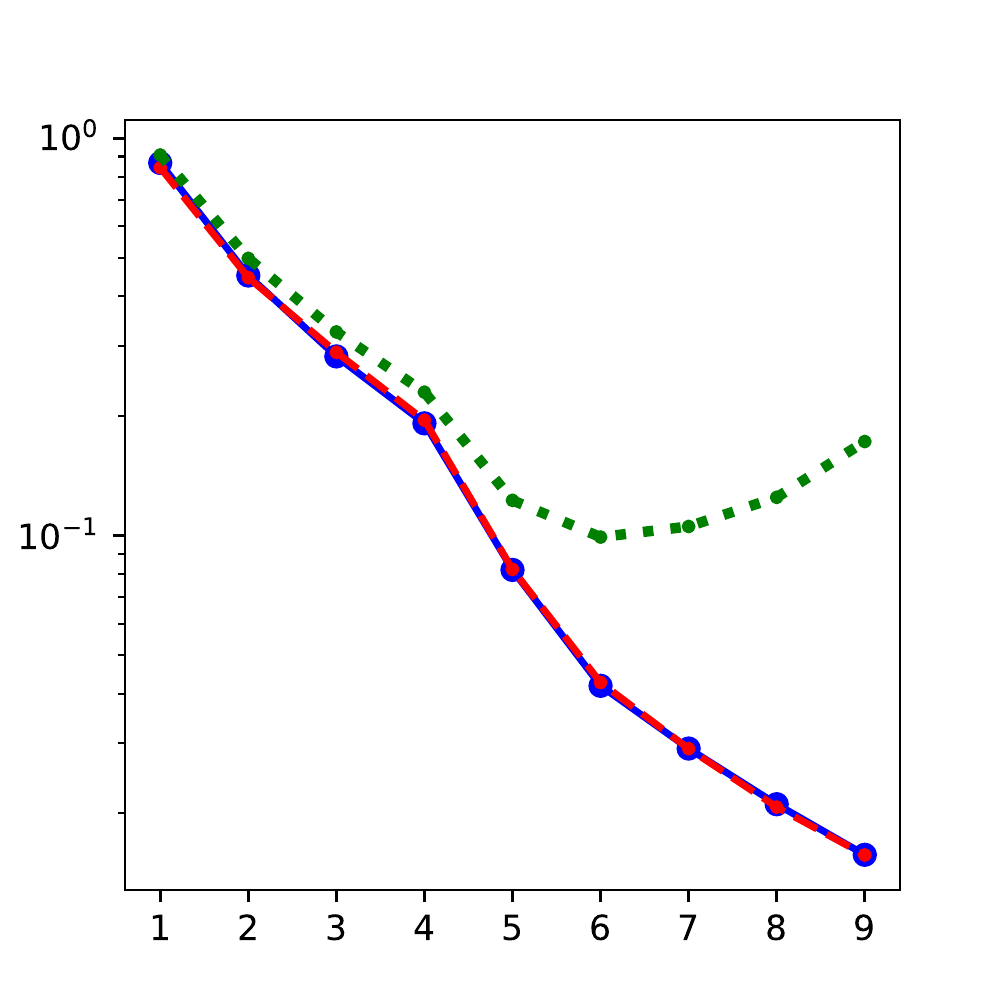}}\\[-2ex]
        {\footnotesize $k$}\\
      \end{tabular}\\
      (b) Trajectories of $Re(u_k(t))$ &  \hspace{6mm}
      (c) Energy $\braket{|u_k|^2}$\\
    \end{tabular}
  \end{center}
  \caption{The results using a linear regression with $p=1,r=1$ as in
    \keqref{eq:narma-one-step-prediction}. They are almost identical as
    those in Fig.~\ref{fig:burgers-solutions} from a nonlinear
    regression using the model in \keqref{eq:recursion} with $p=1,r=1$.}
  \label{fig:burgers-nonlin-comparison}
\end{figure}

\end{document}